\documentclass[aap,preprint]{imsart}
\RequirePackage[OT1]{fontenc}
\RequirePackage{amsthm,amsmath}
\RequirePackage[numbers]{natbib}
\RequirePackage[colorlinks,citecolor=blue,urlcolor=blue]{hyperref}

\allowdisplaybreaks

\arxiv{arXiv:0000.0000}

\startlocaldefs
\numberwithin{equation}{section}
\theoremstyle{plain}

\endlocaldefs

\usepackage[T1]{fontenc}
\usepackage[latin9]{inputenc}
\usepackage{amsmath, amssymb, amsthm,graphicx}
\usepackage{color,bm}
\usepackage{tikz,pgfplots}
\usepackage{bbold}
\usetikzlibrary{decorations.pathreplacing,calc}

\def\trans{ {\sf T} }
\def\tr{ {\rm tr\,} }
\def\RR{\mathbb{R}}
\def\CC{\mathbb{C}}

\newcommand{\asto}{\overset{\rm a.s.}{\longrightarrow}}
\pgfplotsset{compat=1.9}
\newtheorem{remark}{Remark}
\newtheorem{lemma}{Lemma}
\newtheorem{theorem}{Theorem}
\newtheorem{assumption}{Assumption}
\newtheorem{corollary}{Corollary}

\begin{document}

\begin{frontmatter}

\title{Random matrix-improved estimation\\ of covariance matrix distances}
\runtitle{Improved covariance matrix distance estimation}

\begin{aug}

	\author{\fnms{Romain} \snm{Couillet}\thanksref{m1,m2,t1}\ead[label=e1]{romain.couillet@gipsa-lab.grenoble-inp.fr}}
	\author{\fnms{Malik} \snm{Tiomoko}\thanksref{m1}\ead[label=e2]{malik.tiomoko@lss.centralesupelec.fr}},
	\author{\fnms{Steeve} \snm{Zozor}\thanksref{m2}\ead[label=e3]{steeve.zozor@gipsa-lab.grenoble-inp.fr}}
	\and
	\author{\fnms{Eric} \snm{Moisan}\thanksref{m2}\ead[label=e4]{eric.moisan@gipsa-lab.grenoble-inp.fr}}

\thankstext{m1}{L2S, CentraleSup\'elec, University of Paris Saclay, France.}
\thankstext{m2}{GIPSA-lab, University Grenoble Alpes, France.}
\thankstext{t1}{Couillet's work is supported by the ANR Project RMT4GRAPH (ANR-14-CE28-0006) and by the IDEX GSTATS Chair at University Grenoble Alpes.}

\runauthor{R.\@ Couillet et al.}


%

\end{aug}

\begin{abstract}
	Given two sets $x_1^{(1)},\ldots,x_{n_1}^{(1)}$ and $x_1^{(2)},\ldots,x_{n_2}^{(2)}\in\RR^p$ (or $\CC^p$) of random vectors with zero mean and positive definite covariance matrices $C_1$ and $C_2\in\RR^{p\times p}$ (or $\CC^{p\times p}$), respectively, this article provides novel estimators for a wide range of distances between $C_1$ and $C_2$ (along with divergences between some zero mean and covariance $C_1$ or $C_2$ probability measures) of the form $\frac1p\sum_{i=1}^n f(\lambda_i(C_1^{-1}C_2))$ (with $\lambda_i(X)$ the eigenvalues of matrix $X$). These estimators are derived using recent advances in the field of random matrix theory and are asymptotically consistent as $n_1,n_2,p\to\infty$ with non trivial ratios $p/n_1<1$ and $p/n_2<1$ (the case $p/n_2>1$ is also discussed). A first ``generic'' estimator, valid for a large set of $f$ functions, is provided under the form of a complex integral. Then, for a selected set of $f$'s of practical interest (namely, $f(t)=t$, $f(t)=\log(t)$, $f(t)=\log(1+st)$ and $f(t)=\log^2(t)$), a closed-form expression is provided. Beside theoretical findings, simulation results suggest an outstanding performance advantage for the proposed estimators when compared to the classical ``plug-in'' estimator $\frac1p\sum_{i=1}^n f(\lambda_i(\hat C_1^{-1}\hat C_2))$ (with $\hat C_a=\frac1{n_a}\sum_{i=1}^{n_a}x_i^{(a)}x_i^{(a)\trans}$), and this even for very small values of $n_1,n_2,p$.
\end{abstract}



\end{frontmatter}

\section{Introduction}
\label{sec:intro}

In a host of statistical signal processing and machine learning methods, distances between covariance matrices are regularly sought for. These are notably exploited to estimate centroids or distances between clusters of data vectors mostly distinguished through their second order statistics. We may non exhaustively cite brain graph signal processing and machine learning (from EEG datasets in particular) which is a field largely rooted in these approaches \citep{FAL14}, hyperspectral and synthetic aperture radar (SAR) clustering \citep{CHA03,VAS10}, patch-based image processing \citep{HOU17}, etc. For random independent $p$-dimensional real or complex data vectors $x_1^{(a)},\ldots,x_{n_a}^{(a)}$, $a\in\{1,2\}$, having zero mean and covariance matrix $C_a$, and for a distance (or divergence) $D(X,Y)$ between covariance matrices $X$ and $Y$ (or probability measures associated to random variables with these covariances), the natural approach is to estimate $D(C_1,C_2)$ through the ``plug-in'' substitute $D(\hat C_1,\hat C_2)$. For well-behaved functions $D$, this generally happens to be a consistent estimator as $n_1,n_2\to\infty$ in the sense that $D(\hat C_1,\hat C_2)\to D(C_1,C_2)$ almost surely as $n_1,n_2\to\infty$ while $p$ remains fixed. This is particularly the case for all subsequently introduced distances and divergences.

\bigskip

However, in many modern applications, one cannot afford large $n_a$ ($a\in\{1,2\}$) values or, conversely, $p$ may be commensurable, if not much larger, than the $n_a$'s. When $n_a,p\to\infty$ in such a way that $p/n_a$ remains away from zero and infinity, it has been well documented in the random matrix literature, starting with the seminal works of Mar\u{c}enko and Pastur~\citep{MAR67}, that the operator norm $\|\hat C_a-C_a\|$ no longer vanishes. This entails, as a consequence, that the aforementioned estimator $D(\hat C_1,\hat C_2)$ for $D(C_1,C_2)$ is likely to be inconsistent as $p,n_1,n_2\to\infty$ at a commensurable rate.

\bigskip

This said, it is now interesting to note that many standard matrix distances $D(C_1,C_2)$ classically used in the literature can be written under the form of \emph{functionals} of the eigenvalues of $C_1^{-1}C_2$ (assuming at least $C_1$ is invertible). A first important example is the squared Fisher distance between $C_1$ and $C_2$ \cite{COS15} given by
\begin{align*}
	D_{\rm F}(C_1,C_2)^2 &= \frac1p\left\|\log(C_1^{-\frac12}C_2C_1^{-\frac12})\right\|_F^2 = \frac1p\sum_{i=1}^p \log^2(\lambda_i(C_1^{-1}C_2))
\end{align*}
with $\|\cdot\|_F$ the matrix Frobenius norm, $X^{\frac12}$ the unique nonnegative definite square root of $X$, and $\log(X)\equiv U\log(\Lambda)U^\trans$ for symmetric $X=U\Lambda U^\trans$ (and $(\cdot)^\trans$ the matrix transpose operator) in its spectral composition. This estimator, arising from information geometry, corresponds to the length of the geodesic between $C_1$ and $C_2$ in the manifold of positive definite matrices.

Another example is the Bhattacharyya distance \cite{BHA43} between two real Gaussian distributions with zero mean and covariances $C_1$ and $C_2$, respectively (i.e., $\mathcal N(0,C_1)$ and $\mathcal N(0,C_2)$), which measures some similarity between the two laws and reads
\begin{align*}
	D_{\rm B}(C_1,C_2) &= \frac1{2p}\log\det\left( \frac12[C_1+C_2]\right)-\frac1{4p}\log\det C_1-\frac1{4p}\log\det C_2
\end{align*}
which can be rewritten under the form
\begin{align*}
	D_{\rm B}(C_1,C_2) &= \frac1{2p}\log\det(I_p+C_1^{-1}C_2) - \frac1{4p}\log\det (C_1^{-1}C_2) - \frac12\log(2) \\
	&= \frac1{2p}\sum_{i=1}^p \log(1+\lambda_i(C_1^{-1}C_2)) - \frac1{4p} \sum_{i=1}^p \log(\lambda_i(C_1^{-1}C_2)) - \frac12\log(2).
\end{align*}

In a similar manner, the Kullback-Leibler divergence \cite{BAS13b} of the Gaussian distribution $\mathcal N(0,C_2)$ with respect to $\mathcal N(0,C_1)$ is given by
\begin{align*}
	D_{\rm KL} &= \frac1{2p}\tr(C_1^{-1}C_2) -\frac12 + \frac1{2p}\log\det (C_1^{-1}C_2) \\
	&= \frac1{2p}\sum_{i=1}^p \lambda_i(C_1^{-1}C_2) - \frac12 + \frac1{2p}\sum_{i=1}^p\log (\lambda_i(C_1^{-1}C_2)).
\end{align*}

More generally, the R\'enyi divergence \cite{BAS13} of $\mathcal N(0,C_2)$ with respect to $\mathcal N(0,C_1)$ reads, for $\alpha\in\RR\setminus\{1\}$,
\begin{align}
	\label{eq:Renyi}
	D^\alpha_{\rm R} &= -\frac1{2(\alpha-1)} \frac1p\sum_{i=1}^p \log\left( \alpha + (1-\alpha) \lambda_i(C_1^{-1}C_2) \right) + \frac1{2p} \sum_{i=1}^p\log\left( \lambda_i(C_1^{-1}C_2) \right)
\end{align}
(one may check that $\lim_{\alpha\to 1}D^\alpha_{\rm R}=D_{\rm KL}$).

\bigskip

Revolving around recent advances in random matrix theory, this article provides a generic framework to consistently estimate such functionals of the eigenvalues of $C_1^{-1}C_2$ from the samples in the regime where $p,n_1,n_2$ are simultaneously large, under rather mild assumptions. In addition to the wide range of potential applications, as hinted at above, this novel estimator provides in practice a dramatic improvement over the conventional covariance matrix ``plug-in'' approach, as we subsequently demonstrate on a synthetic (but typical) example in Table~\ref{tab:1}. Here, for $C_1^{-\frac12}C_2C_1^{-\frac12}$ (having the same eigenvalues as $C_1^{-1}C_2$) a Toeplitz positive definite matrix, we estimate the (squared) Fisher distance $D_{\rm F}(C_1,C_2)^2$ (averaged over a large number of realizations of zero mean Gaussian $x_i^{(a)}$'s), for $n_1=1\,024$, $n_2=2\,048$ and $p$ varying from $2$ to $512$. A surprising outcome, despite the theoretical request that $p$ must be large for our estimator to be consistent, is that, already for $p=2$ (while $n_1,n_2\sim 10^3$), our proposed estimator largely outperforms the classical approach; for $n_a/p\sim 10$ or less, the distinction in performance between both methods is dramatic with the classical estimator biased by more than $100\%$. 

Our main result, Theorem~\ref{th:contour_integral}, provides a consistent estimator for functionals $f$ of the eigenvalues of $C_1^{-1}C_2$ under the form of a complex integral, valid for all functions $f:\RR\to \RR$ that have natural complex analytic extensions on given bounded regions of $\CC$. This estimator however assumes a complex integral form which we subsequently express explicitly for a family of functions $f$ in a series of corollaries (Corollary~\ref{cor:t} to Corollary~\ref{cor:log2t}). While Corollaries~\ref{cor:t}--\ref{cor:log1st} provide an exact calculus of the form provided in Theorem~\ref{th:contour_integral}, for the case $f(t)=\log^2(t)$, covered in Corollary~\ref{cor:log2t}, the exact calculus leads to an expression involving dilogarithm functions which are not elementary functions. For this reason, Corollary~\ref{cor:log2t} provides a large $p$ approximation of Theorem~\ref{th:contour_integral}, thereby leading to another (equally valid) consistent estimator. This explains why Table~\ref{tab:1} and the figures to come (Figure~\ref{fig:error} and Figure~\ref{fig:error_constant_c's}) display two different sets of estimates.
In passing, it is worth noticing that in the complex Gaussian case, the estimators offer a noticeable improvement on average over their real counterparts; this fact is likely due to a bias in the second-order fluctuations of the estimators, as discussed in Section~\ref{sec:conclusion}.

\begin{table}[h!]
	\centering
	\hspace*{-3cm}\begin{tabular}{l|rrrrrrrrr}
		$p$ & 2 & 4 & 8 & 16 & 32 & 64 & 128 & 256 & 512 \\
		\hline\\
		$D_{\rm F}(C_1,C_2)^2$ & 0.0980 &  0.1456  & 0.1694  &  0.1812 &   0.1872  &  0.1901 &   0.1916 &   0.1924 &   0.1927 \\
		& \\
		\hline\\
		Proposed estimator & 0.0993  &  0.1470  &  0.1708  &  0.1827  &  0.1887  &  0.1918  &  0.1933  &  0.1941  &  0.1953 \\
		Theorem~\ref{th:contour_integral}, $x_i^{(a)}\in\RR^p$ & $\pm$ 0.0242 & $\pm$ 0.0210 & $\pm$ 0.0160 & $\pm$ 0.0120 & $\pm$ 0.0089 & $\pm$ 0.0067 & $\pm$ 0.0051 & $\pm$ 0.0045 & $\pm$ 0.0046 \\
		& \\
		Proposed estimator & {\bf 0.0979}  &  {\bf 0.1455}  &  {\bf 0.1693}  &  {\bf 0.1811}  &  {\bf 0.1871}  &  {\bf 0.1902}  &  {\bf 0.1917}  &  {\bf 0.1926}  &  {\bf 0.1940} \\
		Corollary~\ref{cor:log2t}, $x_i^{(a)}\in\RR^p$ & $\pm$ 0.0242 & $\pm$ 0.0210 & $\pm$ 0.0160 & $\pm$ 0.0120 & $\pm$ 0.0089 & $\pm$ 0.0067 & $\pm$ 0.0051 & $\pm$ 0.0045 & $\pm$ 0.0046 \\
		& \\
		Traditional approach & 0.1024  &   0.1529  &  0.1826 &  0.2063 &  0.2364 &  0.2890  &  0.3954  &  0.6339 &  1.2717 \\
		$x_i^{(a)}\in\RR^p$ & $\pm$ 0.0242 & $\pm$ 0.0210 & $\pm$ 0.0160 & $\pm$ 0.0120 & $\pm$  0.0089 & $\pm$ 0.0068 & $\pm$ 0.0052 &  $\pm$  0.0048 & $\pm$ 0.0056 \\
		\\ \hline \\
		Proposed estimator & {\bf 0.0982}  &  {\bf 0.1455}  &  {\bf 0.1691}  &  {\bf 0.1811}  &  {\bf 0.1877}  &  {\bf 0.1901}  &  {\bf 0.1917}  &  {\bf 0.1922}  &  {\bf 0.1924} \\
		Theorem~\ref{th:contour_integral}, $x_i^{(a)}\in\CC^p$ & $\pm$ 0.0171 & $\pm$ 0.0145 & $\pm$ 0.0114 & $\pm$ 0.0082 & $\pm$ 0.0063 & $\pm$ 0.0046 & $\pm$ 0.0037 & $\pm$ 0.0028 & $\pm$ 0.0028 \\
                & \\
                Proposed estimator &  0.0968  &  0.1441 &  0.1675  &  0.1796 &  0.1861  &  0.1886  &  0.1903  &  0.1913  &  0.1931 \\
		Corollary~\ref{cor:log2t}, $x_i^{(a)}\in\CC^p$ & $\pm$ 0.0171 & $\pm$ 0.0145 & $\pm$ 0.0114 & $\pm$ 0.0082 & $\pm$ 0.0063 & $\pm$ 0.0046 & $\pm$ 0.0037 & $\pm$ 0.0028 & $\pm$ 0.0028 \\
                & \\
                Traditional approach & 0.1012  &   0.1515  &  0.1809 &  0.2048 &  0.2354 &  0.2873  &  0.3937  &  0.6318 &  1.2679 \\
                $x_i^{(a)}\in\CC^p$ & $\pm$ 0.0171 & $\pm$ 0.0146 & $\pm$ 0.0114 & $\pm$ 0.0082 & $\pm$ 0.0064 & $\pm$ 0.0047 & $\pm$ 0.0038 & $\pm$ 0.0030 & $\pm$ 0.0034
        \end{tabular}
	\bigskip
	\caption{Estimation of the Fisher distance $D_{\rm F}(C_1,C_2)$. Simulation example for real $x^{(a)}_i\sim \mathcal N(0,C_a)$ (top part) and complex $x^{(a)}_i\sim \mathcal {CN}(0,C_a)$ (bottom part), $[C_1^{-\frac12}C_2C_1^{-\frac12}]_{ij}=.3^{|i-j|}$, $n_1=1\,024$, $n_2=2\,048$, as a function of $p$. The ``$\pm$'' values correspond to one standard deviation of the estimates. Best estimates stressed in {\bf boldface} characters.}
\label{tab:1}
\end{table}

\bigskip

Technically speaking, our main result unfolds from a three-step approach: (i) relating the limiting (as $p,n_1,n_2\to\infty$ with $p/n_a=O(1)$) eigenvalue distribution of the (sometimes called Fisher) matrix $\hat C_1^{-1}\hat C_2$ to the limiting eigenvalue distribution of $C_1^{-1}C_2$ by means of a functional identity involving their respective {\it Stieltjes transforms} (see definition in Section~\ref{sec:mainresults}), then (ii) expressing the studied matrix distance as a complex integral featuring the Stieltjes transform and proceeding to successive change of variables to exploit the functional identity of (i), and finally (iii) whenever possible, explicitly evaluating the complex integral through complex analysis techniques. This approach is particularly reminiscent of the eigenvalue and eigenvector projection estimates proposed by Mestre in 2008 in a series of seminal articles \citep{MES08,MES08b}. In \citep{MES08b}, Mestre considers a single sample covariance matrix $\hat C$ setting and provides a complex integration approach to estimate the individual eigenvalues as well as eigenvector projections of the population covariance matrix $C$; there, step (i) follows immediately from the popular result \citep{SIL95} on the limiting eigenvalue distribution of large dimensional sample covariance matrices; step (ii) was then the (simple yet powerful) key innovation and step (iii) followed from a mere residue calculus. In a wireless communication-specific setting, the technique of Mestre was then extended in \citep{COU10b} in a model involving the product of Gram random matrices; there, as in the present work, step (i) unfolds from two successive applications of the results of \citep{SIL95}, step (ii) follows essentially the same approach as in \cite{MES08b} (yet largely simplified) and step (iii) is again achieved from residue calculus. As the random matrix $\hat C_1^{-1}\hat C_2$ may be seen as a single sample covariance matrix conditionally on $\hat C_1$, with $\hat C_1$ itself a sample covariance, in the present work, step (i) is obtained rather straightforwardly from applying twice the results from \citep{SIL95,CHO95}; the model $\hat C_1^{-1}\hat C_2$ is actually reminiscent of the so-called multivariate $F$-matrices, of the type $\hat C_1^{-1}\hat C_2$ but with $C_1=C_2$, extensively studied in \citep{SIL85,BAI88,ZHE12}; more recently, the very model under study here (that is, for $C_1\neq C_2$) was analyzed in \citep{WAN17,ZHE17} (step (i) of the present analysis is in particular consistent with Theorem~2.1 of \citep{ZHE17}, yet our proposed formulation is here more convenient to our purposes). 

But the main technical difficulty of the present contribution lies in steps (ii) and (iii) of the analysis. Indeed, as opposed to \citep{MES08b,COU10b}, the complex integrals under consideration here involve rather non-smooth functions, and particularly complex logarithms. Contour integrals of complex logarithms can in general not be treated through mere residue calculus. Instead, we shall resort here to an in-depth analysis of the so-called branch-cuts, corresponding to the points of discontinuity of the complex logarithm, as well as to the conditions under which valid integration contours can be defined. Once integrals are properly defined, an elaborate contour design then turns the study of the complex integral into that of real integrals. As already mentioned, in the particular case of the function $f(t)=\log^2(t)$, these real integrals result in a series of expressions involving the so-called dilogarithm function (see e.g., \citep{ZAG07}), the many properties of which will be thoroughly exploited to obtain our final results.

\bigskip

The remainder of the article presents our main result, that is the novel estimator, first under the form of a generic complex integral and then, for a set of functions met in the aforementioned classical matrix distances, under the form of a closed-form estimator.

\section{Model}

For $a\in\{1,2\}$, let $x^{(a)}_1,\ldots,x^{(a)}_{n_a}\in\RR^p$ (or $\CC^p$) be $n_a$ independent and identically distributed vectors of the form $x^{(a)}_i=C_a^{\frac12}\tilde x^{(a)}_i$ with $\tilde x^{(a)}_i\in\RR^p$ (respectively $\CC^p$) a vector of i.i.d.\@ zero mean, unit variance, and finite fourth order moment entries, and $C_a\in\RR^{p\times p}$ (respectively $\CC^{p\times p}$) positive definite.
We define the sample covariance matrix
\begin{align*}
	\hat{C}_a&\equiv \frac1{n_a}X_aX_a^\trans,\quad X_a=[x^{(a)}_1,\ldots,x^{(a)}_{n_a}].
\end{align*}

We will work under the following set of assumptions.
\begin{assumption}[Growth Rate]
	\label{ass:growth}
	For $a\in\{1,2\}$,
	\begin{enumerate}
		\item\label{ass:c_a} denoting $c_a\equiv \frac{p}{n_a}$, $c_a<1$ and $c_a\to c_a^\infty\in(0,1)$ as $p\to\infty$;\footnote{The reader must here keep in mind that $c_a$ is a function of $p$; yet, for readability and since this has little practical relevance, we do not make any explicit mention of this dependence. One may in particular suppose that $c_a=c_a^\infty$ for all valid $p$.}
		\item\label{ass:lim_eigs_Ca} $\lim\sup_p \max (\|C_a^{-1}\|,\|C_a\|)<\infty$ with $\|\cdot\|$ the operator norm;
		\item\label{ass:nu} there exists a probability measure $\nu$ such that $\nu_p\equiv \frac1p\sum_{i=1}^p{\bm\delta}_{\lambda_i(C_1^{-1}C_2)} \to \nu$ weakly as $p\to\infty$ (with $\lambda_i(A)$ the eigenvalues of matrix $A$ and ${\bm \delta}_x$ the atomic mass at $x$).\footnote{Here again, in practice, one may simply assume that $\nu=\nu_p$, the discrete empirical spectral distribution for some fixed dimension $p$.}
	\end{enumerate}
\end{assumption}

The main technical ingredients at the core of our derivations rely on an accurate control of the eigenvalues of $\hat C_1^{-1}\hat C_2$. In particular, we shall demand that these eigenvalues remain with high probability in a compact set. Item~\ref{ass:lim_eigs_Ca} and the finiteness of the $2+\varepsilon$ moments of the entries of $X_1$ and $X_2$ enforce this request (through the seminal results of Bai and Silverstein on sample covariance matrix models \citep{SIL95} and \citep{SIL98}). Item~\ref{ass:nu} can be relaxed but is mathematically convenient (in particular to ensure that $\mu_p\equiv \frac1p\sum_{i=1}^p{\bm \delta}_{\lambda_i(\hat C_1^{-1}\hat C_2)}$ has an almost sure limit) and practically inconsequential. 

Item~\ref{ass:c_a} deserves a deeper comment. Requesting that $p<n_1$ and $p<n_2$ is certainly demanding in some practical scarce data conditions. Yet, as we shall demonstrate subsequently, our proof approach relies on some key changes of variables largely involving the signs of $1-c_1$ and $1-c_2$. Notably, working under the assumption $c_1>1$ would demand a dramatic change of approach which we leave open to future work (starting with the fact that $\hat C_1^{-1}$ is no longer defined). The scenario $c_2>1$ is more interesting though. As we shall point out in a series of remarks, for some functionals $f$ having no singularity at zero, the value $\int fd\nu_p$ can still be reliably estimated; however, dealing with $f(t)=\log(t)$ or $f(t)=\log^2(t)$ (unfortunately at the core of all aforementioned distances and divergences $D_{\rm F}$, $D_{\rm B}$, $D_{\rm KL}$, and $D_{\rm R}^\alpha$) will not be possible under our present scheme.

\section{Main Results}
\label{sec:mainresults}

\subsection{Preliminaries}

For $f$ some complex-analytically extensible real function, our objective is to estimate
\begin{align*}
	\int f d\nu_p
\end{align*}
where we recall that
\begin{align*}
	\nu_p\equiv \frac1p\sum_{i=1}^p {\bm\delta}_{\lambda_i(C_1^{-1}C_2)}
\end{align*}
from the samples $x_i^{(1)}$'s and $x_i^{(2)}$'s in the regime where $n_1,n_2,p$ are all large and of the same magnitude. In particular, given the aforementioned applications, we will be interested in considering the cases where $f(t)=t$, $f(t)=\log(t)$, $f(t)=\log(1+st)$ for some $s>0$, or $f(t)=\log^2(t)$. Note in passing that, since $C_1^{-1}C_2$ has the same eigenvalues as $C_1^{-\frac12}C_2C_1^{-\frac12}$ (by Sylverster's identity), the eigenvalues of $C_1^{-1}C_2$ are all real positive.

It will be convenient in the following to define 
\begin{align*}
	\lambda_i\equiv \lambda_i(\hat{C}_1^{-1}\hat{C}_2),\quad 1\leq i\leq p
\end{align*}
and
\begin{align*}
	\mu_p &\equiv \frac1p\sum_{i=1}^p {\bm\delta}_{\lambda_i}.
\end{align*}
By the law of large numbers, it is clear that, as $n_1,n_2\to\infty$ with $p$ fixed, $\mu_p\asto \nu_p$ (in law) and thus (up to some support boundedness control for unbounded $f$), $\int fd\mu_p-\int fd\nu_p\asto 0$. The main objective here is to go beyond this simple result accounting for the fact that $n_1,n_2$ may not be large compared to $p$. In this case, under Assumption~\ref{ass:growth}, $\mu_p \not\to \nu$ and it is unlikely that for most $f$, the convergence $\int fd\mu_p - \int fd\nu_p\asto 0$ would still hold.

\bigskip

Our main line of arguments follows results from the field of random matrix theory and complex analysis. We will notably largely rely on the relation linking the \emph{Stieltjes transform} of several measures (such as $\nu_p$ and $\mu_p$) involved in the model. The Stieltjes transform $m_\mu$ of a measure $\mu$ is defined, for $z\in\CC\setminus {\rm Supp}(\mu)$ (with ${\rm Supp}(\mu)$ the support of $\mu$), as
\begin{align*}
	m_\mu(z) &\equiv \int \frac{d\mu(t)}{t-z}
\end{align*}
which is complex-analytic on its definition domain and in particular has complex derivative
\begin{align*}
	m'_\mu(z) &\equiv \int \frac{d\mu(t)}{(t-z)^2}.
\end{align*}
For instance, for a discrete measure $\mu\equiv \sum_{i=1}^p \alpha_i {\bm \delta}_{\lambda_i}$, $m_\mu(z)=\sum_{i=1}^p\frac{\alpha_i}{\lambda_i-z}$ and $m'_\mu(z)=\sum_{i=1}^p\frac{\alpha_i}{(\lambda_i-z)^2}$. 

\subsection{Generic results}

With the notations from the section above at hand, our main technical result is as follows.
\begin{theorem}[Estimation via contour integral]
	\label{th:contour_integral}
	Let Assumption~\ref{ass:growth} hold and $f:\CC\to \CC$ be analytic on $\{z\in\CC,\Re[z]>0\}$. Take also $\Gamma\subset \{z\in\CC,\Re[z]>0\}$ a (positively oriented) contour strictly surrounding $\cup_{p=1}^\infty{\rm Supp}(\mu_p)$ (this set is known to be almost surely compact). For $z\in\CC \setminus {\rm Supp}(\mu_p)$, define the two functions
	\begin{align*}
		\varphi_p(z) &\equiv z+c_1z^2m_{\mu_p}(z) \\
		\psi_p(z) &\equiv 1-c_2-c_2zm_{\mu_p}(z).
	\end{align*}
	Then, the following result holds
	\begin{align*}
		\int fd\nu_p - \frac1{2\pi\imath} \oint_{\Gamma} f\left( \frac{\varphi_p(z)}{\psi_p(z)} \right) \left( \frac{\varphi_p'(z)}{\varphi_p(z)} - \frac{\psi_p'(z)}{\psi_p(z)}\right)\frac{\psi_p(z)}{c_2} dz\asto 0.
	\end{align*}
\end{theorem}

\begin{remark}[Known $C_1$]
	For $C_1$ known, Theorem~\ref{th:contour_integral} is particularized by taking the limit $c_1\to 0$, i.e.,
	\begin{align*}
		\int fd\nu_p - \frac1{2\pi\imath} \oint_{\Gamma} f\left( \frac{z}{\psi_p(z)} \right) \left( \frac1z - \frac{\psi_p'(z)}{\psi_p(z)}\right)\frac{\psi_p(z)}{c_2} dz\asto 0
	\end{align*}
	where now $m_{\mu_p}(z)=\frac1p\sum_{i=1}^p\frac1{\lambda_i(C_1^{-1}\hat C_2)-z}$ and $m'_{\mu_p}(z)=\frac1p\sum_{i=1}^p\frac1{(\lambda_i(C_1^{-1}\hat C_2)-z)^2}$. Basic algebraic manipulations allow for further simplification, leading up to
	\begin{align*}
		\int fd\nu_p - \frac1{2\pi\imath c_2} \oint_{\Gamma} f\left( \frac{-1}{m_{\tilde \mu_p}(z)} \right)m'_{\tilde \mu_p}(z) z dz \asto 0
	\end{align*}
	where $\tilde\mu_p=c_2\mu_p+(1-c_2){\bm\delta}_0$ is the eigenvalue distribution of $\frac1{n_2}X_2^\trans C_1^{-1}X_2$ (and thus $m_{\tilde \mu_p}(z)=c_2m_{\mu_p}(z)-(1-c_2)/z$). Letting $g(z)=f(1/z)$ and $G(z)$ such that $G'(z)=g(z)$, integration by parts of the above expression further gives
	\begin{align*}
		\int fd\nu_p - \frac1{2\pi\imath c_2} \oint_{\Gamma} G\left( - m_{\tilde \mu_p}(z) \right) dz \asto 0.
	\end{align*}
	For instance, for $f(z)=\log^2(z)$, $G(z)=z(\log^2(z)-2\log(z)+2)$.
\end{remark}

\begin{remark}[Extension to the $c_2>1$ case]
	\label{rem:c2>1}
	Theorem~\ref{th:contour_integral} extends to the case $c_2>1$ for all $f:\CC\to\CC$ analytic on the whole set $\CC$. This excludes notably $f(z)=\log^k(z)$ for $k\geq 1$, but also $f(z)=\log(1+sz)$ for $s>0$. Yet, the analyticity request on $f$ can be somewhat relaxed. More precisely, Theorem~\ref{th:contour_integral} still holds when $c_2>1$ if there exists a $\Gamma$ as defined in Theorem~\ref{th:contour_integral} such that $f$ is analytic on the interior of the contour described by $(\varphi/\psi)(\Gamma)$, where $\varphi$ and $\psi$ are the respective almost sure limits of $\varphi_p$ and $\psi_p$ (see Appendix~\ref{sec:integral_form} for details). The main issue with the case $c_2>1$, as thoroughly detailed in Appendix~\ref{sec:contour}, is that, while $\Gamma\subset \{z\in\CC,\Re[z]>0\}$, the interior of $(\varphi/\psi)(\Gamma)$ necessarily contains zero. This poses dramatic limitations to the applicability of our approach for $f(z)=\log^k(z)$ for which we so far do not have a workaround. For $f(z)=\log(1+sz)$ though, we will show that there exist sufficient conditions on $s>0$ to ensure that $-1/s$ (the singularity of $z\mapsto\log(1+sz)$) is not contained within $(\varphi/\psi)(\Gamma)$, thereby allowing for the extension of Theorem~\ref{th:contour_integral} to $f(t)=\log(1+st)$.
\end{remark}

\subsection{Special cases}

While Theorem~\ref{th:contour_integral} holds for all well-behaved $f$ on $\Gamma$, a numerical complex integral is required in practice to estimate $\int fd\nu_p$. It is convenient, when feasible, to assess the approximating complex integral in closed form, which is the objective of this section. When $f$ is analytic in the inside of $\Gamma$, the integral can be estimated merely through a residue calculus. This is the case notably of polynomials $f(t)=t^k$. If instead $f$ exhibits singularities in the inside of $\Gamma$, as for $f(t)=\log^k(t)$, more advanced contour integration arguments are required. Of utmost interest are the following results.

\begin{corollary}[Case $f(t)=t$]
	\label{cor:t}
	Under the conditions of Theorem~\ref{th:contour_integral},
	\begin{align*}
		\int t d\nu_p(t) - (1-c_1)\int td\mu_p(t) &\asto 0.
	\end{align*}
	and in the case where $c_1\to 0$, this is simply $\int t d\nu_p(t) - \int td\mu_p(t) \asto 0$.
\end{corollary}

As such, the classical sample covariance matrix estimator $\int td\mu_p(t)$ needs only be corrected by a product with $(1-c_1)$.

\begin{corollary}[Case $f(t)=\log(t)$]
	\label{cor:logt}
	Under the conditions of Theorem~\ref{th:contour_integral},
	\begin{align*}
		\int \log(t) d\nu_p(t) - \left[ \int \log(t)d\mu_p(t) -\frac{1-c_1}{c_1}\log(1-c_1)+\frac{1-c_2}{c_2}\log(1-c_2) \right] &\asto 0.
	\end{align*}
	When $c_1\to 0$, $\int \log(t) d\nu_p(t) - \left[ \int \log(t)d\mu_p(t) +\frac{1-c_2}{c_2}\log(1-c_2) + 1 \right] \asto 0$.
\end{corollary}

Note interestingly that, for $f(t)=\log(t)$ and $c_1=c_2$, the standard estimator is asymptotically $p,n_1,n_2$-consistent. This is no longer true though for $c_1\neq c_2$ but only a fixed bias is induced.

\begin{corollary}[Case $f(t)=\log(1+st)$]
	\label{cor:log1st}
	Under the conditions of Theorem~\ref{th:contour_integral}, let $s>0$ and denote $\kappa_0$ the unique \emph{negative} solution to $1+s\frac{\varphi_p(x)}{\psi_p(x)}=0$. Then we have 
	\begin{align*}
		&\int \log(1+st) d\nu_p(t)
		-\left[\frac{c_1+c_2-c_1c_2}{c_1c_2} \log\left( \frac{c_1+c_2-c_1c_2}{(1-c_1)(c_2-sc_1\kappa_0)} \right) \right. \\ &\left.+ \frac1{c_2}\log \left( -s\kappa_0(1-c_1) \right) + \int \log\left(1-\frac{t}{\kappa_0}\right)d\mu_p(t) \right] \asto 0.
	\end{align*}
	In the case where $c_1\to 0$, this is simply 
	\begin{align*}
		\int \log(1+st) d\nu_p(t) - \left[ \frac{1+s\kappa_0+\log(-s\kappa_0)}{c_2} + \int\log\left( 1- \frac{t}{\kappa_0}\right)d\mu_p(t) \right] \asto 0.
	\end{align*}
	
\end{corollary}
\begin{remark}[Limit when $s\to \infty$ and $s\to 0$]
	\label{rem:limits_over_s}
	It is interesting to note that, as $s\to \infty$, $\kappa_0\sim -(1-c_2)/s$. Plugging this into the expression above, we find that
	\begin{align*}
		&\int \log(1+st) d\nu_p(t) - \log(s) \nonumber \\
		&\sim_{s\to\infty} \int \log(t)d\mu_p(t) +\frac{1-c_2}{c_2}\log(1-c_2)-\frac{1-c_1}{c_1}\log(1-c_1)
	\end{align*}
	therefore recovering, as one expects, the result from Corollary~\ref{cor:logt}.

	Checking similarly the case where $s\to 0$ demands a second-order expansion of $\kappa_0$. It is first clear that $\kappa_0$ must grow unbounded as $s\to 0$, otherwise there would be a \emph{negative} solution to $\psi_p(x)=0$, which is impossible. This said, we then find that $\kappa_0\sim -\frac1{(1-c_1)s}+\frac{c_1+c_2-c_1c_2}{1-c_1}\frac1p\sum_{i=1}^p\lambda_i+o(1)$. Plugging this into Corollary~\ref{cor:log1st}, we find
	\begin{align*}
		\frac1s \int \log(1+st) d\nu_p(t) \sim_{s\to 0} (1-c_1)\int td\mu_p(t)
	\end{align*}
	therefore recovering the results from Corollary~\ref{cor:t}, again as expected.
\end{remark}

\begin{remark}[Location of $\kappa_0$]
	For numerical purposes, it is convenient to easily locate $\kappa_0$. Using the fact that, by definition,
	\begin{align*}
		\kappa_0 &= \frac{1-c_2-c_2\kappa_0m_{\mu_p}(\kappa_0)}{-s(1+c_1\kappa_0m_{\mu_p}(\kappa_0))}
	\end{align*}
	and the bound $-1<xm_{\mu_p}(x)<0$, for $x<0$, we find that
	\begin{align*}
		-\frac1{s(1-c_1)}<\kappa_0<0.
	\end{align*}
	As such, $\kappa_0$ may be found by a dichotomy search on the set $(-1/(s(1-c_1)),0)$.
\end{remark}


\begin{remark}[The case $c_2>1$]
	\label{rem:log1st_c2>1}
	Theorem~\ref{th:contour_integral} and Corollary~\ref{cor:log1st} may be extended to the case $c_2>1$, however with important restrictions. Precisely, let $\mu^-\equiv \inf \{{\rm Supp}(\mu)\}$ for $\mu$ the almost sure weak limit of the empirical measure $\mu_p$, and let $x^-\equiv \lim_{x\uparrow \mu^-}\frac{\varphi(x)}{\psi(x)}$, for $\varphi(z)\equiv z(1+c_1^\infty zm_\mu(z))$ and $\psi(z)\equiv 1-c_2^\infty -c_2^\infty zm_\mu(z)$ the respective almost sure functional limits of $\varphi_p$ and $\psi_p$. Also redefine $\kappa_0$ in Corollary~\ref{cor:log1st} as the smallest real (non-necessarily negative) solution to $1+s\frac{\varphi_p(x)}{\psi_p(x)}=0$. Then, for all $s>0$ satisfying
	\begin{align*}
		1+sx^->0
	\end{align*}
	we have the following results
	\begin{enumerate}
		\item Theorem~\ref{th:contour_integral} holds true, however for a contour $\Gamma$ having a leftmost real crossing within the set $(\kappa_0,\mu^-)$;
		\item Corollary~\ref{cor:log1st} extends to
			\begin{align*}
				&\int \log(1+st) d\nu_p(t) -\left[\frac{c_1+c_2-c_1c_2}{c_1c_2} \log\left( \frac{c_1+c_2-c_1c_2}{(1-c_1)|c_2-sc_1\kappa_0|} \right) \right. \\ &\left.+ \frac1{c_2}\log \left| -s\kappa_0(1-c_1) \right| + \int \log\left|1-\frac{t}{\kappa_0}\right|d\mu_p(t) \right] \asto 0.
			\end{align*}
	\end{enumerate}
	For instance, for $C_1=C_2$, i.e., $C_1^{-1}C_2=I_p$, we have
	\begin{align*}
		\mu^-=\left(\frac{1 - \sqrt{c_1^\infty +c_2^\infty -c_1^\infty c_2^\infty } }{1-c_1^\infty }\right)^2,\quad x^-=\frac{1 - \sqrt{c_1^\infty +c_2^\infty -c_1^\infty c_2^\infty } }{1-c_1^\infty }
	\end{align*}
	so that $x^-<0$ when $c_2^\infty >1$; so, there, Theorem~\ref{th:contour_integral} and Corollary~\ref{cor:log1st} hold true as long as
	\begin{align}
		\label{eq:bound_s}
		0< s < \frac{1-c_1^\infty }{\sqrt{c_1^\infty +c_2^\infty -c_1^\infty c_2^\infty }-1}.
	\end{align}
\end{remark}

It is clear from the second part of the remark that, as $c_2$ increases, the range of valid $s$ values vanishes. Perhaps surprinsingly though, as $c_1^\infty\uparrow 1$, the range of valid $s$ values converges to the fixed set $(0,2/(c_2-1))$ and thus does not vanish.


\bigskip

As opposed to the previous scenarios, for the case $f(t)=\log^2(t)$, the exact form of the integral from Theorem~\ref{th:contour_integral} is non-trivial and involves dilogarithm functions (see its expression in \eqref{eq:log2t_exact} in the Appendix). This involved expression can nonetheless be significantly simplified using a large-$p$ approximation, resulting in an estimate only involving usual functions, as shown subsequently. 

\begin{corollary}[Case $f(t)=\log^2(t)$]
	\label{cor:log2t}
	Let $0<\eta_1<\ldots<\eta_p$ be the eigenvalues of $\Lambda - \frac{\sqrt{\lambda}\sqrt{\lambda}^\trans}{p-n_1}$ and $0<\zeta_1<\ldots<\zeta_p$ the eigenvalues of $\Lambda - \frac{\sqrt{\lambda}\sqrt{\lambda}^\trans}{n_2}$, where $\lambda=(\lambda_1,\ldots,\lambda_p)^\trans$, $\Lambda={\rm diag}(\lambda)$, and $\sqrt{\lambda}$ is understood entry-wise. Then, under the conditions of Theorem~\ref{th:contour_integral},
	\begin{align*}
		&\int \log^2(t)d\nu_p(t) \nonumber \\
		&- \left[ \frac1p\sum_{i=1}^p\log^2( (1-c_1)\lambda_i ) + 2\frac{c_1+c_2-c_1c_2}{c_1c_2} \left\{ \left( \Delta_\zeta^\eta \right)^\trans M \left( \Delta_\lambda^\eta \right) + (\Delta_{\lambda}^\eta)^\trans r  \right\}  \right. \\
		&\left. - \frac2p\left( \Delta_\zeta^\eta \right)^\trans N1_p - 2\frac{1-c_2}{c_2}\left\{ \frac12\log^2( (1-c_1)(1-c_2) ) + (\Delta_{\zeta}^\eta)^\trans r\right\} \right] \asto 0
	\end{align*}
	where we defined $\Delta_a^b$ the vector with $(\Delta_a^b)_i=b_i-a_i$ and, for $i,j\in\{1,\ldots,p\}$, $r_i = \frac{\log( (1-c_1)\lambda_i )}{\lambda_i}$ and
	\begin{align*}
		M_{ij} = \left\{ \begin{array}{ll} \frac{\frac{\lambda_i}{\lambda_j}-1-\log\left(\frac{\lambda_i}{\lambda_j}\right)}{(\lambda_i-\lambda_j)^2}&,~i\neq j \\ \frac1{2\lambda_i^2}&,~i=j \end{array}\right., \quad
		N_{ij} = \left\{ \begin{array}{ll} \frac{\log\left(\frac{\lambda_i}{\lambda_j}\right)}{\lambda_i-\lambda_j}&,~i\neq j \\ \frac1{\lambda_i}&,~i=j. \end{array}\right.
	\end{align*}
	In the limit $c_1\to 0$ (i.e., for $C_1$ known), this becomes
	\begin{align*}
		&\int \log^2(t)d\nu_p(t) - \left[ \frac1p\sum_{i=1}^p\log^2( \lambda_i ) + \frac2p\sum_{i=1}^p \log( \lambda_i ) - \frac2p\left( \Delta_\zeta^\lambda \right)^\trans Q1_p \nonumber\right. \\ &\left. - 2\frac{1-c_2}{c_2}\left\{ \frac12\log^2( 1-c_2 ) + \left( \Delta_\zeta^\lambda \right)^\trans q\right\} \right] \asto 0
	\end{align*}
	with
	\begin{align*}
		Q_{ij} = \left\{ \begin{array}{ll} \frac{ \lambda_i \log \left( \frac{\lambda_i}{\lambda_j} \right) - (\lambda_i-\lambda_j) }{(\lambda_i-\lambda_j)^2}&,~i\neq j \\ \frac1{2\lambda_i}&,~i=j \end{array}\right., \textmd{ and }
		q_i = \frac{\log( \lambda_i )}{\lambda_i}.
	\end{align*}
\end{corollary}

\section{Experimental Results}

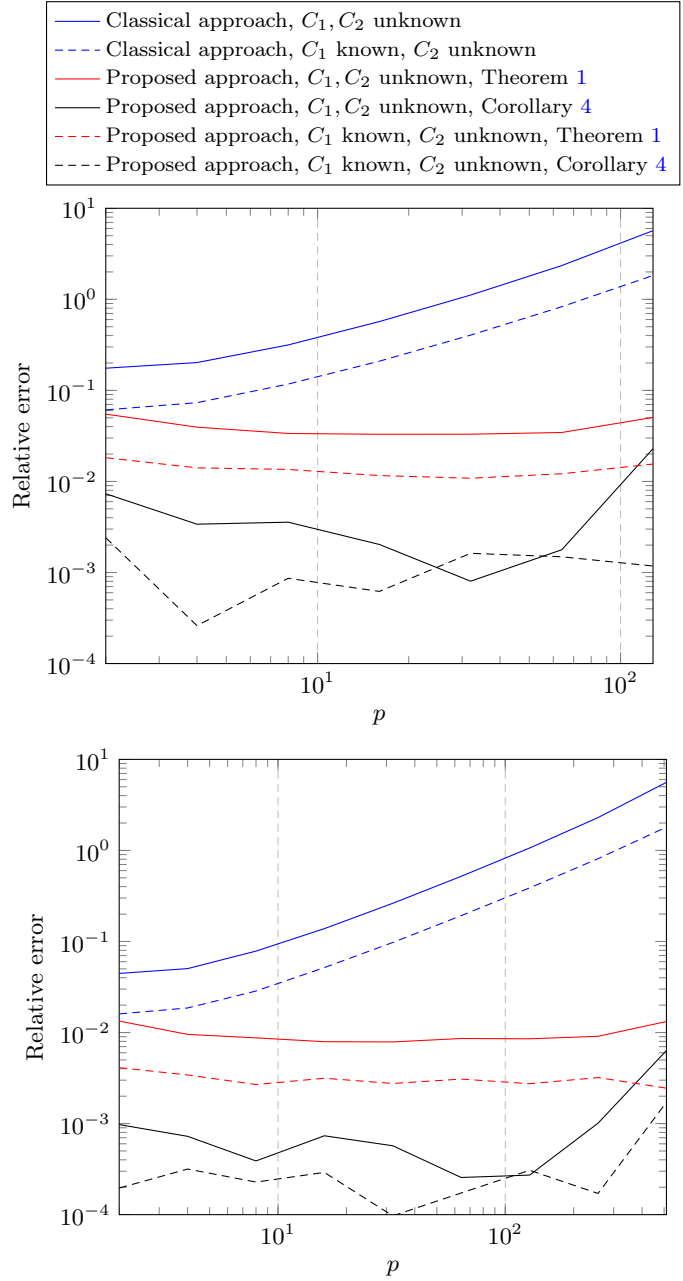
\begin{figure}
	\centering
	\begin{tabular}{cc}

	\begin{tikzpicture}[font=\footnotesize]
		\renewcommand{\axisdefaulttryminticks}{4} 
		\tikzstyle{every major grid}+=[style=densely dashed]       
		\tikzstyle{every axis legend}+=[cells={anchor=west},fill=white,
		at={(1.05,1.05)}, anchor=south east, font=\scriptsize ]
		\begin{loglogaxis}[
				width=.7\linewidth,
				xmin=2,
				ymin=1e-4,
				xmax=128,
				ymax=10,
				grid=major,
				ymajorgrids=false,
				scaled ticks=true,
				xlabel={$p$},
				ylabel={Relative error}
			]
			\addplot[blue] plot coordinates{
				(2,0.175335)(4,0.202115)(8,0.315768)(16,0.569918)(32,1.114142)(64,2.349436)(128,5.675005)

			};
			\addplot[blue,densely dashed] plot coordinates{
				(2,0.061083)(4,0.073352)(8,0.117337)(16,0.209262)(32,0.404634)(64,0.830053)(128,1.828267)
			};
			\addplot[red] plot coordinates{
				(2,0.054886)(4,0.039448)(8,0.033671)(16,0.032942)(32,0.033024)(64,0.034492)(128,0.050486)
			};
			\addplot[black] plot coordinates{
				(2,0.007327)(4,0.003393)(8,0.003570)(16,0.002032)(32,0.000803)(64,0.001778)(128,0.022893)
			};
			\addplot[red,densely dashed] plot coordinates{
				(2,0.018218)(4,0.014107)(8,0.013543)(16,0.011604)(32,0.010851)(64,0.012145)(128,0.015535)
			};
			\addplot[black,densely dashed] plot coordinates{
				(2,0.002408)(4,0.000261)(8,0.000864)(16,0.000620)(32,0.001624)(64,0.001488)(128,0.001178)
			};
			\legend{ {Classical approach, $C_1,C_2$ unknown},{Classical approach, $C_1$ known, $C_2$ unknown},{Proposed approach, $C_1,C_2$ unknown, Theorem~\ref{th:contour_integral}},{Proposed approach, $C_1,C_2$ unknown, Corollary~\ref{cor:log2t}},{Proposed approach, $C_1$ known, $C_2$ unknown, Theorem~\ref{th:contour_integral}},{Proposed approach, $C_1$ known, $C_2$ unknown, Corollary~\ref{cor:log2t}} }
		\end{loglogaxis}
	\end{tikzpicture}
	\\
	\begin{tikzpicture}[font=\footnotesize]
		\renewcommand{\axisdefaulttryminticks}{4} 
		\tikzstyle{every major grid}+=[style=densely dashed]       
		\tikzstyle{every axis legend}+=[cells={anchor=west},fill=white,
		at={(1,1.02)}, anchor=south east, font=\scriptsize ]
		\begin{loglogaxis}[
				width=.7\linewidth,
				xmin=2,
				ymin=1e-4,
				xmax=512,
				ymax=10,
				grid=major,
				ymajorgrids=false,
				scaled ticks=true,
				xlabel={$p$},
				ylabel={Relative error}
			]
			\addplot[blue] plot coordinates{
				(2,0.044596)(4,0.050369)(8,0.078438)(16,0.138463)(32,0.262826)(64,0.519936)(128,1.063538)(256,2.294995)(512,5.597887)
			};
			\addplot[blue,densely dashed] plot coordinates{
				(2,0.016020)(4,0.018618)(8,0.028583)(16,0.051710)(32,0.097550)(64,0.192284)(128,0.388004)(256,0.811822)(512,1.803616)
			};
			\addplot[red] plot coordinates{
				(2,0.013367)(4,0.009545)(8,0.008740)(16,0.007935)(32,0.007889)(64,0.008603)(128,0.008536)(256,0.009099)(512,0.013216)
			};
			\addplot[black] plot coordinates{
				(2,0.000975)(4,0.000725)(8,0.000389)(16,0.000735)(32,0.000571)(64,0.000256)(128,0.000272)(256,0.001010)(512,0.006352)
			};
			\addplot[red,densely dashed] plot coordinates{
				(2,0.004114)(4,0.003421)(8,0.002688)(16,0.003148)(32,0.002760)(64,0.003082)(128,0.002743)(256,0.003203)(512,0.002450)
			};
			\addplot[black,densely dashed] plot coordinates{
				(2,0.000196)(4,0.000317)(8,0.000228)(16,0.000291)(32,0.000097)(64,0.000173)(128,0.000307)(256,0.000171)(512,0.001716)
			};
		\end{loglogaxis}
	\end{tikzpicture}
\end{tabular}
\caption{Estimation of the squared Fisher distance $D^2_{\rm F}(C_1,C_2)$. Relative estimation error in logarithmic scale for $x_i^{(a)}\sim \mathcal N(0,C_a)$ with $[C_1^{-\frac12}C_2C_1^{-\frac12}]_{ij}=.3^{|i-j|}$; (top) $n_1=256$, $n_2=512$, and (bottom) $n_1=1024$, $n_2=2048$, varying $p$. }
	\label{fig:error}
\end{figure}

This section presents a series of experimental verifications, completing Table~\ref{tab:1} in Section~\ref{sec:intro}. In all cases presented below, we consider the squared Fisher distance (i.e., the case $f(t)=\log^2(t)$) between $C_1$ and $C_2$, where $C_1^{-\frac12}C_2C_1^{-\frac12}$ is a Toeplitz matrix with entry $(i,j)$ equal to $.3^{|i-j|}$. Figure~\ref{fig:error} displays the normalized estimation error (i.e., the absolute difference between genuine and estimated squared distance over the genuine squared distance) for $x_i^{(a)}$ real $\mathcal N(0,C_a)$, for varying values of the size $p$ of $C_1$ and $C_2$ but for $n_1$ and $n_2$ fixed ($n_1=256$ and $n_2=512$ in the left-hand display, $n_1=1024$, $n_2=2048$ in the right-hand display); results are averaged over $10\,000$ independent realizations. We compare the classical ($n_1,n_2$-consistent estimator) to the proposed $n_1,n_2,p$-consistent estimators obtained from both Theorem~\ref{th:contour_integral} (or equivalently here from the closed-form expression \eqref{eq:log2t_exact} form the Appendix) and from Corollary~\ref{cor:log2t}. We also add in dashed lines the distance estimate for $C_1$ a priori known. 

It is seen from the figures, provided in log scale, that the relative error of the standard $n_1,n_2$-consistent approach diverges with $p$, while the error of the proposed estimator remains at a low, rather constant, value. As already observed in Table~\ref{tab:1}, the large-$p$ approximation formula from Corollary~\ref{cor:log2t} even exhibits a better behavior than the exact integral expression from Theorem~\ref{th:contour_integral} in this real Gaussian setup.

\bigskip

\begin{figure}
	\centering
	\begin{tikzpicture}[font=\footnotesize]
		\renewcommand{\axisdefaulttryminticks}{4} 
		\tikzstyle{every major grid}+=[style=densely dashed]       
		\tikzstyle{every axis legend}+=[cells={anchor=west},fill=white,
		at={(1.05,1.05)}, anchor=south east, font=\scriptsize ]
		\begin{loglogaxis}[
				width=.7\linewidth,
				xmin=2,
				ymin=1e-3,
				xmax=256,
				ymax=100,
				grid=major,
				ymajorgrids=false,
				scaled ticks=true,
				xlabel={$p$},
				ylabel={Relative error}
			]
			\addplot[blue] plot coordinates{
				(2,18.981856)(4,10.262743)(8,7.822954)(16,6.785695)(32,6.321351)(64,6.111228)(128,6.005908)(256,5.957551)
			};
			\addplot[blue,densely dashed] plot coordinates{
				(2,15.430718)(4,8.222204)(8,6.219797)(16,5.365422)(32,4.992780)(64,4.809957)(128,4.723409)(256,4.685768)
			};
			\addplot[red] plot coordinates{
				(2,7.191610)(4,2.259119)(8,0.970607)(16,0.436049)(32,0.200623)(64,0.100498)(128,0.051152)(256,0.028008)
			};
			\addplot[black] plot coordinates{
				(2,1.740714)(4,0.197307)(8,0.010512)(16,0.033481)(32,0.033419)(64,0.016403)(128,0.007322)(256,0.001202)
			};
			\addplot[red,densely dashed] plot coordinates{
				(2,5.514679)(4,1.763428)(8,0.763615)(16,0.342722)(32,0.165627)(64,0.078972)(128,0.039809)(256,0.023684)
			};
			\addplot[black,densely dashed] plot coordinates{
				(2,1.346736)(4,0.058204)(8,0.053929)(16,0.065687)(32,0.039560)(64,0.024329)(128,0.012128)(256,0.002312)
			};
			\legend{ {Classical approach, $C_1,C_2$ unknown},{Classical approach, $C_1$ known, $C_2$ unknown},{Proposed approach, $C_1,C_2$ unknown, Theorem~\ref{th:contour_integral}},{Proposed approach, $C_1,C_2$ unknown, Corollary~\ref{cor:log2t}},{Proposed approach, $C_1$ known, $C_2$ unknown, Theorem~\ref{th:contour_integral}},{Proposed approach, $C_1$ known, $C_2$ unknown, Corollary~\ref{cor:log2t}} }
			\draw[-] (axis cs:32,0.200623) -| (axis cs:64,0.100498);
			\node at (axis cs:60,0.3) {slope $\simeq -1$};
		\end{loglogaxis}
	\end{tikzpicture}
	\caption{Estimation of the squared Fisher distance $D^2_{\rm F}(C_1,C_2)$. Relative estimation error for $x_i^{(a)}\sim \mathcal N(0,C_a)$ with $[C_1^{-\frac12}C_2C_1^{-\frac12}]_{ij}=.3^{|i-j|}$, $\frac{p}{n_1}=\frac14$, $\frac{p}{n_2}=\frac12$, varying $p$. }
	\label{fig:error_constant_c's}
\end{figure}
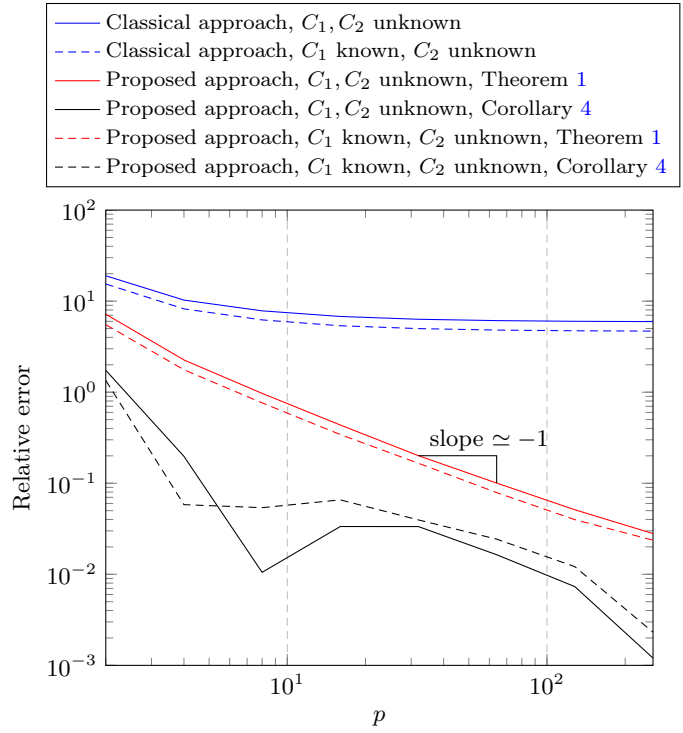

In the same setting, Figure~\ref{fig:error_constant_c's} subsequently shows the performances in terms of relative error for \emph{simultaneously} growing $n_1,n_2,p$, with a constant ratio $p/n_1=1/4$ and $p/n_2=1/2$. As expected from our theoretical analysis, the proposed estimators display a vanishing error as $n_1,n_2,p\to\infty$. The standard estimator, on the opposite, shows a high saturating relative error.

\section{Concluding Remarks}
\label{sec:conclusion}

As pointed out in the introduction, the series of estimators derived in this article allow the practitioner to assess a large spectrum of covariance matrix distances or divergences based on sample estimators when the number $n_1,n_2$ of available data is of similar order of magnitude as the size $p$ of the individual samples. For instance, the R\'enyi divergence \eqref{eq:Renyi} between $\mathcal N(0,C_2)$ and $\mathcal N(0,C_1)$ can be expressed as a linear combination of the estimator of Corollary~\ref{cor:log1st} with $s=\frac{1-\alpha}{\alpha}$ and of the estimator of Corollary~\ref{cor:logt}. In particular, the results from Corollaries~\ref{cor:t}--\ref{cor:log2t} are sufficient to cover all matrix distances and divergences discussed in the introductory Section~\ref{sec:intro}. For other distances involving ``smooth'' linear functionals of the eigenvalues of $C_1^{-1}C_2$, the generic result from Theorem~\ref{th:contour_integral} can be used to retrieve (at least) a numerical estimate of the sought-for distance or divergence. Applications of the present estimators to specific contexts of signal and data engineering at large, and specifically to the areas of statistical signal processing, artificial intelligence, and machine learning, are numerous.

\bigskip

The analysis of the slope in the log-log plot of Figure~\ref{fig:error_constant_c's} of the relative error for the proposed estimator reveals that, for fixed ratios $p/n_1$ and $p/n_2$, the error behaves as $p^{-1}$ (or equivalent as $n_1^{-1}$ or $n_2^{-1}$). This suggests that, as opposed to classical $n$-consistent estimators for which optimal fluctuations of the estimators are usually expected to be of order $n^{-\frac12}$, the proposed $n_1,n_2,p$-consistent estimator exhibits a quadratically faster convergence speed. This observation is in fact consistent with previous findings, such as \citep{YAO11} which demonstrated a central limit theorem with speed $p^{-1}$ for the eigenvalue functional estimators proposed by Mestre in \citep{MES08b}. This is also reminiscent from the numerous central limit theorems of eigenvalue functionals derived in random matrix theory \citep{HAC08,NAJ13,ZHE17} since the early findings from Bai and Silverstein \citep{BAI04}. This fast convergence speed partly explains the strong advantage of the proposed estimator over standard large-$n$ alone estimates, even for very small ratios $p/n$, as shown in the various simulations provided in the article.

For practical purposes, it would certainly be convenient to obtain a result similar to \citep{YAO11} in the present context, that is a central limit theorem for the fluctuations of the estimator obtained in Theorem~\ref{th:contour_integral}. This would allow practitioners to access both a consistent estimator for their sought-for matrix distance \emph{as well as} a confidence margin. This investigation demands even more profound calculi (as can be seen from the detailed derivations of \citep{YAO11}) and is left to future work.

Also, recall from Table~\ref{tab:1} that the estimator from Theorem~\ref{th:contour_integral} has a much better behavior on average in the complex Gaussian rather than in the real Gaussian case. This observation is likely due to a systematic bias of order $O(1/p)$ which is absent in the complex case. Indeed, the results from \citep{ZHE17} show that the difference $p(m_{\mu_p}(z)-m_\mu(z))$ satisfies a central limit theorem with in general non-zero mean apart from the complex Gaussian case (or, to be more exact, apart from the cases where the $x_i^{(a)}$'s have complex independent entries with zero mean, unit variance, and zero kurtosis). Coupling this result with \citep{YAO11} strongly suggests that the proposed estimators in the present article exhibit a systematic order-$p^{-1}$ bias but in the complex Gaussian case. As this bias is likely itself prone to estimation from the raw data, this observation opens the door to a further improvement of the estimator in Theorem~\ref{th:contour_integral} that would discard the bias. Proceeding as such would not change the order of magnitude of the error (still of order $1/p$) but reduces the systematic part of the error.

\bigskip

A last important item to be discussed at this point lies in the necessary condition $n_1>p$ and $n_2>p$ in the analysis. We have shown in the proof of Theorem~\ref{th:contour_integral} (in Appendices~\ref{sec:integral_form} and \ref{sec:contour}) that the requirements $c_1<1$ and $c_2<1$ are both mandatory for our estimation approach to remain valid on a range of functions $f$ analytic on $\{z\in\CC,\Re[z]>0\}$ (which notably includes here logarithm functions). Yet, as discussed throughout the article, while $n_1>p$ is mandatory for our proof approach to remain valid, the constraint $n_2>p$ can be relaxed to some extent. Yet, this excludes functions $f$ that are not analytic in a neighborhood of zero, thereby excluding functions such as powers of $1/z$ or of $\log(z)$. More advanced considerations, and possibly a stark change of approach, are therefore demanded to retrieve consistent estimators when $p>n_2$ for these functions. If one resorts to projections, dimension reductions, or regularization techniques to obtain an invertible ersatz for $\hat C_1$, one may even allow for $p>n_1$, but this would dramatically change the present analysis. As such, the quest for $n_1,n_2,p$-consistent estimators of the matrix distances when either $n_1<p$ or $n_2<p$ also remains an interesting open research avenue.

\appendix

\section*{Appendices}

We provide here the technical developments for the proof of Theorem~\ref{th:contour_integral} as well as all subsequent corollaries (Corollaries~\ref{cor:t}--\ref{cor:log2t}). 

The appendix is structured as follows: Appendix~\ref{sec:integral_form} provides the proof of Theorem~\ref{th:contour_integral} following the same approach as in \citep{COU10b}, relying mostly on the results from \citep{SIL95,CHO95}. Appendix~\ref{sec:integral_calculus} then provides the technical details of the calculi behind Corollaries~\ref{cor:t}--\ref{cor:log2t}; this is undertaken through a first thorough characterization of the singular points of $\varphi_p$ and $\psi_p$ and functionals of these (these singular points are hereafter denoted $\lambda_i$, $\eta_i$, $\zeta_i$ and $\kappa_i$), allowing for a proper selection of the integration contour, and subsequently through a detailed calculus for all functions $f(t)$ under study. Appendix~\ref{sec:contour} discusses in detail the question of the position of the complex contours when affected by change of variables. Finally, Appendix~\ref{sec:case_c2>1} provides some analysis of the extension of Theorem~\ref{th:contour_integral} and the corollaries to the $c_2>1$ scenario.

\setcounter{section}{0}
\section{Integral Form}
\label{sec:integral_form}

\subsection{Relating $m_\nu$ to $m_\mu$}
We start by noticing that we may equivalently assume the following setting:
\begin{itemize}
	\item $x^{(1)}_1,\ldots,x^{(1)}_{n_1}\in\RR^p$ vectors of i.i.d.\@ zero mean and unit variance entries
	\item $x^{(2)}_1,\ldots,x^{(2)}_{n_2}\in\RR^p$ of the form $x^{(2)}_i=C^{\frac12}\tilde x^{(2)}_i$ with $\tilde x^{(2)}_i\in\RR^p$ a vector of i.i.d.\@ zero mean and unit variance entries
\end{itemize}
where $C\equiv C_1^{-\frac12}C_2C_1^{-\frac12}$. 

Indeed, with our first notations, $\hat{C}_1^{-1}\hat{C}_2=\frac1{n_1}C_1^{-\frac12}\tilde X_1\tilde X_1^\trans C_1^{-\frac12}\frac1{n_2}C_2^{\frac12}\tilde X_2\tilde X_2^\trans C_2^{\frac12}$ (here $\tilde X_a=[\tilde x_1^{(a)},\ldots,\tilde x_{n_a}^{(a)}]$), which has the same spectrum as the matrix $(\frac1{n_1}\tilde X_1\tilde X_1^\trans) (\frac1{n_2} C_1^{-\frac12}C_2^{\frac12}\tilde X_2\tilde X_2^\trans C_2^{\frac12}C_1^{-\frac12})$ and we may then consider that the $x_i^{(1)}$'s actually have covariance $I_p$, while the $x_i^{(2)}$'s have covariance $C=C_1^{-\frac12}C_2C_1^{-\frac12}$, without altering the spectra under study. With these new definitions, we first condition with respect to the $x^{(2)}_i$'s, and study the spectrum of $\hat{C}_1^{-1}\hat{C}_2$, which is the same as that of $\hat{C}_2^{\frac12}\hat{C}_1^{-1}\hat{C}_2^{\frac12}$. A useful remark is the fact that $\hat{C}_2^{\frac12}\hat{C}_1^{-1}\hat{C}_2^{\frac12}$ is the ``inverse spectrum'' of $\hat{C}_2^{-\frac12}\hat{C}_1\hat{C}_2^{-\frac12}$, which is itself the same spectrum as that of $\frac1{n_1}X_1^\trans \hat{C}_2^{-1}X_1$ except for $n_1-p$ additional zero eigenvalues.

Denoting $\tilde{\mu}_{p}^{-1}$ the eigenvalue distribution of $\frac1{n_1}\tilde{X}_1^\trans \hat{C}_{2}^{-1}\tilde{X}_{1}$, we first know from \citep{SIL95} that, under Assumption~\ref{ass:growth}, as $p\to\infty$, $\tilde{\mu}_{p}^{-1}\asto \tilde{\mu}^{-1}$, where $\tilde{\mu}^{-1}$ is the probability measure with Stieltjes transform $m_{\tilde \mu^{-1}}$ defined as the unique (analytical function) solution to
\begin{equation*}
	m_{\tilde{\mu}^{-1}}(z)=\left(-z+c_1^\infty\int \frac{td\xi_{2}^{-1}(t)}{1+tm_{\tilde{\mu}^{-1}}(z)}\right)^{-1}
\end{equation*}
with $\xi_2$ the almost sure limiting spectrum distribution of $\hat{C}_2$ and $m_{\xi_2}$ its associated Stieltjes transform (note importantly that, from \citep{SIL95} and Assumption~\ref{ass:growth}, $\xi_2$ has bounded support and is away from zero). Recognizing a Stieltjes transform from the right-hand side integral, this can be equivalently written
\begin{equation}
	\label{eq:baisilverstein1}
	m_{\tilde{\mu}^{-1}}(z)=\left(-z+\frac{c_1^\infty}{m_{\tilde{\mu}^{-1}}(z)} - \frac{c_1^\infty}{m_{\tilde{\mu}^{-1}}(z)^2} m_{\xi_2^{-1}}\left( -\frac1{m_{\tilde{\mu}^{-1}}(z)} \right)\right)^{-1}.
\end{equation}

Accounting for the aforementioned additional zero eigenvalues, $\tilde{\mu}^{-1}$ is related to $\mu^{-1}$, the almost sure limiting spectrum distribution of $\hat{C}_2^{-1/2}\hat{C}_1\hat{C}_2^{-1/2}$, through the relation $\tilde{\mu}^{-1}=c^\infty_1\mu^{-1}+(1-c_1^\infty){\bm\delta}_{0}$ with ${\bm\delta}_x$ the Dirac measure at $x$ and we have
\begin{equation*}
m_{\tilde{\mu}^{-1}}(z)=c^\infty_1m_{\mu^{-1}}(z)-(1-c_1^\infty)\frac{1}{z}.
\end{equation*}
Plugging this last relation in \eqref{eq:baisilverstein1} leads then to
\begin{equation}
	\label{eq:baisilverstein}
	m_{\xi_2^{-1}}\left(\frac{z}{1-c_1^\infty-c_1^\infty zm_{\mu^{-1}}(z)}\right)=m_{\mu^{-1}}(z)\left(1-c_1^\infty-c_1^\infty zm_{\mu^{-1}}(z)\right).
\end{equation}

Now, with the convention that, for a probability measure $\theta$, $\theta^{-1}$ is the measure defined through $\theta^{-1}\left([a,b]\right)=\theta\left([\frac1a,\frac1b]\right)$, we have the Stieltjes transform relation
\begin{equation*}
	m_{\theta^{-1}}(z)=-\frac1z-\frac1{z^2}m_{\theta}\left(\frac1z\right).
\end{equation*}
Using this relation in \eqref{eq:baisilverstein1}, we then deduce
\begin{align}
	\nonumber
	zm_{\mu}(z)&=(z+c_1^\infty z^2m_\mu(z))m_{\xi_2}(z+c_1^\infty z^2m_\mu(z)) \\
	\label{eq:link_mmu_mxi2}
	&=\varphi(z)m_{\xi_2}(\varphi(z))
\end{align}
where we recall that $\varphi(z)=z(1+c_1^\infty zm_\mu(z))$.
It will come in handy in the following to differentiate this expression along $z$ to obtain
\begin{align*}
	m_{\xi_2}'(\varphi(z)) &= \frac1{\varphi(z)} \left( \frac{m_\mu(z)+zm'_\mu(z)}{\varphi'(z)} - m_{\xi_2}(\varphi(z)) \right) 
\end{align*}
which might be conveniently rewritten as
\begin{align}
	\label{eq:link_mpmu_mxi2}
	m_{\xi_2}'(\varphi(z)) &= \frac1{\varphi(z)} \left( -\frac{\psi'(z)}{c_2^\infty\varphi'(z)} - m_{\xi_2}(\varphi(z)) \right). 
\end{align}

We next determine $m_{\xi_2}$ as a function of $\nu$. Since $\hat{C}_2$ is itself a sample covariance matrix, we may apply again the results from \citep{SIL95}. Denoting $\tilde{\xi}_2$ the almost sure limiting spectrum distribution of $\frac{1}{n_2}\tilde{X}_2^\trans C\tilde{X}_2$, we first have
\begin{equation}
	\label{eq:baisilverstein2}
	m_{\tilde{\xi}_2}(z)=\left(-z+c_2^\infty\int \frac{tdd\nu(t)}{1+tm_{\tilde{\xi}_2}(z)}\right)^{-1}.
\end{equation}
Similar to previously, we have the Stieltjes transform relation $m_{\tilde{\xi}_2}(z)=c_2^\infty m_{\xi_2}(z)-\frac{(1-c_2^\infty)}{z}$ which yields, when plugged in \eqref{eq:baisilverstein2}
\begin{align}
	\label{eq:link_mnu_mxi2}
	m_\nu\left( -\frac{z}{c_2^\infty zm_{\xi_2}(z)-(1-c_2^\infty)} \right) &= -m_{\xi_2}(z) \left(c_2^\infty zm_{\xi_2}(z)-(1-c_2^\infty)\right).
\end{align}

The two relations \eqref{eq:link_mmu_mxi2} and \eqref{eq:baisilverstein2} will be instrumental to relating $\int fd\nu$ to the observation measure $\mu_p$, as described in the next section.

\begin{remark}[The case $c_2>1$]
	The aforementioned reasoning carries over to the case $c_2>1$. Indeed, since the equation \eqref{eq:baisilverstein1} is now meaningless (as the support of $\xi_2$ contains the atom $\{0\}$), consider the model $\hat{C}_1^{-1}(\hat C_2+\varepsilon I_p)=\hat{C}_1^{-1}\hat C_2+\varepsilon \hat{C}_1^{-1}$ for some small $\varepsilon>0$. Then \eqref{eq:link_mmu_mxi2} holds with now $\xi_2$ the limiting empirical spectral distribution of $\hat C_2+\varepsilon I_p$. Due to $\varepsilon$, Equation~\eqref{eq:baisilverstein2} now holds with $m_{\tilde \xi_2}(z)$ replaced by $m_{\tilde \xi_2}(z+\varepsilon)$. By continuity in the small $\varepsilon$ limit, we then have that \eqref{eq:link_mmu_mxi2} and \eqref{eq:link_mnu_mxi2} still hold in the small $\varepsilon$ limit. Now, since $\hat{C}_1^{-1}(\hat C_2+\varepsilon I_p)-\hat{C}_1^{-1}\hat C_2=\varepsilon \hat{C}_1^{-1}$, the operator norm of which amost surely vanishes as $\varepsilon\to 0$ (as per the almost sure boundedness of $\limsup_p \|\hat{C}_1^{-1}\|$), we deduce that $\mu_p\to \mu$ defined through \eqref{eq:link_mmu_mxi2} and \eqref{eq:link_mnu_mxi2}, almost surely, also for $c_2^\infty>1$.
\end{remark}

\subsection{Integral formulation over $m_\nu$}

With the formulas above, we are now in position to derive the proposed estimator. We start by using Cauchy's integral formula to obtain
\begin{align*}
	\int f d\nu &= -\frac1{2\pi\imath} \oint_{\Gamma_\nu} f(z)m_\nu(z)dz
\end{align*}
for $\Gamma_\nu$ a complex contour surrounding the support of $\nu$ but containing no singularity of $f$ in its inside. This contour is carefully chosen as the image of the mapping $\omega\mapsto z=-\omega/(c_2^\infty\omega m_{\xi_2}(\omega)-(1-c_2^\infty))$ of another contour $\Gamma_{\xi_2}$ surrounding the limiting support of $\xi_2$; the details of this (non-trivial) contour change are provided in Appendix~\ref{sec:contour} (where it is seen that the assumption $c_2^\infty<1$ is crucially exploited). We shall admit here that this change of variable is licit.

Operating the aforementioned change of variable gives
\begin{align}
	\label{eq:integral_xi2}
	\int f d\nu &=\frac1{2\pi\imath} \oint_{\Gamma_{\xi_2}} \frac{f\left( \frac{-\omega}{c_2^\infty \omega m_{\xi_2}(\omega)-(1-c_2^\infty )} \right) m_{\xi_2}(\omega)\left( c_2^\infty \omega^2m'_{\xi_2}(\omega)+(1-c_2^\infty ) \right)}{c_2^\infty \omega m_{\xi_2}(\omega)-(1-c_2^\infty )} d\omega
\end{align}
where we used \eqref{eq:link_mnu_mxi2} to eliminate $m_\nu$.

To now eliminate $m_{\xi_2}$ and obtain an integral form only as a function of $m_\mu$, we next proceed to the variable change $u\mapsto \omega = \varphi(u)= u+c_1^\infty u^2m_\mu(u)$. Again, this involves a change of contour, which is valid as long as $\Gamma_{\xi_2}$ is the image by $\varphi$ of a contour $\Gamma_\mu$ surrounding the support of $\mu$, which is only possible if $c_1^\infty <1$ (see Appendix~\ref{sec:contour} for further details). With this variable change, we can now exploit the relations \eqref{eq:link_mmu_mxi2} and \eqref{eq:link_mpmu_mxi2} to obtain, after basic algebraic calculus (using in particular the relation $um_\mu(u)=(-\psi(u)+1-c_2^\infty )/c_2^\infty $)
\begin{align*}
	\int f d\nu &= \frac{-1}{2\pi\imath} \oint_{\Gamma_\mu} f\left( \frac{\varphi(u)}{\psi(u)} \right) \frac1{c_2^\infty \varphi(u)} \left[\varphi(u)\psi'(u)-\psi(u)\varphi'(u) \right]du \\
		    &= \frac1{2\pi\imath} \oint_{\Gamma_\mu} f\left( \frac{\varphi(u)}{\psi(u)} \right) \frac{\psi(u)}{c_2^\infty } \left[ \frac{\varphi'(u)}{\varphi(u)} - \frac{\psi'(u)}{\psi(u)} \right]du.
\end{align*}
It then remains to use the convergence $\nu_p\to\nu$ and $m_{\mu_p}\asto m_\mu$, along with the fact that the eigenvalues $\hat C_1^{-1}\hat C_2$ almost surely do not escape the limiting support $\mu$ as $p\to\infty$ (this is ensured from \citep{SIL98}, Item~\ref{ass:lim_eigs_Ca} of Assumption~\ref{ass:growth} and the analyticity of the involved functions), to retrieve Theorem~\ref{th:contour_integral} by uniform convergence on the compact contour (see also \citep{COU10b} for a similar detailed derivation).

\begin{remark}[Case $C_1$ known]
	The case where $C_1$ is known is equivalent to setting $c_1^\infty \to 0$ above, leading in particular to $m_\mu=m_{\xi_2}$ and to the unique functional equation
	\begin{align*}
		m_\nu\left( \frac{z}{1-c_2^\infty-c_2^\infty zm_\mu(z)} \right) &= m_\nu(z) \left( 1-c_2^\infty-c_2^\infty zm_\nu(z) \right).
	\end{align*}
	In particular, if $C_1=C_2$, this reduces to
	\begin{align*}
		1 &= -m_\mu(z) (z-\psi(z))
	\end{align*}
	with $\psi(z)=1-c_2^\infty -c_2^\infty zm_\mu(z)$, which is the functional Stieltjes-tranform equation of the popular Mar\u{c}enko--Pastur law \citep{MAR67}. 
\end{remark}

\section{Integral Calculus}
\label{sec:integral_calculus}

To compute the complex integral, note first that, depending on $f$, several types of singularities in the integral may arise. Of utmost interest (but not always exhaustively, as we shall see for $f(t)=\log(1+st)$) are: (i) the eigenvalues $\lambda_i$ of $\hat{C}_1^{-1}\hat{C}_2$, (ii) the values $\eta_i$ such that $\varphi_p(\eta_i)=0$, (iii) the values $\zeta_i$ such that $\psi_p(\zeta_i)=0$.

In the following, we first introduce a sequence of intermediary results of interest for most of the integral calculi.

\subsection{Rational expansion}

At the core of the subsequent analysis is the function $\left(\frac{\varphi_p'(z)}{\varphi_p(z)}-\frac{\psi_p'(z)}{\psi_p(z)}\right)\frac{\psi_p(z)}{c_2}$. As this is a mere rational function, we first obtain the following important expansion, that will be repeatedly used in the sequel:
\begin{align}
	\label{eq:expansion}
	&\left(\frac{\varphi_p'(z)}{\varphi_p(z)}-\frac{\psi_p'(z)}{\psi_p(z)} \right)\frac{\psi_p(z)}{c_2} \\ &= \left( \frac1p - \frac{c_1+c_2-c_1c_2}{c_1c_2} \right) \sum_{j=1}^p \frac1{z-\lambda_j} + \frac{1-c_2}{c_2} \frac1z + \frac{c_1+c_2-c_1c_2}{c_1c_2}\sum_{j=1}^p \frac1{z-\eta_j}.\nonumber
\end{align}
This form is obtained by first observing that the $\lambda_j$'s, $\eta_j$'s and $0$ are the poles of the left-hand side expression. Then, pre-multiplying the left-hand side by $(z-\lambda_j)$, $z$, or $(z-\eta_j)$ and taking the limit when these terms vanish, we recover the right-hand side, using in particular the following estimates (which easily entail from the definitions of $\varphi_p$ and $\psi_p$):
\begin{align*}
	\varphi_p(z) &= \frac{c_1}p \frac{\lambda_i^2}{\lambda_i-z} - 2c_1\frac{\lambda_i}p + \lambda_i + \frac{c_1}p \sum_{j\neq i} \frac{\lambda_i^2}{\lambda_j-\lambda_i} + O(\lambda_i-z ) \\
	\varphi_p'(z) &= \frac{c_1}p \frac{\lambda_i^2}{(\lambda_i-z)^2} + O(1) \\
	\psi_p(z) &= -\frac{c_2}p \frac{\lambda_i}{\lambda_i-z} + \frac{c_2}p + 1 -c_2 - \frac{c_2}p \sum_{j\neq i} \frac{\lambda_i}{\lambda_j-\lambda_i} + O(\lambda_i-z) \\
	\psi_p'(z) &= - \frac{c_2}p \frac{\lambda_i}{(\lambda_i-z)^2} + O(1)
\end{align*}
in the vicinity of $\lambda_i$, along with $\psi_p(\eta_i)=\frac{c_1+c_2-c_1c_2}{c_1}$ and $\psi_p(0)=1-c_2$.

\bigskip

From this expression, we have the following immediate corollary.
\begin{remark}[Residue for $f$ analytic at $\lambda_i$]
	\label{rem:Res_lambda}
	If $f \circ (\varphi_p/\psi_p)$ is analytic in a neighborhood of $\lambda_i$, i.e., if $f$ is analytic in a neighborhood of $-(c_1/c_2)\lambda_i$, then $\lambda_i$ is a first order pole for the integrand, leading to the residue
\begin{align*}
	{\rm Res}(\lambda_i) &= - f\left( -\frac{c_1}{c_2}\lambda_i \right) \left[\frac{c_1+c_2-c_1c_2}{c_1c_2} - \frac1p\right].
\end{align*}
\end{remark}

\subsection{Characterization of $\eta_i$ and $\zeta_i$, and $\varphi_p/\psi_p$}

First note that the $\eta_i$ (the zeros of $\varphi_p(z)$) and $\zeta_i$ (the zeros of $\psi_p(z)$) are all real as one can verify that, for $\Im[z]\neq 0$, $\Im[\varphi_p(z)]\Im[z]>0$ and $\Im[\psi_p(z)]\Im[z]<0$.

\begin{figure}
	\includegraphics[width=\linewidth]{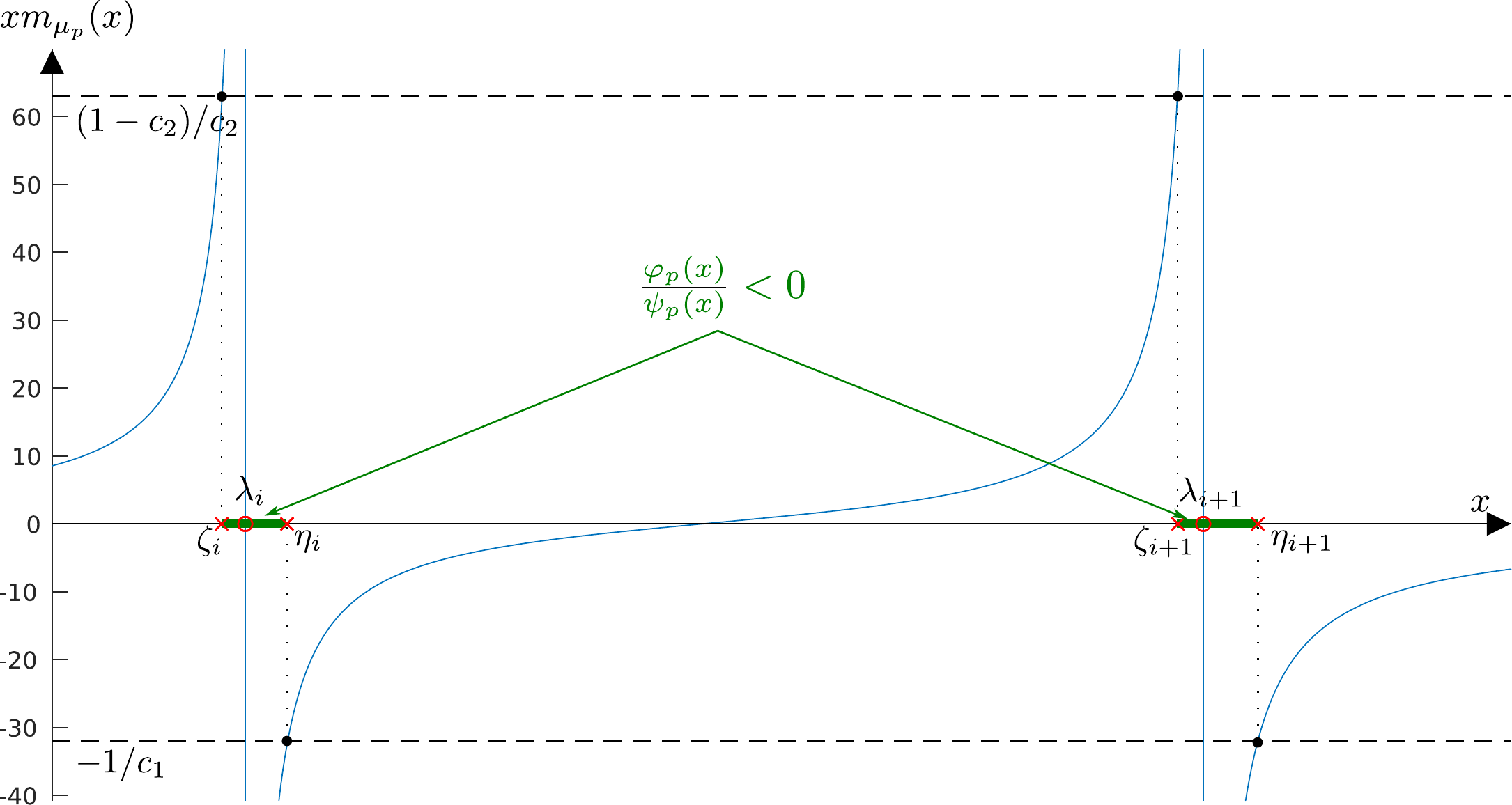}
	\caption{Visual representation of $x\mapsto M_p(x)=xm_{\mu_p}(x)$; here for $p=4$, $n_1=8$, $n_2=16$. Solutions to $M_p(x)=-1/c_1$ (i.e., $\eta_i$'s) and to $M_p(x)=(1-c_2)/c_2$ (i.e., $\zeta_i$'s) indicated in red crosses. Green solid lines indicate sets of negative $\varphi_p/\psi_p$.}
	\label{fig:M_p}
\end{figure}

Before establishing the properties of $\varphi_p$ and $\psi_p$ in the vicinity of $\eta_i$ and $\zeta_i$, let us first locate these values. A study of the function $M_p:\RR\to \RR$, $x\mapsto xm_{\mu_p}(x)$ (see Figure~\ref{fig:M_p}) reveals that $M_p$ is increasing (since $x/(\lambda_i-x)=-1+1/(\lambda_i-x)$) and has asymptotes at each $\lambda_i$ with $\lim_{x\uparrow \lambda_i} M_p(x)=\infty$ and $\lim_{x\downarrow \lambda_i} M_p(x)=-\infty$. As a consequence, since $\varphi_p(x)=0\Leftrightarrow M_p(x)=-\frac1{c_1}<-1$, there exists exactly one solution to $\varphi_p(x)=0$ in the set $(\lambda_i,\lambda_{i+1})$. This solution will be subsequently called $\eta_i$. Since $M_p(x)\to -1$ as $x\to\infty$, there exists a last solution to $\varphi_p(x)=0$ in $(\lambda_p,\infty)$, hereafter referred to as $\eta_p$. Similarly, $\psi_p(x)=0\Leftrightarrow M_p(x)=(1-c_2)/c_2>0$ and thus there exists exactly one solution, called $\zeta_i$ in $(\lambda_{i-1},\lambda_i)$. When $x\to 0$, $M_p(x)\to 0$ so that a further solution is found in $(0,\lambda_1)$, called $\zeta_1$. Besides, due to the asymptotes at every $\lambda_i$, we have that $\zeta_1<\lambda_i<\eta_1<\zeta_2<\ldots<\eta_p$.

As such, the set $\Gamma$ defined in Theorem~\ref{th:contour_integral} exactly encloses all $\eta_i$, $\lambda_i$, and $\zeta_i$, for $i=1,\ldots,p$, possibly to the exception of the leftmost $\zeta_1$ and the rightmost $\eta_p$ (as those are not comprised in a set of the form $[\lambda_{i+1},\lambda_i]$). To ensure that the latter do asymptotically fall within the interior of $\Gamma$, one approach is to exploit Theorem~\ref{th:contour_integral} for the elementary function $f(t)=1$. There we find that
\begin{align*}
	\frac1{2\pi\imath}\oint_{\Gamma_\nu} m_\nu(z)dz - \frac1{2\pi\imath}\oint_{\Gamma} \left( \frac{\varphi'_p(z)}{\varphi_p(z)} - \frac{\psi'_p(z)}{\psi_p(z)} \right) \frac{\psi_p(z)}{c_2}dz \asto 0.
\end{align*}
The left integral is easily evaluated by residue calculus and equals $-1$ (each $\lambda_i(C_1^{-1}C_2)$, $1\leq i\leq p$, is a pole with associated residue $-1/p$), while the right integral can be computed from \eqref{eq:expansion} again by residue calculus and equals $-1+\frac{c_1+c_2-c_1c_2}{c_1c_2}(p-\#\{\eta_i \in \Gamma^\circ\})$ with $\Gamma^\circ$ the ``interior'' of $\Gamma$. As such, since both integrals are (almost surely) arbitrarily close in the large $p$ limit, we deduce that $\#\{\eta_i \in \Gamma^\circ\}=p$ for all large $p$ and thus, in particular, $\eta_p$ is found in the interior of $\Gamma$. To obtain the same result for $\zeta_1$, note that, from the relation $\psi_p(z)=\frac{c_1+c_2-c_1c_2}{c_1}-\frac{c_2}{c_1}\frac{\varphi_p(z)}z$ along with the fact that $\frac{\varphi'_p(z)}{\varphi_p(z)} - \frac{\psi'_p(z)}{\psi_p(z)}$ is an exact derivative (of $\log(\varphi_p/\psi_p)$), the aforementioned convergence can be equivalently written
\begin{align*}
	\frac1{2\pi\imath}\oint_{\Gamma_\nu} m_\nu(z)dz - \frac{-1}{2\pi\imath}\oint_{\Gamma} \left( \frac{\varphi'_p(z)}{\varphi_p(z)} - \frac{\psi'_p(z)}{\psi_p(z)} \right) \frac{\varphi_p(z)}{zc_1}dz \asto 0.
\end{align*}
Reproducing the same line of argument (with an expansion of $( \frac{\varphi'_p(z)}{\varphi_p(z)} - \frac{\psi'_p(z)}{\psi_p(z)}) \frac{\varphi_p(z)}{zc_1}$ equivalent to \eqref{eq:expansion}), the same conclusion arises and we then proved that both $\zeta_1$ and $\eta_p$ (along with all other $\zeta_i$'s and $\eta_i$'s) are asymptotically found within the interior of $\Gamma$.

\medskip

One can also establish that, on its restriction to $\RR^+$, $\varphi_p$ is everywhere positive but on the set $\cup_{i=1}^p(\lambda_i,\eta_i)$. Similarly, $\psi_p$ is everywhere positive but on the set $\cup_{i=1}^p(\zeta_i,\lambda_i)$. As a consequence, the ratio $\varphi_p/\psi_p$ is everywhere positive on $\RR^+$ but on the set $\cup_{i=1}^p(\zeta_i,\eta_i)$.

\begin{figure}
	\centering
	\includegraphics[width=.6\linewidth]{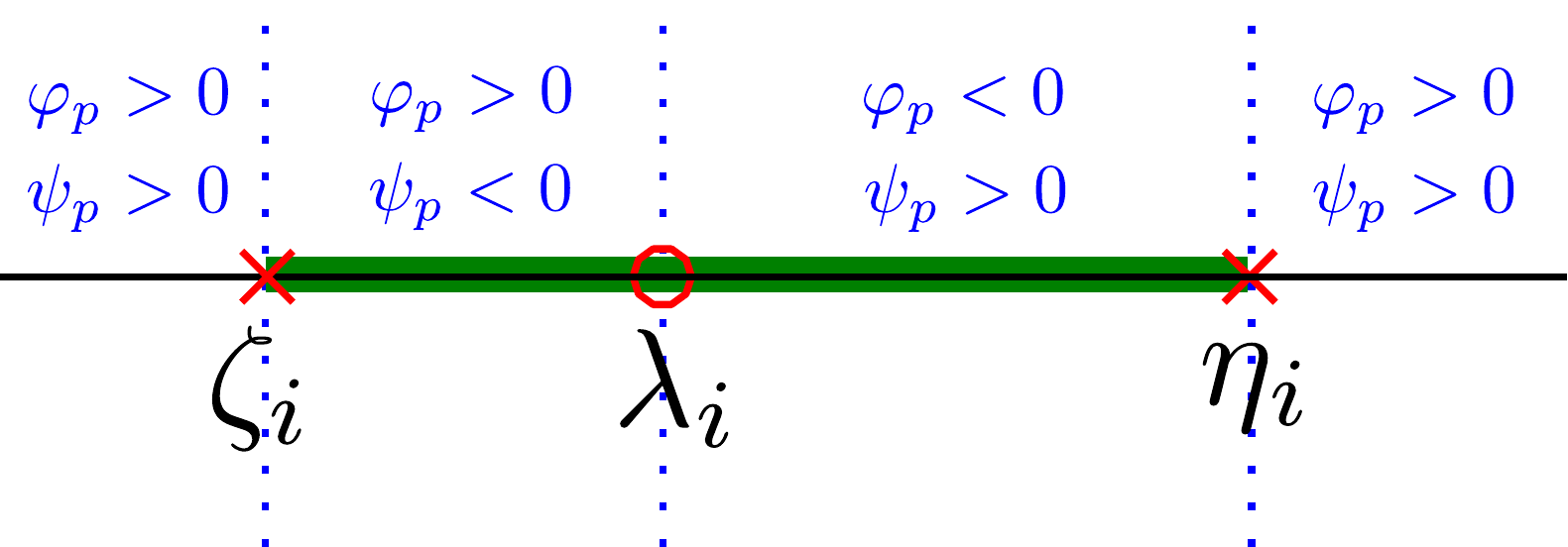}
	\caption{Visual representation of the signs of $\varphi_p$ and $\psi_p$ around singularities.}
	\label{fig:sign_phi_over_psi}
\end{figure}

These observations are synthesized in Figure~\ref{fig:sign_phi_over_psi}.

\bigskip

In terms of monotonicity on their restrictions to the real axis, since $\psi_p(x)=1-c_2\int \frac{t}{t-x}d\mu_p(t)$, $\psi_p$ is decreasing. As for $\varphi_p$, note that
\begin{align*}
	\varphi_p'(x) &= 1 + 2c_1\int \frac{x}{t-x}d\mu_p(t) + c_1\int \frac{x^2}{(t-x)^2}d\mu_p(t) \\ &= \int \frac{t^2-2(1-c_1)xt+(1-c_1)x^2}{(t-x)^2}d\mu_p(t).
\end{align*}
Since $c_1<1$, we have $1-c_1>(1-c_1)^2$, and therefore 
\begin{align*}
	\varphi_p'(x) &> \int \frac{ (t-(1-c_1)x)^2}{(t-x)^2}d\mu_p(t) >0
\end{align*}
ensuring that $\varphi_p$ is increasing on its restriction to $\RR$.

Showing that $x\mapsto\varphi_p(x)/\psi_p(x)$ is increasing is important for the study of the case $f(t)=\log(1+st)$ but is far less immediate. This unfolds from the following remark, also of key importance in the following.

\begin{remark}[Alternative form of $\varphi_p$ and $\psi_p$]
	\label{rem:prod_form}
	It is interesting to note that, in addition to the zero found at $z=0$ for $\varphi_p$, we have enumerated all zeros and poles of the rational functions $\varphi_p$ and $\psi_p$ (this can be ensured from their definition as rational functions) and it thus comes that
\begin{align}
	\label{eq:phi_in_product}
	\varphi_p(z) &= (1-c_1) z\frac{\prod_{j=1}^p(z-\eta_j)}{\prod_{j=1}^p (z-\lambda_j)} \\
	\label{eq:psi_in_product}
	\psi_p(z)    &= \frac{\prod_{j=1}^p(z-\zeta_j)}{\prod_{j=1}^p (z-\lambda_j)}
\end{align}
where the constants $1-c_1$ and $1$ are found by observing that, as $z=x\in\RR\to\infty$, $\varphi_p(x)/x\to 1-c_1$ and $\psi_p(x)\to 1$. In particular
\begin{align}
	\label{eq:phi_over_psi_in_product}
	\frac{\varphi_p(z)}{\psi_p(z)} &= (1-c_1)z \frac{\prod_{j=1}^p(z-\eta_j)}{\prod_{j=1}^p(z-\zeta_j)}.
\end{align}
\end{remark}

A further useful observation is that the $\eta_i$'s are the eigenvalues of
\begin{align*}
	\Lambda - \frac1{p-n_1}{\sqrt{\lambda}\sqrt{\lambda}^\trans}
\end{align*}
where $\Lambda={\rm diag}(\{\lambda_i\}_{i=1}^p)$ and $\lambda=(\lambda_1,\ldots,\lambda_p)^\trans$. Indeed, these eigenvalues are found by solving 
\begin{align*}
	0&=\det\left(\Lambda - \frac{\sqrt{\lambda}\sqrt{\lambda}^\trans}{p-n_1}-xI_p\right) \\
	&=\det(\Lambda-xI_p)\det\left(I_p-(\Lambda-xI_p)^{-1}\frac{\sqrt{\lambda}\sqrt{\lambda}^\trans}{p-n_1}\right) \\
	&= \det(\Lambda-xI_p)\left(1-\frac1{p-n_1}\sqrt{\lambda}^\trans(\Lambda-xI_p)^{-1}\sqrt{\lambda}\right) \\
	&= \det(\Lambda-xI_p)\left(1-\frac1{p-n_1}\sum_{i=1}^p \frac{\lambda_i}{\lambda_i-x}\right)
\end{align*}
which, for $x$ away from the $\lambda_i$ (not a solution to $\varphi_p(x)=0$), reduces to $\frac1p\sum_{i=1}^p \frac{\lambda_i}{\lambda_i-x}=1-\frac1{c_1}$, which is exactly equivalent to $m_{\mu_p}(x)=-\frac1{c_1x}$, i.e., $\varphi_p(x)=0$.

Similarly, the $\zeta_i$'s are the eigenvalues of the matrix
\begin{align*}
	\Lambda - \frac1{n_2}{\sqrt{\lambda}\sqrt{\lambda}^\trans}.
\end{align*} 

\bigskip

These observations allow for the following useful characterization of $\varphi_p/\psi_p$:
\begin{align*}
	\frac{\varphi_p(z)}{\psi_p(z)} &= (1-c_1)z \frac{\det\left(zI_p-\Lambda-\frac1{n_1-p}\sqrt{\lambda}\sqrt{\lambda}^\trans\right)}{\det\left(zI_p-\Lambda+\frac1{n_2}\sqrt{\lambda}\sqrt{\lambda}^\trans\right)} \\
	&= (1-c_1)z \left( 1 - \frac{n_1+n_2-p}{n_2(n_1-p)}\sqrt{\lambda}^\trans \left(zI_p-\Lambda+\frac1{n_2}\sqrt{\lambda}\sqrt{\lambda}^\trans\right)^{-1}\sqrt{\lambda} \right)
\end{align*}
(after factoring out the matrix in denominator from the determinant in the numerator) the derivative of which is, after simplification,
\begin{align*}
	\left(\frac{\varphi_p(z)}{\psi_p(z)}\right)' &= (1-c_1)\left( 1 + \frac{n_1+n_2-p}{n_2(n_1-p)}\sqrt{\lambda}^\trans Q \left(\Lambda-\frac1{n_2}{\sqrt{\lambda}\sqrt{\lambda}^\trans}\right)Q \sqrt{\lambda} \right).
\end{align*}
for $Q=(zI_p-\Lambda+\frac1{n_2}{\sqrt{\lambda}\sqrt{\lambda}^\trans})^{-1}$.
Since $\Lambda-\frac1{n_2}\sqrt{\lambda}\sqrt{\lambda}^\trans$ is positive definite (its eigenvalues being the $\zeta_i$'s), on the real axis the derivative is greater than $1-c_1>0$ and the function $x\mapsto \varphi_p(x)/\psi_p(x)$ is therefore increasing. 

\begin{figure}
	\centering
	\includegraphics[width=\linewidth]{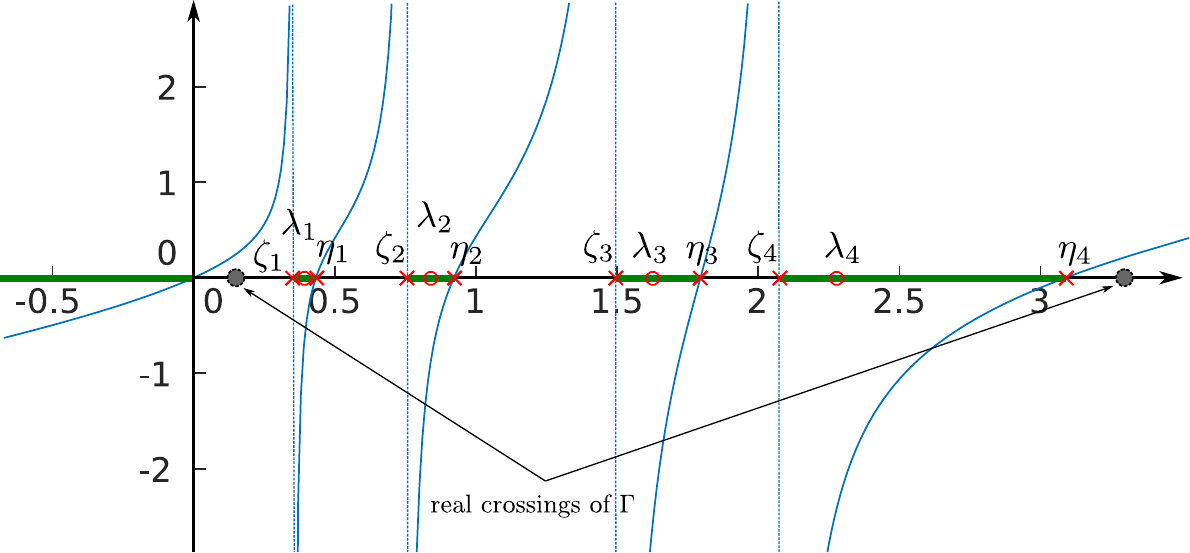}
	\caption{Example of visual representation of $\varphi_p/\psi_p:\RR\to\RR$, $x\mapsto \varphi_p(x)/\psi_p(x)$; here for $p=4$, $n_1=8$, $n_2=16$. In green solid lines are stressed the sets over which $\varphi_p(x)/\psi_p(x)<0$ (which correspond to branch cuts in the study of $f(z)=\log^k(z)$). Possible real crossings of the contour $\Gamma$ are indicated, notably showing that no branch cut is passed through when $f(z)=\log^k(z)$.}
	\label{fig:phi_over_psi}
\end{figure}

Figure~\ref{fig:phi_over_psi} displays the behavior of $\varphi_p/\psi_p$ when restricted to the real axis.

\bigskip

Since we now know that the contour $\Gamma$ from Theorem~\ref{th:contour_integral} encloses exactly all $\eta_i$'s and $\zeta_i$'s, it is sensible to evaluate the residues for these values when $f(z)$ is analytic in their neighborhood.
\begin{remark}[Residue for $f$ analytic at $\eta_i$ and $\zeta_i$]
	\label{rem:Res_eta_zeta}
If $f$ is analytic with no singularity at zero, then the integral has a residue at $\eta_i$ easily found to be
\begin{align*}
	{\rm Res}(\eta_i) &= f(0) \frac{c_1+c_2-c_1c_2}{c_2}.
\end{align*}
Similarly, if $f(\omega)$ has a well defined limit as $|\omega|\to \infty$, then no residue is found at $\zeta_i$. 
\end{remark}

As a consequence of Remarks~\ref{rem:Res_lambda} and \ref{rem:Res_eta_zeta}, we have the following immediate corollary.
\begin{remark}[The case $f(t)=t$]
	\label{rem:f=t}
In the case where $f(t)=t$, a singularity appears at $\zeta_i$, which is nonetheless easily treated by noticing that the integrand then reduces to
\begin{align*}
	f\left( \frac{\varphi_p(z)}{\psi_p(z)} \right) \left( \frac{\varphi_p'(z)}{\varphi_p(z)} - \frac{\psi_p'(z)}{\psi_p(z)}\right)\frac{\psi_p(z)}{c_2} = \frac{\varphi_p'(z)}{c_2} - \frac{\psi_p'(z)\varphi_p(z)}{c_2\psi_p(z)}
\end{align*}
and thus, with $\psi_p(z)=(z-\zeta_i)\psi'_p(\zeta_i)+O((z-\zeta_i)^2)$, we easily find the residue
\begin{align*}
	{\rm Res}_{\{f(t)=t\}}(\zeta_i) &= -\frac1{c_2\varphi_p(\zeta_i)} = -\zeta_i\frac{c_1+c_2-c_1c_2}{c_2^2}.
\end{align*}
Together with Remarks~\ref{rem:Res_lambda} and \ref{rem:Res_eta_zeta}, along with the fact that $\Gamma$ encloses all $\eta_i$ and $\lambda_i$, for $i=1,\ldots,p$, we then find that
\begin{align*}
	\int t d\nu_p(t) - \left[ \frac{c_1+c_2-c_1c_2}{c_2^2} \sum_{i=1}^p (\lambda_i - \zeta_i) - \frac{c_1}{c_2} \frac1p\sum_{i=1}^p \lambda_i \right] \asto 0.
\end{align*}
By then noticing that $\sum_i \zeta_i=\tr (\Lambda-\frac1{n_2}\sqrt{\lambda}\sqrt{\lambda}^\trans)=(1-c_2/p)\sum_i \lambda_i $, we retrieve Corollary~\ref{cor:t}.
\end{remark}

\subsection{Development for \lowercase{$f(t)=\log(t)$}}

The case $f(t)=\log(t)$ leads to an immediate simplification as, then, $\log\det( \hat{C}_1^{-1}\hat{C}_2)=\log\det(\hat{C}_2)-\log\det(\hat{C}_1)$; one may then use previously established results from the random matrix literature (e.g., the G-estimators in \citep{GIR87} or more recently \citep{KAM11}) to obtain the sought-for estimate. Nonetheless, the full explicit derivation of the contour integral in this case is quite instructive and, being simpler than the subsequent cases where $f(t)=\log^2(t)$ or $f(t)=\log(1+st)$ that rely on the same key ingredients, we shall here conduct a thorough complex integral calculus.

\medskip

\begin{figure}
	\centering
	\includegraphics[width=\linewidth]{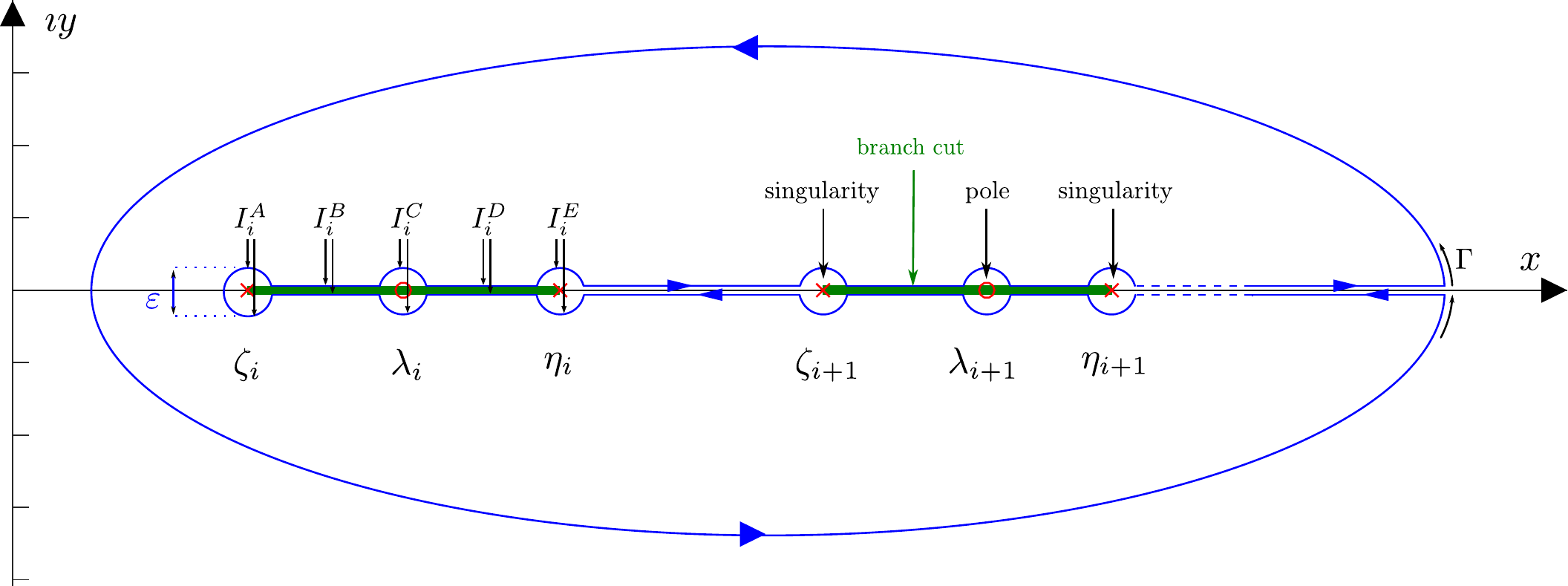}
	\caption{Chosen integration contour. The set $I_i^B$ is the disjoint union of the segments $[\zeta_i+\varepsilon+0^+\imath,\lambda_i-\varepsilon+0^+\imath]$ and $[\zeta_i+\varepsilon+0^-\imath,\lambda_i-\varepsilon+0^-\imath]$. Similarly the set $I_i^D$ is the disjoint union of the segments $[\lambda_i+\varepsilon+0^+\imath,\eta_i-\varepsilon+0^+\imath]$ and $[\lambda_i+\varepsilon+0^-\imath,\eta_i-\varepsilon+0^-\imath]$. The sets $I_i^A$, $I_i^C$ and $I_i^E$ are the disjoint unions of semi-circles (in the upper- or lower-half complex plane) of diameters $\varepsilon$ surrounding $\zeta_i$, $\lambda_i$ and $\eta_i$ respectively.}
	\label{fig:contour}
\end{figure}

For $z\in\CC$, define first $f(z)=\log(z)$ where $\log(z)=\log(|z|)e^{i\arg(z)}$, with $\arg(z)\in(-\pi,\pi]$. For this definition of the complex argument, since $\varphi_p(x)/\psi_p(x)$ is everywhere positive but on $\cup_{i=1}^p(\zeta_i,\eta_i)$, we conclude that $\arg(\varphi_p(z)/\psi_p(z))$ abruptly moves from $\pi$ to $-\pi$ as $z$ moves from $x+0^+\imath$ to $x+0^-\imath$ for all $x\in \cup_{i=1}^p (\zeta_i,\eta_i)$. This creates a set of $p$ branch cuts $[\zeta_i,\eta_i]$, $i=1,\ldots,p$ as displayed in Figure~\ref{fig:phi_over_psi}. This naturally leads to computing the complex integral estimate of $\int fd\nu$ based on the contour displayed in Figure~\ref{fig:contour}, which avoids the branch cuts.

This contour encloses no singularity of the integrand and therefore has a null integral. With the notations of Figure~\ref{fig:contour}, the sought-for integral (over $\Gamma$) therefore satisfies
\begin{align*}
	0 &= \oint_{\Gamma} + \sum_{i=1}^p \left( \int_{I^A_i} + \int_{I^B_i} + \int_{I^C_i} + \int_{I^D_i} + \int_{I^E_i} \right).
\end{align*}	

We start by the evaluation of the integrals over $I^B_i$ and $I^D_i$, which can be similarly handled. To this end, note that, since $\arg(\frac{\varphi_p}{\psi_p})$ moves from $\pi$ to $-\pi$ across the branch cut, we have
\begin{align*}
	\frac{1}{2\pi \imath} \int_{I^B_i} &= \frac{1}{2\pi \imath}\int_{\zeta_i+\varepsilon}^{\lambda_i-\varepsilon} \left[\log \left( - \frac{\varphi_p(x)}{\psi_p(x)}\right)+\imath \pi -\log\left(-\frac{\varphi_p(x)}{\psi_p(x)}\right) +\imath\pi\right] \nonumber \\
	&\times\left( \frac{\varphi_p'(x)}{\varphi_p(x)}-\frac{\psi_p'(x)}{\psi_p(x)}\right)\frac{\psi_p(x)}{c_{2}}dx\\
	&=\int_{\zeta_i+\varepsilon}^{\lambda_i-\varepsilon} \left( \frac{\varphi_p'(x)}{\varphi_p(x)}-\frac{\psi_p'(x)}{\psi_p(x)}\right)\frac{\psi_p(x)}{c_{2}}dx.
\end{align*}


We first exploit the rational form expansion \eqref{eq:expansion} of $(\frac{\varphi_p'(z)}{\varphi_p(z)}-\frac{\psi_p'(z)}{\psi_p(z)})\frac{\psi_p(x)}{c_{2}}$ to obtain the integral over $I_i^B$
\begin{align*}
	\frac{1}{2\pi \imath} \int_{I^B_i} &= \left(\frac1p - \frac{c_1+c_2-c_1c_2}{c_1c_2} \right) \left( \sum_{j\neq i} \log\left|\frac{\lambda_i-\lambda_j}{\zeta_i-\lambda_j} \right| + \log \left| \frac{\varepsilon}{\zeta_i-\lambda_i} \right| \right) \nonumber \\
	&+ \frac{1-c_2}{c_2}\log\frac{\lambda_i}{\zeta_i} + \frac{c_1+c_2-c_1c_2}{c_1c_2}\sum_{j=1}^p\log\left| \frac{\lambda_i-\eta_j}{\zeta_i-\eta_j}\right| + o(\varepsilon).
\end{align*}

The treatment is similar for the integral over $I_i^D$ which results, after summation of both integrals, to
\begin{align*}
	\frac{1}{2\pi \imath} \int_{I^B_i\cup I_i^D} &= \left(\frac1p - \frac{c_1+c_2-c_1c_2}{c_1c_2} \right) \sum_{j=1}^p \log\left|\frac{\eta_i-\lambda_j}{\zeta_i-\lambda_j} \right| + \frac{1-c_2}{c_2}\log\frac{\eta_i}{\zeta_i} \nonumber \\
	&+ \frac{c_1+c_2-c_1c_2}{c_1c_2}\left(\sum_{j\neq i} \log\left| \frac{\eta_i-\eta_j}{\zeta_i-\eta_j}\right|+\log \left|\frac{\varepsilon}{\zeta_i-\eta_i} \right|\right) + o(\varepsilon).
\end{align*}
Note here the asymmetry in the behavior of the integrand in the neighborhood of $\zeta_i$ ($+\varepsilon$) and $\eta_i$ ($-\varepsilon$); in the former edge, the integral is well defined while in the latter it diverges as $\log\varepsilon$ which must then be maintained.

Summing now over $i\in\{1,\ldots,p\}$, we recognize a series of identities. In particular, note that from the product form \eqref{eq:phi_over_psi_in_product},
\begin{align*}
	\sum_{j=1}^p \sum_{i=1}^p \log\left|\frac{\eta_i-\lambda_j}{\zeta_i-\lambda_j} \right| &= \sum_{j=1}^p \log\left| \frac{\frac{\psi_p}{\varphi_p}(\lambda_j)}{(1-c_1)\lambda_j} \right| \\
	&= \log\left( \frac{c_1}{c_2(1-c_1)} \right) \\
	\sum_{j=1}^p \log \frac{\eta_i}{\zeta_i} &= \lim_{z\to 0} \log \left( \frac{\frac{\psi_p}{\varphi_p}(z)}{(1-c_1)z} \right) \\ 
	&= -\log\left( (1-c_1)(1-c_2) \right) \\
	\sum_{i=1}^p \sum_{j\neq i} \log\left| \frac{\eta_i-\eta_j}{\zeta_i-\eta_j} \right| + \sum_{i=1}^p \log\left| \frac1{\zeta_i-\eta_i}\right| &= \lim_{z\to \eta_i}\sum_{j=1}^p \log \left| \frac{\frac{\psi_p}{\varphi_p}(z)}{(1-c_1)z(z-\eta_j)} \right| \\ &= \sum_{j=1}^p \log\left| \frac{\left(\frac{\psi_p}{\varphi_p}\right)'(\eta_j)}{(1-c_1)\eta_j} \right|.
\end{align*}

As such, we now find that
\begin{align*}
	\frac{1}{2\pi \imath} \sum_{i=1}^p \int_{I^B_i\cup I_i^D} &= \log\left( \frac{c_1}{c_2(1-c_1)} \right)-\frac{1-c_2}{c_2}\log\left((1-c_1)(1-c_2)\right) \\
	&- \frac{c_1+c_2-c_1c_2}{c_1c_2} \left( p \log \left( \frac{c_1}{c_2(1-c_1)} \right) - \sum_{j=1}^p \log\left| \frac{\left(\frac{\psi_p}{\varphi_p}\right)'(\eta_j)}{(1-c_1)\eta_j} \right|  \right) \\
	&+ \frac{c_1+c_2-c_1c_2}{c_1c_2}p\log\varepsilon+o(\varepsilon).
\end{align*}

The diverging term in $\log\varepsilon$ is compensated by the integral over $I_i^E$. Indeed, letting $z=\eta_{i}+\varepsilon e^{i\theta}$, we may write
\begin{align*}
&\frac{1}{2\pi \imath}\int_{I_i^E} \log\left(\frac{\varphi_p(z)}{\psi_p(z)}\right)\left(\frac{\varphi'_p(z)}{\varphi_p(z)}-\frac{\psi_p'(z)}{\psi_p(z)}\right)\frac{\psi(z)}{c_{2}}dz \\
&=\frac{\varepsilon}{2\pi c_2}\left[\int_{\pi}^{0^+} + \int_{0^-}^{-\pi}\right] \log\left( (1-c_1)(\eta_i+\varepsilon e^{i\theta})\frac{\prod_{j=1}^{p}(\eta_{i}-\eta_{j}+\varepsilon e^{i\theta})}{\prod_{j=1}^{p}(\eta_{i}-\zeta_{j}+\varepsilon e^{i\theta})}\right) \nonumber \\
&\times \left(  \sum_{j=1}^p \frac{\frac1p - \frac{c_1+c_2-c_1c_2}{c_1c_2}}{\eta_i+\varepsilon e^{\imath\theta}-\lambda_j} + \frac{\frac{1-c_2}{c_2} }{\eta_i+\varepsilon e^{\imath\theta}} + \sum_{j=1}^p \frac{\frac{c_1+c_2-c_1c_2}{c_1c_2}}{\eta_i+\varepsilon e^{\imath\theta}-\eta_j} \right) e^{i \theta} d\theta.
\end{align*}
To evaluate the small $\varepsilon$ limit of this term, first remark importantly that, for small $\varepsilon$, the term in the logarithm equals
\begin{align*}
	(1-c_1)\frac{\eta_i}{\eta_i-\zeta_i} \frac{\prod_{j\neq i}(\eta_i-\eta_j)}{\prod_{j\neq i}(\eta_i-\zeta_j)} \varepsilon e^{i\theta} + o(\varepsilon)
\end{align*}
the argument of which equals that of $\theta$. As such, on the integral over $(\pi,0)$, the log term reads $\log|\cdot|+\imath \theta+o(\varepsilon)$, while on $(0,-\pi)$, it reads $\log|\cdot|-\imath \theta+o(\varepsilon)$.
With this in mind, keeping only the non-vanishing terms in the small $\varepsilon$ limit (that is: the term in $\log\varepsilon$ and the term in $\frac1{\varepsilon}$) leads to
\begin{align*}
	&\frac{1}{2\pi \imath}\int_{I_i^E} \log\left(\frac{\varphi_p(z)}{\psi_p(z)}\right)\left(\frac{\varphi'_p(z)}{\varphi_p(z)}-\frac{\psi_p'(z)}{\psi_p(z)}\right)\frac{\psi(z)}{c_{2}}dz \\
	&=\frac{c_1+c_2-c_1c_2}{c_1c_2} \log \varepsilon + \frac{c_1+c_2-c_1c_2}{c_1c_2} \log \left| \left(\frac{\varphi_p}{\psi_p} \right)'(\eta_i) \right| + o(\varepsilon)
\end{align*}
where we used the fact that $\lim_{\varepsilon\to 0} \frac1{\varepsilon e^{\imath\theta}}\left(\frac{\varphi_p}{\psi_p} \right)(\eta_i+\varepsilon e^{\imath\theta})=\left(\frac{\varphi_p}{\psi_p} \right)'(\eta_i)$.

We proceed similarly to handle the integral over $I_i^C$
\begin{align*}
&\frac{1}{2\pi \imath}\int_{I_i^C} \log\left(\frac{\varphi_p(z)}{\psi_p(z)}\right)\left(\frac{\varphi'_p(z)}{\varphi_p(z)}-\frac{\psi_p'(z)}{\psi_p(z)}\right)\frac{\psi(z)}{c_{2}}dz \\
&=\frac{\varepsilon}{2\pi c_2}\left[\int_{\pi}^{0^+} + \int_{0^-}^{-\pi}\right] \log\left( (1-c_1)(\lambda_i+\varepsilon e^{i\theta})\frac{\prod_{j=1}^{p}(\lambda_{i}-\eta_{j}+\varepsilon e^{i\theta})}{\prod_{j=1}^{p}(\lambda_{i}-\zeta_{j}+\varepsilon e^{i\theta})}\right) \nonumber \\
&\times \left(  \sum_{j=1}^p \frac{\frac1p - \frac{c_1+c_2-c_1c_2}{c_1c_2}}{\lambda_i+\varepsilon e^{\imath\theta}-\lambda_j} + \frac{\frac{1-c_2}{c_2} }\lambda_i+\varepsilon e^{\imath\theta} + \sum_{j=1}^p \frac{\frac{c_1+c_2-c_1c_2}{c_1c_2}}{\lambda_i+\varepsilon e^{\imath\theta}-\eta_j} \right) e^{i \theta} d\theta.
\end{align*}
Here, for small $\varepsilon$, the angle of the term in the argument of the logarithm is that of
\begin{align*}
	&\frac{\varphi_p}{\psi_p}(\lambda_i)+\left(\frac{\varphi_p}{\psi_p}\right)'(\lambda_i)\varepsilon e^{\imath \theta}+o(\varepsilon) \\
	&= -\frac{c_1}{c_2}\lambda_i + \varepsilon e^{\imath \theta} \frac{c_1}{c_2}\left(p\frac{c_1+c_2-c_1c_2}{c_1c_2} -1\right) + o(\varepsilon).
\end{align*}
That is, for all large $p$, the argument equals $\pi+o(\varepsilon)<\pi$ uniformly on $\theta\in(0,\pi)$ and $-\pi+o(\varepsilon)>-\pi$ uniformly on $\theta\in(-\pi,0)$; thus the complex logarithm reads $\log|\cdot|+\imath \theta+o(\varepsilon)$ on $(\pi,0)$, while on $(0,-\pi)$, it reads $\log|\cdot|-\imath \theta+o(\varepsilon)$. Proceeding as previously for the integral over $I_i^E$, we then find after calculus that
\begin{align*}
	&\frac{1}{2\pi \imath}\int_{I_i^C} \log\left(\frac{\varphi_p(z)}{\psi_p(z)}\right)\left(\frac{\varphi'_p(z)}{\varphi_p(z)}-\frac{\psi_p'(z)}{\psi_p(z)}\right)\frac{\psi(z)}{c_{2}}dz \\
	&=\left(\frac{c_1+c_2-c_1c_2}{c_1c_2} - \frac1p \right)\log \left( \frac{c_1}{c_2}\lambda_i \right).
\end{align*}
Note that this expression is reminiscent of a ``residue'' at $\lambda_i$ (with negatively oriented contour), according to Remark~\ref{rem:Res_lambda}, however for the function $\log|\cdot|$ and not for the function $\log(\cdot)$, due to the branch cut passing through $\lambda_i$.

The final integral over $I_i^A$ is performed similarly. However, here, it is easily observed that the integral is of order $O(\varepsilon\log \varepsilon)$ in the small $\varepsilon$ limit, and thus vanishes.

Finally, summing up all contributions, we have
\begin{align*}
	\frac1{2\pi\imath}\oint_\Gamma &= - \sum_{i=1}^p \frac1{2\pi\imath} \left( \int_{I_i^A} + \int_{I_i^B} + \int_{I_i^C} + \int_{I_i^D} + \int_{I_i^E} \right) \\
	&= -\log\left( \frac{ c_1(1-c_2)}{c_2}\right) +\frac1{c_2}\log( (1-c_1)(1-c_2) ) + \frac1p\sum_{i=1}^p \log\left( \frac{c_1}{c_2}\lambda_i\right)\\
	&- \frac{c_1+c_2-c_1c_2}{c_1c_2} \left(p\log(1-c_1)+\sum_{i=1}^p\log \lambda_i - \sum_{i=1}^p \log( (1-c_1)\eta_i )\right)  \\
	&=-\log(1-c_2) +\frac1{c_2}\log( (1-c_1)(1-c_2) ) + \frac1p\sum_{i=1}^p\log\lambda_i \\
	&+ \frac{c_1+c_2-c_1c_2}{c_1c_2} \sum_{i=1}^p \log \left( \frac{\eta_i}{\lambda_i} \right) \\
	&=\frac1p\sum_{i=1}^p\log\lambda_i  + \frac{1-c_2}{c_2}\log(1-c_2) - \frac{1-c_1}{c_1}\log(1-c_1)
\end{align*}
where in the last equality we used, among other algebraic simplifications, the fact that $\sum_{i=1}^p \log ( \frac{\eta_i}{\lambda_i} )=\lim_{x\to 0}\log (\frac{\psi_p(x)}{(1-c_1)x})=-\log(1-c_1)$. This is the sought-for result.

\subsection{Development for \lowercase{$f(t)=\log(1+st)$}}

\begin{figure}
	\centering
	\includegraphics[width=\linewidth]{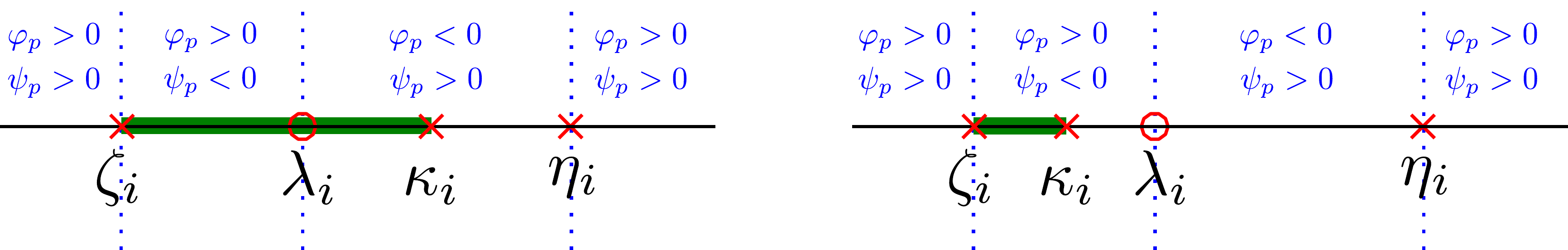}
	\caption{Visual representation of the signs of $\varphi_p$ and $\psi_p$ around singularities for the function $f(t)=\log(1+st)$. Left: case where $\kappa_i>\lambda_i$. Right: case where $\kappa_i>\lambda_i$.}
	\label{fig:sign_phi_over_psi_log1st}
\end{figure}

The development for $f(t)=\log(1+st)$ is quite similar to that of $f(t)=\log(t)$, with some noticeable exceptions with respect to the position of singularity points.

A few important remarks are in order to start with this scenario. First note from Figure~\ref{fig:sign_phi_over_psi} and the previous discussions that the function $z\mapsto \log(1+s\varphi_p(z)/\psi_p(z))$ has a singularity at $z=\kappa_i$, $i=1,\ldots,p$, for some $\kappa_i\in(\zeta_i,\eta_i)$ solution to $1+s\varphi_p(x)/\psi_p(x)=0$ (indeed, $\varphi_p(x)/\psi_p(x)$ is increasing on $(\zeta_i,\eta_i)$ with opposite asymptotes and thus $\kappa_i$ exists and is uniquely defined). In addition, $\log(1+s\varphi_p(z)/\psi_p(z))$ has a further singularity satisfying $1+s\varphi_p(x)/\psi_p(x)=0$ in the interval $(-\infty,0)$ which we shall denote $\kappa_0$. 

A few identities regarding $\kappa_i$ are useful. Using the relation between $\varphi_p$ and $\psi_p$, we find in particular that
\begin{align*}
	\varphi_p(\kappa_i) &= -\frac1s\frac{c_1+c_2-c_1c_2}{c_2}\frac{\kappa_i}{-\frac1s + \frac{c_1}{c_2}\kappa_i} \\
	\psi_p(\kappa_i) &= \frac{c_1+c_2-c_1c_2}{c_2}\frac{\kappa_i}{-\frac1s + \frac{c_1}{c_2}\kappa_i} \\
	(\psi_p+s\varphi_p)\left( \frac{c_2}{c_1s} \right) &= \frac{c_1+c_2-c_1c_2}{c_1}.
\end{align*}

With the discussions above, we also find that
\begin{align}
	\label{eq:log1st_prod}
	1+s\frac{\varphi_p(z)}{\psi_p(z)} &= (1-c_1)s (z-\kappa_0) \frac{\prod_{i=1}^p(z-\kappa_i)}{\prod_{i=1}^p(z-\zeta_i)} \\
	\label{eq:log1st_prod2}
	\psi_p(z)+s\varphi_p(z) &= s(1-c_1)(z-\kappa_0) \frac{\prod_{i=1}^p (z-\kappa_i) }{\prod_{i=1}^p(z-\lambda_i)}.
\end{align}

Note now importantly that $\lambda_i>\frac{c_1}{c_2s}$ is equivalent to $-\frac{c_2}{c_1}\lambda_i<-\frac1s$ which is also $\varphi_p(\lambda_i)/\psi_p(\lambda_i)<\varphi_p(\kappa_i)/\psi_p(\kappa_i)$; then, as $\varphi_p/\psi_p$ is increasing, $\lambda_i>\frac{c_1}{c_2s}$ is equivalent to $\lambda_i<\kappa_i$. On the opposite, for $\lambda_i<\frac{c_1}{c_2s}$, we find $\lambda_i>\kappa_i$. As such, to evaluate the contour integral in this setting, one must isolate two sets of singularities (see Figure~\ref{fig:sign_phi_over_psi_log1st}): (i) those for which $\kappa_i>\lambda_i$ (which are all the largest indices $i$ for which $\lambda_i>\frac{c_1}{c_2s}$) and (ii) those for which $\kappa_i<\lambda_i$. This affects the relative position of the branch cut with respect to $\lambda_i$ and therefore demands different treatments. In particular, the integrals over $I_i^B$ and $I_i^D$ may be restricted to integrals over shorter (possibly empty) segments. Nonetheless, the calculus ultimately reveals that, since the branch cut does not affect the local behavior of the integral around $\lambda_i$, both cases entail the same result. In particular, in case (i) where $\lambda_i>\kappa_i$, recalling \eqref{eq:expansion}, one only has to evaluate
\begin{align*}
	&\int_{\zeta_i+\varepsilon}^{\kappa_i-\varepsilon} \left(\frac{\varphi_p'(x)}{\varphi_p(x)}-\frac{\psi_p'(x)}{\psi_p(x)} \right)\frac{\psi_p(x)}{c_2} dx \nonumber \\
	&= \int_{\zeta_i+\varepsilon}^{\kappa_i-\varepsilon}\left( \frac1p - \frac{c_1+c_2-c_1c_2}{c_1c_2} \right) \sum_{j=1}^p \frac1{x-\lambda_j} + \frac{1-c_2}{c_2} \frac1x \\
	&+ \frac{c_1+c_2-c_1c_2}{c_1c_2}\sum_{j=1}^p \frac1{x-\eta_j} dx \\
	&= \frac1p\sum_{j=1}^p \log \left| \frac{\kappa_i-\lambda_j}{\zeta_i-\lambda_j} \right| + \frac{c_1+c_2-c_1c_2}{c_1c_2}\sum_{j=1}^p \left( \log\left| \frac{\kappa_i-\eta_j}{\kappa_i-\lambda_j} \right| - \log\left| \frac{\zeta_i-\eta_j}{\zeta_i-\lambda_j} \right| \right) \\
	&+ \frac{1-c_2}{c_2}\log \left| \frac{\kappa_i}{\zeta_i} \right|+o(\varepsilon).
\end{align*}
In case~(ii), subdividing the integral as $\int_{\zeta_i+\varepsilon}^{\lambda_i-\varepsilon}+\int_{\lambda_i+\varepsilon}^{\kappa_i-\varepsilon}$ brings immediate simplification of the additional terms in $\lambda_i$ and thus the result remains the same.

The integral over $I_i^C$ is slightly more delicate to handle. In case~(i), in the limit of small $\varepsilon$,
\begin{align*}
	1+ s \frac{\varphi_p}{\psi_p}(\lambda_i+\varepsilon e^{i\theta}) &=  1 - s\frac{c_1}{c_2}\lambda_i + \varepsilon s \frac{c_1}{c_2}\left( p\frac{c_1+c_2-c_1c_2}{c_1c_2} - 1 \right) e^{i\theta} + o(\varepsilon)
\end{align*}
the angle of which is $0+o(\varepsilon)$ uniformly on $\theta\in(-\pi,\pi]$ (since $1 - s\frac{c_1}{c_2}\lambda_i>0$). As such, for all small $\varepsilon$, the sum of the integrals over $(-\pi,0)$ and $(0,\pi]$ reduces to the integral over $(-\pi,\pi]$, leading up to a mere residue calculus, and
\begin{align*}
	&\frac1{2\pi\imath}\oint_{I_i^C} \log\left( 1+s\frac{\varphi_p(z)}{\psi_p(z)} \right) \left(\frac{\varphi_p'(z)}{\varphi_p(z)}-\frac{\psi_p'(z)}{\psi_p(z)} \right)\frac{\psi_p(z)}{c_2} dz \\
	&= \log\left( 1 - s\frac{c_1}{c_2}\lambda_i \right) \left( \frac{c_1+c_2-c_1c_2}{c_1c_2} - \frac1p \right)+o(\varepsilon).
\end{align*}
In case~(ii), $1 - s\frac{c_1}{c_2}\lambda_i<0$ and thus the angle of $1+ s \frac{\varphi_p}{\psi_p}(\lambda_i+\varepsilon e^{i\theta})$ is close to $\pi$; for $\theta\in (0,\pi)$, this leads to an argument equal to $\pi+o(\varepsilon)<\pi$ and for $\theta\in (-\pi,0)$ to an argument equal to $-\pi+o(\varepsilon)>-\pi$. All calculus made, we then find that in either case~(i) or (ii)
\begin{align*}
	&\frac1{2\pi\imath}\oint_{I_i^C} \log\left( 1+s\frac{\varphi_p(z)}{\psi_p(z)} \right) \left(\frac{\varphi_p'(z)}{\varphi_p(z)}-\frac{\psi_p'(z)}{\psi_p(z)} \right)\frac{\psi_p(z)}{c_2} dz \\
	&= \log\left| 1 - s\frac{c_1}{c_2}\lambda_i \right| \left( \frac{c_1+c_2-c_1c_2}{c_1c_2} - \frac1p \right)+o(\varepsilon).
\end{align*}

As in the case of $f(t)=\log(t)$, the integral over $I_i^A$ is of order $o(\varepsilon)$ and vanishes. As a consequence, summing over $i\in\{1,\ldots,p\}$, we find that
\begin{align*}
\frac1{2\pi\imath}\oint_\Gamma &= - \frac1p \sum_{i,j=1}^p \log \left| \frac{\kappa_i-\lambda_j}{\zeta_i-\lambda_j} \right| - \frac{1-c_2}{c_2}\sum_{i=1}^p\log\frac{\kappa_i}{\zeta_i} \\
&+ \frac{c_1+c_2-c_1c_2}{c_1c_2}\sum_{i,j=1}^p \left( \log \left| \frac{\zeta_i-\eta_j}{\zeta_i-\lambda_j}\right| - \log \left| \frac{\kappa_i-\eta_j}{\kappa_i-\lambda_j} \right| \right) \nonumber \\
&-\left( \frac{c_1+c_2-c_1c_2}{c_1c_2}-\frac1p \right)\sum_{i=1}^p\log\left| 1-s\frac{c_1}{c_2}\lambda_j \right|+o(\varepsilon).
\end{align*}

Before reaching the final result, note that, from \eqref{eq:log1st_prod},
\begin{align*}
	&\sum_{i=1}^p \frac1p\sum_{j=1}^p\log \frac{|\kappa_i-\lambda_j|}{|\zeta_i-\lambda_j|} \\ &= \frac1p \sum_{j=1}^p \log \left| \left(1+s\frac{\varphi_p(\lambda_j)}{\psi_p(\lambda_j)} \right)\frac1{\lambda_j-\kappa_0}\frac1{(1-c_1)s} \right| \\
	&= \frac1p\sum_{j=1}^p \log \left|1-\frac{c_1}{c_2}s\lambda_i\right| - \frac1p\sum_{j=1}^p\log(\lambda_j-\kappa_0)-\log((1-c_1)s)
\end{align*}
and similarly 
\begin{align*}
	\sum_{i,j=1}^p \log\left| \frac{\zeta_i-\eta_j}{\zeta_i-\lambda_j}\right| &= \sum_{i=1}^p \log \left| \frac{\varphi_p(\zeta_i)}{(1-c_1)\zeta_i} \right| = p \log \left| \frac{c_1+c_2-c_1c_2}{c_2(1-c_1)} \right| \\
	\sum_{i,j=1}^p \log\left| \frac{\kappa_i-\eta_j}{\kappa_i-\lambda_j}\right| &=  \sum_{i=1}^p \log \left| \frac{\varphi_p(\kappa_i)}{(1-c_1)\kappa_i} \right| = \sum_{i=1}^p\log\left| \frac{c_1+c_2-c_1c_2}{c_2(1-c_1)}\frac1{1-\frac{c_1}{c_2s}\kappa_i} \right| \\
	\sum_{i=1}^p \log \frac{\kappa_i}{\zeta_i} &= \log\left( \frac{1+s\frac{\varphi_p}{\psi_p}(0)}{-(1-c_1)s\kappa_0} \right) = -\log\left(-(1-c_1)s\kappa_0\right).
\end{align*}

Using now \eqref{eq:log1st_prod2}, we find that
	\begin{align*}
		&\sum_{i=1}^p \log\left( \frac{1-s\frac{c_1}{c_2}\lambda_i}{1-s\frac{c_1}{c_2}\kappa_i} \right) = \sum_{i=1}^p \log\left( \frac{\frac{c_2}{c_1s}-\lambda_i}{\frac{c_2}{c_1s}-\kappa_i} \right) = \log \left( \frac{\psi_p\left( \frac{c_2}{c_1s} \right)+s\varphi_p\left( \frac{c_2}{c_1s} \right)}{s(1-c_1)\left( \frac{c_2}{c_1s}-\kappa_0 \right)} \right) \\
		&= \log \left( \frac{c_1+c_2-c_1c_2}{sc_1(1-c_1)\left( \frac{c_2}{c_1s}-\kappa_0 \right)} \right).
	\end{align*}

	Combining the previous results and remarks then leads to
\begin{align*}
	\frac1{2\pi\imath}\oint_{\Gamma} &= \frac{c_1+c_2-c_1c_2}{c_1c_2} \log\left( \frac{c_1+c_2-c_1c_2}{(1-c_1)(c_2-sc_1\kappa_0)} \right) \\
	&+ \frac{1-c_2}{c_2}\log \left( -s\kappa_0(1-c_1) \right) + \log( (1-c_1)s ) + \frac1p\sum_{i=1}^p\log(\lambda_i-\kappa_0) \\
	&= \frac{c_1+c_2-c_1c_2}{c_1c_2} \log\left( \frac{c_1+c_2-c_1c_2}{(1-c_1)(c_2-sc_1\kappa_0)} \right) + \frac1{c_2}\log \left( -s\kappa_0(1-c_1) \right) \\
	&+ \frac1p\sum_{i=1}^p\log\left(1-\frac{\lambda_i}{\kappa_0}\right).
\end{align*}

This concludes the proof for the case $c_1>0$. In the limit where $c_1\to 0$, it suffices to use the Taylor expansion of the leftmost logarithm in the small $c_1$ limit (i.e., $\log(c_2(1-c_1)+c_1)\sim \log(c_2(1-c_1))+c_1/(c_2(1-c_1))$ and $\log(c_2(1-c_1)-sc_1\kappa_0(1-c_1))\sim  \log(c_2(1-c_1))-sc_1\kappa_0/c_2$).

\subsection{Development for $f(t)=\log^2(t)$}

The function $f(t)=\log^2(t)$ is at the core of the Fisher distance and is thus of prime importance in many applications. The evaluation of the complex integral in Theorem~\ref{th:contour_integral} for this case is however quite technical and calls for the important introduction of the dilogarithm function. We proceed with this introduction first and foremost. 

\subsubsection{The dilogarithm function}

The (real) dilogarithm is defined as the function
\begin{align*}
	{\rm Li}_2(x) &= -\int_0^x \frac{\log(1-u)}{u}du.
\end{align*}
for $x\in(-\infty,1]$.


The dilogarithm function will intervene in many instances of the evaluation of the contour integral of Theorem~\ref{th:contour_integral}, through the subsequently defined function $F(X,Y;a)$. This function assumes different formulations depending on the relative position of $X,Y,a$ on the real axis.

\begin{lemma}[Dilogarithm integrals]
	We have the following results and definition 
	\begin{align*}
		(X,Y\geq a>0)\quad & \int_Y^X \frac{\log(x-a)}{x}dx \equiv F(X,Y;a) \\
		&= {\rm Li}_2\left( \frac{a}X \right) - {\rm Li}_2\left( \frac{a}Y \right) + \frac12\left[\log^2(X) - \log^2(Y) \right] \\
		(X,Y>0>a)\quad & \int_Y^X \frac{\log(x-a)}{x}dx \equiv F(X,Y;a) \\
		&= -{\rm Li}_2\left( \frac{X}a \right) + {\rm Li}_2\left( \frac{Y}a \right) + \log\left( \frac{X}Y \right)\log(-a) \\
		(a>X,Y,0~\&~XY>0)\quad & \int_Y^X \frac{\log(a-x)}{x}dx \equiv F(-X,-Y;-a) \\
		&= -{\rm Li}_2\left( \frac{X}a \right) + {\rm Li}_2\left( \frac{Y}a \right) + \log\left( \frac{X}Y \right)\log(a) \\
		(X,Y>0)  \quad & \int_Y^X \frac{\log(x)}{x}dx \equiv F(X,Y;0) \\
		&= \frac12\log^2(X) - \frac12\log^2(Y).
	\end{align*}
\end{lemma}

\begin{lemma}[Properties of Dilogarithm functions {\cite[Section~I-2]{ZAG07}}]
	\label{lem:prop_dilog}
	The following relations hold
	\begin{align*}
		(x<0)\quad &{\rm Li}_2\left( \frac1x \right) + {\rm Li}_2(x) = - \frac12\log^2(-x) - \frac{\pi^2}6 \\
		(0<x<1)\quad &{\rm Li}_2(1-x) + {\rm Li}_2(x) = - \log(x)\log(1-x) + \frac{\pi^2}6 \\
		(0<x<1) \quad &{\rm Li}_2(1-x) + {\rm Li}_2\left( 1-\frac1x \right) = - \frac12\log^2(x).
	\end{align*}
	Besides, for $x<1$ and $\varepsilon>0$ small,
	\begin{align*}
		{\rm Li}_2(x+\varepsilon) &= {\rm Li}_2(x) - \varepsilon \frac{\log(1-x)}x + \varepsilon^2 \frac{(1-x)\log(1-x)+x}{2(1-x)x^2} + O(\varepsilon^3).
	\end{align*}
\end{lemma}

\subsubsection{Integral evaluation}

As in the case where $f(t)=\log(t)$, we shall evaluate the complex integral based on the contour displayed in Figure~\ref{fig:contour}. The main difficulty here arises in evaluating the real integrals over the segments $I_i^B$ and $I_i^D$.

Again, we start from the Equation~\eqref{eq:expansion}. In particular, the integral over $I_i^B$ reads
\begin{align*}
	&\frac1{2\pi\imath}\int_{I_i^B} \log^2\left( \frac{\varphi_p(z)}{\psi_p(z)} \right)  \left(\frac{\varphi_p'(z)}{\varphi_p(z)}-\frac{\psi_p'(z)}{\psi_p(z)} \right)\frac{\psi_p(z)}{c_2}dz \nonumber \\
	&=2\int_{\zeta_i+\varepsilon}^{\lambda_i-\varepsilon} \log\left( - \frac{\varphi_p(x)}{\psi_p(x)} \right)  \left(\frac{\varphi_p'(x)}{\varphi_p(x)}-\frac{\psi_p'(x)}{\psi_p(x)} \right)\frac{\psi_p(x)}{c_2}dx \\
&=2 \int_{\zeta_i+\varepsilon}^{\lambda_i-\varepsilon} \left( \log(1-c_1) + \log (x) + \sum_{l<i}\log(x-\eta_l) \right. \\
&\left. + \sum_{l>i}\log(\eta_l-x) + \log(\eta_i-x) - \sum_{l\leq i} \log(x-\zeta_l) - \sum_{l> i} \log(\zeta_l-x) \right) \\
	&\times \left(  \left( \frac1p - \frac{c_1+c_2-c_1c_2}{c_1c_2} \right) \sum_{j=1}^p \frac1{x-\lambda_j} + \frac{1-c_2}{c_2} \frac1x + \frac{c_1+c_2-c_1c_2}{c_1c_2}\sum_{j=1}^p \frac1{x-\eta_j} \right) dx.
\end{align*}
Note that above we have specifically chosen to write the logarithms in such a way that every integral is a well-defined real integral.

\bigskip

Using now the fact that
\begin{align*}
	\int_Y^X\frac{\log(x-a)}{x-b}dx &= F(X-b,Y-b;a-b) \\
	\int_Y^X\frac{\log(a-x)}{x-b}dx &= F(b-X,b-Y;b-a)
\end{align*}
that we apply repetitively (and very carefully) to the previous equality, we find that the sum of the integral of $I_i^B$ and $I_i^D$ gives
\begin{align*}
	&\frac1{2\pi\imath}\left[\int_{I_i^B}+\int_{I_i^D} \right] \log^2\left( \frac{\varphi_p(z)}{\psi_p(z)} \right)  \left(\frac{\varphi_p'(z)}{\varphi_p(z)}-\frac{\psi_p'(z)}{\psi_p(z)} \right)\frac{\psi_p(z)}{c_2}dz \nonumber \\
	&=2\left[\int_{\zeta_i+\varepsilon}^{\lambda_i-\varepsilon} + \int_{\lambda_i+\varepsilon}^{\eta_i-\varepsilon}\right]\log\left( - \frac{\varphi_p(x)}{\psi_p(x)} \right)  \left(\frac{\varphi_p'(x)}{\varphi_p(x)}-\frac{\psi_p'(x)}{\psi_p(x)} \right)\frac{\psi_p(x)}{c_2}dx \\
&= 2 \left( \frac1p - \frac{c_1+c_2-c_1c_2}{c_1c_2} \right) \left( \log(1-c_1) \sum_{j\neq i} \log\left| \frac{\eta_i-\lambda_j}{\zeta_i-\lambda_j} \right|  \right.\\
& \left. + F(-\varepsilon,\zeta_i-\lambda_i+\varepsilon;-\lambda_i) + F(\eta_i-\lambda_i-\varepsilon,\varepsilon;-\lambda_i)\right.\\
	&\left. + \sum_{j\neq i} F(\eta_i-\lambda_j,\zeta_i-\lambda_j;-\lambda_j) \right.\\
	&\left. +  \sum_{l< i}\sum_{j\neq (i,l)} F(\eta_i-\lambda_j,\zeta_i-\lambda_j;\eta_l-\lambda_j)-F(\eta_i-\lambda_j,\zeta_i-\lambda_j;\zeta_l-\lambda_j) \right. \\
	&\left. +  \sum_{l< i} F(\eta_i-\lambda_l,\zeta_i-\lambda_l;\eta_l-\lambda_l)-F(\eta_i-\lambda_l,\zeta_i-\lambda_l;\zeta_l-\lambda_l) \right. \\
	&\left. + \sum_{l< i} F(-\varepsilon,\zeta_i-\lambda_i+\varepsilon;\eta_l-\lambda_i)+F(\eta_i-\lambda_i-\varepsilon,\varepsilon;\eta_l-\lambda_i) \right. \\
	&\left. - \sum_{l< i} F(-\varepsilon,\zeta_i-\lambda_i+\varepsilon;\zeta_l-\lambda_i)+F(\eta_i-\lambda_i-\varepsilon,\varepsilon;\zeta_l-\lambda_i) \right. \\
	&\left. +  \sum_{l> i}\sum_{j\neq (i,l)} F(-\eta_i+\lambda_j,-\zeta_i+\lambda_j;-\eta_l+\lambda_j)-F(-\eta_i+\lambda_j,-\zeta_i+\lambda_j;-\zeta_l+\lambda_j) \right. \\
	&\left. +  \sum_{l> i} F(-\eta_i+\lambda_l,-\zeta_i+\lambda_l;-\eta_l+\lambda_l)-F(-\eta_i+\lambda_l,-\zeta_i+\lambda_l;-\zeta_l+\lambda_l) \right. \\
	&\left. + \sum_{l> i} F(\varepsilon,-\zeta_i+\lambda_i-\varepsilon;-\eta_l+\lambda_i)+F(-\eta_i+\lambda_i+\varepsilon,-\varepsilon;-\eta_l+\lambda_i) \right. \\
	&\left. - \sum_{l> i} F(\varepsilon,-\zeta_i+\lambda_i-\varepsilon;-\zeta_l+\lambda_i)+F(-\eta_i+\lambda_i+\varepsilon,-\varepsilon;-\zeta_l+\lambda_i) \right. \\
	&\left. + \sum_{j\neq i} F(-\eta_i+\lambda_j,-\zeta_i+\lambda_j;-\eta_i+\lambda_j)-F(\eta_i-\lambda_j,\zeta_i-\lambda_j;\zeta_i-\lambda_j)\right. \\
	&\left. + F(\varepsilon,\lambda_i-\zeta_i-\varepsilon;\lambda_i-\eta_i)+F(\lambda_i-\eta_i+\varepsilon,-\varepsilon;\lambda_i-\eta_i) \right. \\
	&\left. - F(-\varepsilon,\zeta_i-\lambda_i+\varepsilon;\zeta_i-\lambda_i)-F(\eta_i-\lambda_i-\varepsilon,\varepsilon;\zeta_i-\lambda_i) \right) \\
	&+2\frac{1-c_2}{c_2}\left( \log(1-c_1) \log\left( \frac{\eta_i}{\zeta_i} \right) + F(\eta_i,\zeta_i;0) + F(-\eta_i,-\zeta_i;-\eta_i)-F(\eta_i,\zeta_i;\zeta_i) \right. \\
	&\left. + \sum_{l<i} F(\eta_i,\zeta_i;\eta_l)-F(\eta_i,\zeta_i;\zeta_l) + \sum_{l>i} F(-\eta_i,-\zeta_i;-\eta_l)-F(-\eta_i,-\zeta_i;-\zeta_l) \right) \\
	&+2\frac{c_1+c_2-c_1c_2}{c_1c_2} \left( \log(1-c_1)\sum_{j\neq i} \log\left( \frac{\eta_i-\eta_j}{\zeta_i-\eta_j} \right) + \sum_{j\neq i} F(\eta_i-\eta_j,\zeta_i-\eta_j;-\eta_j) \right. \\ 
	&\left. + \sum_{j\neq (i,l)} \sum_{l<i} F(\eta_i-\eta_j,\zeta_i-\eta_j;\eta_l-\eta_j) + \sum_{j\neq (i,l)} \sum_{l>i} F(-\eta_i+\eta_j,-\zeta_i+\eta_j;-\eta_l+\eta_j) \right. \\
	&\left. + \sum_{j\neq i} F(-\eta_i+\eta_j,-\zeta_i+\eta_j;-\eta_i+\eta_j) - \sum_{j\neq (i,l)} \sum_{l<i} F(\eta_i-\eta_j,\zeta_i-\eta_j;\zeta_l-\eta_j) \right. \\
	&\left.- \sum_{j\neq (i,l)} \sum_{l>i} F(-\eta_i+\eta_j,-\zeta_i+\eta_j;-\zeta_l+\eta_j) - \sum_{j\neq i} F(\eta_i-\eta_j,\zeta_i-\eta_j;\zeta_i-\eta_j)\right.\\
	&\left. +\log(1-c_1)\log\left(\frac{\varepsilon}{\eta_i-\zeta_i}\right)+\sum_{l<i}F(-\varepsilon,\zeta_i-\eta_i+\varepsilon;\eta_l-\eta_i)-F(-\varepsilon,\zeta_i-\eta_i+\varepsilon;\zeta_l-\eta_i)\right.\\
	&\left. +\sum_{l>i} F(\varepsilon,\eta_i-\zeta_i-\varepsilon;\eta_i-\eta_l)-F(\varepsilon,\eta_i-\zeta_i-\varepsilon;\eta_i-\zeta_l)\right.\\
	&\left. +F(\varepsilon,-\zeta_i+\eta_i-\varepsilon,\varepsilon)-F(-\varepsilon,\zeta_i-\eta_i+\varepsilon,\zeta_i-\eta_i)+F(-\varepsilon,\zeta_i-\eta_i+\varepsilon,-\eta_i) \right.\\
	&\left. +\sum_{l<i} F(\eta_i-\eta_l,\zeta_i-\eta_l;0)+\sum_{l>i}F(-\eta_i+\eta_l,-\zeta_i+\eta_l;0) \right. \\
	&\left. -\sum_{l<i}F(\eta_i-\eta_l,\zeta_i-\eta_l;\zeta_l-\eta_l)-\sum_{l>i}F(-\eta_i+\eta_l,-\zeta_i+\eta_l;-\zeta_l+\eta_l)
 \right) \\
	&+2\left( \frac1p - \frac{c_1+c_2-c_1c_2}{c_1c_2} \right) \log(1-c_1)\log\left( \frac{\eta_i-\lambda_i}{\lambda_i-\zeta_i} \right) + o_\varepsilon(1).
\end{align*}
To retrieve the expression above, particular care was taken on the relative positions of the $\lambda_{i,j,l}$, $\eta_{i,j,l}$ and $\zeta_{i,j,l}$ to obtain the proper form of the $F$ function; besides, to avoid further complications, a small $\varepsilon$ approximation was used whenever the $F$ function has a finite limit when $\varepsilon\to 0$ (hence the trailing $o_\varepsilon(1)$ in the formula).

To go further, we now make use of the following additional identities obtained from Lemma~\ref{lem:prop_dilog} (these are easily proved).
\begin{lemma}[Properties of the function $F$]
	\label{lem:prop_F}
	We have the following properties of the function $F$:
	\begin{align*}
		(X\geq Y>0) \quad &F(-X,-Y;-X)-F(X,Y;Y) = \frac12\log^2\left(\frac{X}Y\right) \\
		(Y\geq X>0) \quad &F(X,Y;X)-F(-X,-Y;-Y) = \frac12\log^2\left(\frac{X}Y\right) \\
		(X,Y>0) \quad &F(-X+\varepsilon,-\varepsilon;-X)+F(\varepsilon,Y-\varepsilon;-X)\nonumber \\
		&-F(X-\varepsilon,\varepsilon;-Y)-F(-\varepsilon,-Y+\varepsilon;-Y) \nonumber \\
		&= -\frac{\pi^2}2+\frac12\log^2\left(\frac{X}Y\right) + o_\varepsilon(1) \\
		(T,Z\geq X,Y>0) \quad & F(-\varepsilon,-Y+\varepsilon;-T) + F(X-\varepsilon,\varepsilon;-T) \nonumber \\
		&-F(-\varepsilon,-Y+\varepsilon;-Z)-F(X-\varepsilon,\varepsilon;-Z)\\
		&+F(T,Z;-X)-F(T,Z;Y) \nonumber \\
		&= \log\left( \frac{X}{Y} \right)\log\left( \frac{T}Z \right) + o_\varepsilon(1) \\
		(T,Z\geq X,Y>0) \quad  &F(\varepsilon,-Y-\varepsilon;-T) + F(-X+\varepsilon,-\varepsilon;-T) \nonumber \\ 
		&-F(\varepsilon,X-\varepsilon;-Z)-F(-Y+\varepsilon,-\varepsilon;-Z)\\
		&+F(T,Z;X)-F(T,Z;-Y) \nonumber \\
		&= \log\left( \frac{X}{Y} \right)\log\left( \frac{T}Z \right) + o_\varepsilon(1) \\
		(X,Y,Z,T>0) \quad & F(X,Y;T)+F(T,Z;X)-F(X,Y;Z)-F(T,Z;Y) \nonumber \\
		&= \log\left(\frac{X}Y\right)\log\left( \frac{T}Z \right) \\
		(XY>0 ~\&~ ZT>0) \quad & F(X,Y;-Z)+F(Z,T;-X)-F(X,Y;-T)-F(Z,T;-Y) \nonumber \\
		&= \log\left(\frac{X}Y\right)\log\left( \frac{T}Z \right)\\
		(Y,Z>0 ~\&~ \varepsilon X>0) \quad & F(\varepsilon,X;-Y)+F(\varepsilon,X;-Z)-F(Y,Z;-X) \nonumber \\
		&= \log\left(\frac{\varepsilon}X\right)\log\left( \frac{Y}Z \right)+\frac12 \log^{2}(Z)-\frac12 \log^{2}(Y).
	\end{align*}
\end{lemma}

Exploiting the relations from the previous lemma, we have the following first result:
\begin{align*}
	&\sum_{i=1}^p F(-\varepsilon,\zeta_i-\lambda_i+\varepsilon;-\lambda_i) + F(\eta_i-\lambda_i-\varepsilon,\varepsilon;-\lambda_i) \\
	&=\sum_{i=1}^p{\rm Li}_2\left(1-\frac{\zeta_i}{\lambda_i}\right)-{\rm Li}_2\left(1-\frac{\eta_i}{\lambda_i}\right)+\log(\lambda_i)\log\left(\frac{\eta_i-\lambda_i}{\lambda_i-\zeta_i}\right)
\end{align*}
The terms involving double or triple sums (over $i,l$ or $i,j,l$) are more subtle to handle. By observing that $\sum_i\sum_{l>i} G_{il}=\sum_l\sum_{i<l} G_{il}$ which, up to a switch in the notation $(i,l)$ into $(l,i)$, is the same as $\sum_i\sum_{l<i} G_{li}$, we have that
\begin{align*}
	\sum_i\sum_{l>i} G_{il} + \sum_i\sum_{l<i} G_{il} &= \sum_i\sum_{l>i} G_{il}+G_{li}.
\end{align*}
Using this observation to gather terms together, we find notably from Lemma~\ref{lem:prop_F} that
\begin{align*}
&\sum_i\sum_{l< i}\sum_{j\neq (i,l)} F(\eta_i-\lambda_j,\zeta_i-\lambda_j;\eta_l-\lambda_j)-F(\eta_i-\lambda_j,\zeta_i-\lambda_j;\zeta_l-\lambda_j)\\
&+\sum_i\sum_{l> i}\sum_{j\neq (i,l)} F(-\eta_i+\lambda_j,-\zeta_i+\lambda_j;-\eta_l+\lambda_j)-F(-\eta_i+\lambda_j,-\zeta_i+\lambda_j;-\zeta_l+\lambda_j) \\
&= \sum_i\sum_{l<i}\sum_{j\neq (i,l)}\log\left(\frac{\lambda_j-\eta_i}{\lambda_j-\zeta_i}\right)\log\left(\frac{\lambda_j-\eta_l}{\lambda_j-\zeta_l}\right)\\
&=\frac12\sum_{i}\sum_{l\neq i}\sum_{j\neq (i,l)}\log\left(\frac{\lambda_j-\eta_i}{\lambda_j-\zeta_i}\right)\log\left(\frac{\lambda_j-\eta_l}{\lambda_j-\zeta_l}\right)
\end{align*}
Similarly,
\begin{align*}
	&\sum_i \sum_{l< i} F(-\varepsilon,\zeta_i-\lambda_i+\varepsilon;\eta_l-\lambda_i)+F(\eta_i-\lambda_i-\varepsilon,\varepsilon;\eta_l-\lambda_i) \\
	&-\sum_i \sum_{l< i} F(-\varepsilon,\zeta_i-\lambda_i+\varepsilon;\zeta_l-\lambda_i)+F(\eta_i-\lambda_i-\varepsilon,\varepsilon;\zeta_l-\lambda_i) \nonumber\\
	&+ \sum_i \sum_{l> i} F(-\eta_i+\lambda_l,-\zeta_i+\lambda_l;-\eta_l+\lambda_l)-F(-\eta_i+\lambda_l,-\zeta_i+\lambda_l;-\zeta_l+\lambda_l) \\
	&=\sum_i \sum_{l< i} F(-\varepsilon,\zeta_i-\lambda_i+\varepsilon;\eta_l-\lambda_i)+F(\eta_i-\lambda_i-\varepsilon,\varepsilon;\eta_l-\lambda_i) \\
	&-\sum_i \sum_{l< i} F(-\varepsilon,\zeta_i-\lambda_i+\varepsilon;\zeta_l-\lambda_i)+F(\eta_i-\lambda_i-\varepsilon,\varepsilon;\zeta_l-\lambda_i) \nonumber\\
	&+ \sum_i \sum_{l<i} F(-\eta_l+\lambda_i,-\zeta_l+\lambda_i;-\eta_i+\lambda_i)-F(-\eta_l+\lambda_i,-\zeta_l+\lambda_i;-\zeta_i+\lambda_i) \\
	&=\sum_i \sum_{l< i} \log\left( \frac{\eta_i-\lambda_i}{\zeta_i-\lambda_i} \right)\log\left( \frac{\lambda_i-\eta_l}{\lambda_i-\zeta_l} \right)+o_\varepsilon(1)
\end{align*}
and, symmetrically,
\begin{align*}
	& \sum_i\sum_{l> i} F(\varepsilon,-\zeta_i+\lambda_i-\varepsilon;-\eta_l+\lambda_i)+F(-\eta_i+\lambda_i+\varepsilon,-\varepsilon;-\eta_l+\lambda_i) \nonumber \\
	&-\sum_i \sum_{l> i} F(\varepsilon,-\zeta_i+\lambda_i-\varepsilon;-\zeta_l+\lambda_i)+F(-\eta_i+\lambda_i+\varepsilon,-\varepsilon;-\zeta_l+\lambda_i) \nonumber \\
	&+\sum_i \sum_{l< i} F(\eta_i-\lambda_l,\zeta_i-\lambda_l;\eta_l-\lambda_l)-F(\eta_i-\lambda_l,\zeta_i-\lambda_l;\zeta_l-\lambda_l) \\
	&=\sum_i\sum_{l> i} F(\varepsilon,-\zeta_i+\lambda_i-\varepsilon;-\eta_l+\lambda_i)+F(-\eta_i+\lambda_i+\varepsilon,-\varepsilon;-\eta_l+\lambda_i) \nonumber \\
	&-\sum_i \sum_{l> i} F(\varepsilon,-\zeta_i+\lambda_i-\varepsilon;-\zeta_l+\lambda_i)+F(-\eta_i+\lambda_i+\varepsilon,-\varepsilon;-\zeta_l+\lambda_i) \nonumber \\
	&+\sum_i \sum_{l> i} F(\eta_l-\lambda_i,\zeta_l-\lambda_i;\eta_i-\lambda_i)-F(\eta_l-\lambda_i,\zeta_l-\lambda_i;\zeta_i-\lambda_i) \\
	&=\sum_i \sum_{l> i} \log\left( \frac{\eta_i-\lambda_i}{\zeta_i-\lambda_i} \right)\log\left( \frac{\lambda_i-\eta_l}{\lambda_i-\zeta_l} \right)+o_\varepsilon(1).
\end{align*}

Also, using Items 1 and 2 of Lemma~\ref{lem:prop_F}, we find
\begin{align*}
&\sum_i\sum_{j\neq i} F(-\eta_i+\lambda_j,-\zeta_i+\lambda_j;-\eta_i+\lambda_j)-F(\eta_i-\lambda_j,\zeta_i-\lambda_j;\zeta_i-\lambda_j)\\
&=\sum_i\sum_{j\neq i}\frac12 \log^{2}\left(\frac{\lambda_j-\eta_i}{\lambda_j-\zeta_i}\right).
\end{align*}
Again from Lemma~\ref{lem:prop_F}, we also have
\begin{align*}
	&F(\varepsilon,\lambda_i-\zeta_i-\varepsilon;\lambda_i-\eta_i)+F(\lambda_i-\eta_i+\varepsilon,-\varepsilon;\lambda_i-\eta_i)\nonumber \\
	&- F(-\varepsilon,\zeta_i-\lambda_i+\varepsilon;\zeta_i-\lambda_i)-F(\eta_i-\lambda_i-\varepsilon,\varepsilon;\zeta_i-\lambda_i) \\
	&= -\frac{\pi^2}2 + \frac12\log^2\left( \frac{\eta_i-\lambda_i}{\lambda_i-\zeta_i} \right).
\end{align*}
Before going further, remark that the last four established relations can be assembled to reach
\begin{align*}
	&\left. \sum_{j\neq i} F(\eta_i-\lambda_j,\zeta_i-\lambda_j;-\lambda_j) +  \sum_{l< i}\sum_{j\neq (i,l)} F(\eta_i-\lambda_j,\zeta_i-\lambda_j;\eta_l-\lambda_j)\right. \\
	&\left. -F(\eta_i-\lambda_j,\zeta_i-\lambda_j;\zeta_l-\lambda_j) \right. \\
	&\left. +  \sum_{l< i} F(\eta_i-\lambda_l,\zeta_i-\lambda_l;\eta_l-\lambda_l)-F(\eta_i-\lambda_l,\zeta_i-\lambda_l;\zeta_l-\lambda_l) \right. \\
	&\left. + \sum_{l< i} F(-\varepsilon,\zeta_i-\lambda_i+\varepsilon;\eta_l-\lambda_i)+F(\eta_i-\lambda_i-\varepsilon,\varepsilon;\eta_l-\lambda_i) \right. \\
	&\left. - \sum_{l< i} F(-\varepsilon,\zeta_i-\lambda_i+\varepsilon;\zeta_l-\lambda_i)+F(\eta_i-\lambda_i-\varepsilon,\varepsilon;\zeta_l-\lambda_i) \right. \\
	&\left. +  \sum_{l> i}\sum_{j\neq (i,l)} F(-\eta_i+\lambda_j,-\zeta_i+\lambda_j;-\eta_l+\lambda_j)-F(-\eta_i+\lambda_j,-\zeta_i+\lambda_j;-\zeta_l+\lambda_j) \right. \\
	&\left. +  \sum_{l> i} F(-\eta_i+\lambda_l,-\zeta_i+\lambda_l;-\eta_l+\lambda_l)-F(-\eta_i+\lambda_l,-\zeta_i+\lambda_l;-\zeta_l+\lambda_l) \right. \\
	&\left. + \sum_{l> i} F(\varepsilon,-\zeta_i+\lambda_i-\varepsilon;-\eta_l+\lambda_i)+F(-\eta_i+\lambda_i+\varepsilon,-\varepsilon;-\eta_l+\lambda_i) \right. \\
	&\left. - \sum_{l> i} F(\varepsilon,-\zeta_i+\lambda_i-\varepsilon;-\zeta_l+\lambda_i)+F(-\eta_i+\lambda_i+\varepsilon,-\varepsilon;-\zeta_l+\lambda_i) \right. \\
	&\left. + \sum_{j\neq i} F(-\eta_i+\lambda_j,-\zeta_i+\lambda_j;-\eta_i+\lambda_j)-F(\eta_i-\lambda_j,\zeta_i-\lambda_j;\zeta_i-\lambda_j)\right. \\
&\left. + F(\varepsilon,\lambda_i-\zeta_i-\varepsilon;\lambda_i-\eta_i)+F(\lambda_i-\eta_i+\varepsilon,-\varepsilon;\lambda_i-\eta_i)\right.\\
	&\left. - F(-\varepsilon,\zeta_i-\lambda_i+\varepsilon;\zeta_i-\lambda_i)-F(\eta_i-\lambda_i-\varepsilon,\varepsilon;\zeta_i-\lambda_i)\right.\\
	&=\frac12 p\log^{2}\left(\frac{c_1}{c_2(1-c_1)}\right)+o_\varepsilon(1).
\end{align*}

Next, we have
\begin{align*}
	F(\eta_i,\zeta_i;0) &= \frac12\log^2(\eta_i)-\frac12\log^2(\zeta_i)
\end{align*}
and, again by Lemma~\ref{lem:prop_F},
\begin{align*}
&\sum_i \sum_{l\neq i} F(\eta_i,\zeta_i;\eta_l)-F(\eta_i,\zeta_i;\zeta_l) \\
&= \sum_{i>l} F(\eta_i,\zeta_i;\eta_l)+F(\eta_l,\zeta_l;\eta_i)-F(\eta_i,\zeta_i;\zeta_l) -F(\eta_l,\zeta_l;\zeta_i) \\
&= \sum_{i>l} \log\left( \frac{\eta_i}{\zeta_i} \right)\log\left( \frac{\eta_l}{\zeta_l} \right) \\
&= \frac12 \sum_i\sum_{l\neq i} \log\left( \frac{\eta_i}{\zeta_i} \right)\log\left( \frac{\eta_l}{\zeta_l} \right) \\
&\sum_i F(-\eta_i,-\zeta_i;-\eta_i)-F(\eta_i,\zeta_i;\zeta_i) = \frac12 \sum_i \log^2\left( \frac{\eta_i}{\zeta_i} \right)
\end{align*}
so that
\begin{align*}
	&\sum_i \sum_{l\neq i} F(\eta_i,\zeta_i;\eta_l)-F(\eta_i,\zeta_i;\zeta_l) + \sum_i F(-\eta_i,-\zeta_i;-\eta_i)-F(\eta_i,\zeta_i;\zeta_i) \\
	&= \frac12 \left( \sum_i \log\left( \frac{\eta_i}{\zeta_i} \right) \right)^2.
\end{align*}
Recall now the already established identity $\sum_i \log ( \frac{\eta_i}{\zeta_i}) = -\log((1-c_1)(1-c_2))$ from which
\begin{align*}
	&\sum_i \sum_{l\neq i} F(\eta_i,\zeta_i;\eta_l)-F(\eta_i,\zeta_i;\zeta_l) + \sum_i F(-\eta_i,-\zeta_i;-\eta_i)-F(\eta_i,\zeta_i;\zeta_i) \\
	&= \frac12 \log^2( (1-c_1)(1-c_2) ).
\end{align*}

\bigskip

Continuing, we also have
\begin{align*}
	\sum_{j\neq i} \log\left( \frac{\eta_i-\lambda_j}{\zeta_i-\lambda_j} \right) &= \log \left( -\frac{c_1}{c_2(1-c_1)} \right) - \log\left( \frac{\eta_i-\lambda_i}{\zeta_i-\lambda_i} \right) = \log\left( \frac{c_1}{c_2(1-c_1)}\frac{\eta_i-\lambda_i}{\lambda_i-\zeta_i} \right).
\end{align*}
By Lemma~\ref{lem:prop_F} again, we next find
\begin{align*}
&\sum_i\sum_{l>i} F(-\eta_i,-\zeta_i;-\eta_l)-F(-\eta_i,-\zeta_i;-\zeta_l) +\sum_{l<i} F(\eta_i,\zeta_i;\eta_l)-F(\eta_i,\zeta_i;\zeta_l) \\
&=\sum_{i}\sum_{l<i}\log\left(\frac{\eta_i}{\zeta_i}\right)\log\left(\frac{\eta_l}{\zeta_l}\right)
\end{align*}
from which we deduce that 
\begin{align*}
&\sum_i F(-\eta_i,-\zeta_i,-\eta_i)-F(\eta_i,\zeta_i;\zeta_i)+\sum_i\sum_{l>i} F(-\eta_i,-\zeta_i;-\eta_l) \nonumber \\
&-F(-\eta_i,-\zeta_i;-\zeta_l) +\sum_{l<i} F(\eta_i,\zeta_i;\eta_l)-F(\eta_i,\zeta_i;\zeta_l)\\
&=\frac12\sum_{i,l}\log\left(\frac{\eta_i}{\zeta_i}\right)\log\left(\frac{\eta_l}{\zeta_l}\right)\\
&=\frac12\log^{2}\left((1-c_1)(1-c_2)\right).
\end{align*}

The next term also simplifies through the definition of $\varphi_p/\psi_p$:
\begin{align*}
	\sum_i \sum_{j\neq i}\log\left|\frac{\eta_i-\eta_j}{\zeta_i-\eta_j}\right|&=\sum_{j}\log\left(\frac{c_1}{c_1+c_2-c_1c_2}\varphi_p'(\eta_j)\frac{(\zeta_j-\eta_j)}{(1-c_1)\eta_j}\right).
\end{align*}

Still from Lemma~\ref{lem:prop_F} and with the same connection to $\varphi_p/\psi_p$, we have
\begin{align*}
&\sum_j \sum_{l<i} F(\eta_i-\eta_j,\zeta_i-\eta_j;\eta_l-\eta_j)+ \sum_j \sum_{l>i} F(-\eta_i+\eta_j,-\zeta_i+\eta_j;-\eta_l+\eta_j)  \\
&+ \sum_{j\neq i} F(-\eta_i+\eta_j,-\zeta_i+\eta_j;-\eta_i+\eta_j) - \sum_j \sum_{l<i} F(\eta_i-\eta_j,\zeta_i-\eta_j;\zeta_l-\eta_j) \\
&- \sum_j \sum_{l>i} F(-\eta_i+\eta_j,-\zeta_i+\eta_j;-\zeta_l+\eta_j) - \sum_{j\neq i} F(\eta_i-\eta_j,\zeta_i-\eta_j;\zeta_i-\eta_j)\\
&=\sum_i\sum_{j\neq (i,l)}\sum_{l<i}\log\left(\frac{\eta_i-\eta_j}{\zeta_i-\eta_j}\right)\log\left(\frac{\eta_l-\eta_j}{\zeta_l-\eta_j}\right)+\frac12\log^{2}\left(\frac{\eta_j-\eta_i}{\eta_j-\zeta_i}\right)\\
&=\frac12\sum_i\sum_{j\neq (i,l)}\sum_{l\neq i}\log\left(\frac{\eta_i-\eta_j}{\zeta_i-\eta_j}\right)\log\left(\frac{\eta_l-\eta_j}{\zeta_l-\eta_j}\right)+\frac12\log^{2}\left(\frac{\eta_j-\eta_i}{\eta_j-\zeta_i}\right)\\
&=\frac12\sum_i \sum_l\sum_{j\neq (i,l)}\log\left(\frac{\eta_i-\eta_j}{\zeta_i-\eta_j}\right)\log\left(\frac{\eta_l-\eta_j}{\zeta_l-\eta_j}\right)\\
&=\frac12\sum_{j}\log^{2}\left(\frac{c_1}{c_1+c_2-c_1c_2}\varphi'(\eta_j)\frac{(\zeta_j-\eta_j)}{(1-c_1)\eta_j}\right).
\end{align*}
Again from Lemma~\ref{lem:prop_F}, we next have
\begin{align*}
&\sum_{i}\sum_{l<i}F(-\varepsilon,\zeta_i-\eta_i+\varepsilon;\eta_l-\eta_i)-F(-\varepsilon,\zeta_i-\eta_i+\varepsilon;\zeta_l-\eta_i) \\
&-\sum_{l>i}F(-\eta_i+\eta_l,-\zeta_i+\eta_l;-\zeta_l+\eta_l)\\
&=\sum_{i}\sum_{l<i}\log\left(\frac{\varepsilon}{\eta_i-\zeta_i}\right)\log\left(\frac{\eta_i-\eta_l}{\eta_i-\zeta_l}\right)-\frac12\log^{2}(\eta_i-\eta_l)+\frac12\log^{2}(\eta_i-\zeta_l)+o_\varepsilon(1).
\end{align*}
We also have the following relations
\begin{align*}
&\sum_{i}\sum_{l>i} F(\varepsilon,\eta_i-\zeta_i-\varepsilon;\eta_i-\eta_l)-F(\varepsilon,\eta_i-\zeta_i-\varepsilon;\eta_i-\zeta_l)\\
&-\sum_{l<i}F(\eta_i-\eta_l,\zeta_i-\eta_l;\zeta_l-\eta_l)\\
&=\sum_{i}\sum_{l>i}\frac12 \log^{2}(\zeta_l-\eta_i)-\frac12 \log^{2}(\eta_l-\eta_i)+\log\left(\frac{\eta_l-\eta_i}{\zeta_l-\eta_i}\right)\log\left(\frac{\varepsilon}{\eta_i-\zeta_i}\right)+o_\varepsilon(1)
\end{align*}
and
\begin{align*}
&\sum_{l<i}F(-\varepsilon,\zeta_i-\eta_i+\varepsilon;\eta_l-\eta_i)-F(-\varepsilon,\zeta_i-\eta_i+\varepsilon;\zeta_l-\eta_i)\\
&-\sum_{l>i}F(-\eta_i+\eta_l,-\zeta_i+\eta_l;-\zeta_l+\eta_l)\\
&=\sum_i\sum_{l< i}\log\left(\frac{\varepsilon}{\eta_i-\zeta_i}\right)\log\left(\frac{\eta_i-\eta_l}{\eta_i-\zeta_l}\right)-\frac12\log^{2}(\eta_i-\eta_l)+\frac12\log^{2}(\eta_i-\zeta_l)+o_\varepsilon(1)
\end{align*}
which together gives
\begin{align*}
&\sum_{i}\sum_{l>i} F(\varepsilon,\eta_i-\zeta_i-\varepsilon;\eta_i-\eta_l)-F(\varepsilon,\eta_i-\zeta_i-\varepsilon;\eta_i-\zeta_l) \\ 
&-\sum_{l<i}F(\eta_i-\eta_l,\zeta_i-\eta_l;\zeta_l-\eta_l) +\sum_{l<i}F(-\varepsilon,\zeta_i-\eta_i+\varepsilon;\eta_l-\eta_i) \\
&-F(-\varepsilon,\zeta_i-\eta_i+\varepsilon;\zeta_l-\eta_i)-\sum_{l>i}F(-\eta_i+\eta_l,-\zeta_i+\eta_l;-\zeta_l+\eta_l)\\
&=\sum_i\sum_{l\neq i}\log\left(\frac{\varepsilon}{\eta_i-\zeta_i}\right)\log\left(\frac{\eta_i-\eta_l}{\eta_i-\zeta_l}\right).
\end{align*}

The next term is 
\begin{align*}
F(\varepsilon,-\zeta_i+\eta_i-\varepsilon,\varepsilon)-F(-\varepsilon,\zeta_i-\eta_i+\varepsilon,\zeta_i-\eta_i)&=\frac12 \log^{2}\left(\frac{\varepsilon}{\eta_i-\zeta_i}\right)
\end{align*}
and finally the last term gives
\begin{align*}
	&\sum_{l<i} F(\eta_i-\eta_l,\zeta_i-\eta_l;0)+\sum_{l>i}F(-\eta_i+\eta_l,-\zeta_i+\eta_l;0) \\
	&=\frac12\sum_i \sum_{l\neq i} \log^{2}\left|\eta_i-\eta_l\right|-\log^{2}\left|\zeta_i-\eta_l\right|.
\end{align*}

Putting all results above together, we obtain
\begin{align*}
	&\sum_{i=1}^p 2\left[\int_{\zeta_i+\varepsilon}^{\lambda_i-\varepsilon} + \int_{\lambda_i+\varepsilon}^{\eta_i-\varepsilon}\right]\log\left( - \frac{\varphi_p(x)}{\psi_p(x)} \right)  \left(\frac{\varphi_p'(x)}{\varphi_p(x)}-\frac{\psi_p'(x)}{\psi_p(x)} \right)\frac{\psi_p(x)}{c_2}dx \\
&= 2 \left( \frac1p - \frac{c_1+c_2-c_1c_2}{c_1c_2} \right) \left( p\log(1-c_1) \log\left(\frac{c_1}{c_{2}(1-c_1)}\right) + \sum_{i}{\rm Li}_2\left(1-\frac{\zeta_i}{\lambda_i}\right) \right. \\
	&\left.-{\rm Li}_2\left(1-\frac{\eta_i}{\lambda_i}\right)+\log(\lambda_i)\log\left(\frac{\eta_i-\lambda_i}{\lambda_i-\zeta_i}\right)\right. \\
	&\left. + \sum_i\sum_{j\neq i} F(\eta_i-\lambda_j,\zeta_i-\lambda_j;-\lambda_j)+\frac12p\log^{2}\left(\frac{c_1}{c_2(1-c_1)}\right)-\frac{\pi^2}{2}p\right) \\
&+2\frac{1-c_2}{c_2}\left(-\log(1-c_1) \log\left((1-c_1)(1-c_2)\right) + \sum_{i}\frac12\log^2(\eta_i)\right. \\
	&\left. -\frac12\log^2(\zeta_i)  +\frac12\log^{2}\left((1-c_1)(1-c_2)\right)  \right) \\
	&+2\frac{c_1+c_2-c_1c_2}{c_1c_2} \left( \log(1-c_1)\sum_i\sum_{j\neq i} \log\left( \frac{\eta_i-\eta_j}{\zeta_i-\eta_j} \right) + \sum_i\sum_{j\neq i} F(\eta_i-\eta_j,\zeta_i-\eta_j;-\eta_j)\right.\\
	&\left. +\frac12\sum_i \sum_l\sum_{j\neq (i,l)}\log\left(\frac{\eta_i-\eta_j}{\zeta_i-\eta_j}\right)\log\left(\frac{\eta_l-\eta_j}{\zeta_l-\eta_j}\right) +\sum_i\log(1-c_1)\log\left(\frac{\varepsilon}{\eta_i-\zeta_i}\right)\right. \\
	&\left.+\sum_{i}\sum_{l\neq i}\log\left|\frac{\eta_l-\eta_i}{\zeta_l-\eta_i}\right|\log\left|\frac{\varepsilon}{\eta_i-\zeta_i}\right| \right.\\
	&\left.+\sum_i {\rm Li}_2\left(1-\frac{\zeta_i}{\eta_i}\right)+\log\left(\frac{\varepsilon}{\eta_i-\zeta_i}\right)\log(\eta_i)+\sum_i \frac12 \log^{2}\left(\frac{\varepsilon}{\eta_i-\zeta_i}\right) -p\frac{\pi^2}{6}\right) + o_\varepsilon(1).
\end{align*}

The integral over the contour $I_i^E$ can be computed using the same reasoning as for the $f(t)=\log(t)$ function and is easily obtained as
\begin{align*}
	\frac{1}{2\pi \imath}\int_{I_i^E} &= \frac{c_1+c_2-c_1c_2}{c_1c_2} \left( -\log^{2}\left( \varphi_p'(\eta_i) \frac{c_1}{c_1+c_2-c_1c_2} \right) -\log^{2}(\varepsilon) \right. \\ &\left. -2\log(\varepsilon)\left[\left(\frac{\varphi}{\psi}\right)'(\eta_i)\right]  -\frac{\pi^2}{3}\right) +O(\varepsilon\log^{2}(\varepsilon)).
\end{align*}

Adding up the ``residue'' at $\lambda_i$ (i.e., the integral over $I_i^C$), we end up with the following expression for the sought-for integral
\begin{align*}
	&\frac1{2\pi\imath}\oint_\Gamma\log^{2}\left( \frac{\varphi_p(z)}{\psi_p(z)} \right)  \left(\frac{\varphi_p'(z)}{\varphi_p(z)}-\frac{\psi_p'(z)}{\psi_p(z)} \right)\frac{\psi_p(z)}{c_2}dz \\
	&= -2 \left( \frac1p - \frac{c_1+c_2-c_1c_2}{c_1c_2} \right) \left( p\log(1-c_1) \log\left(\frac{c_1}{c_{2}(1-c_1)}\right) \right. \\ 
	&\left.+ \sum_{i}{\rm Li}_2\left(1-\frac{\zeta_i}{\lambda_i}\right)-{\rm Li}_2\left(1-\frac{\eta_i}{\lambda_i}\right)+\log(\lambda_i)\log\left(\frac{\eta_i-\lambda_i}{\lambda_i-\zeta_i}\right)\right. \\
	&\left. + \sum_i\sum_{j\neq i} F(\eta_i-\lambda_j,\zeta_i-\lambda_j;-\lambda_j)+p\log^{2}\left(\frac{c_1}{c_2(1-c_1)}\right)-\frac{\pi^2}{2}p\right) \\
	&-2\frac{1-c_2}{c_2}\left(-\log(1-c_1) \log\left((1-c_1)(1-c_2)\right) + \sum_{i}\frac12\log^2(\eta_i)-\frac12\log^2(\zeta_i) \right. \\
	&\left.+\frac12\log^{2}\left((1-c_1)(1-c_2)\right)  \right) \\
	&-2\frac{c_1+c_2-c_1c_2}{c_1c_2} \left( \log(1-c_1)\sum_{j}\log\left(\frac{c_1}{c_1+c_2-c_1c_2}\varphi'(\eta_j)\frac{(\zeta_j-\eta_j)}{(1-c_1)\eta_j }\right) \right. \\
&\left.+ \sum_i\sum_{j\neq i} F(\eta_i-\eta_j,\zeta_i-\eta_j;-\eta_j)\right.\\
	&\left. \frac12\sum_{j}\log^{2}\left(\frac{c_1}{c_1+c_2-c_1c_2}\varphi'(\eta_j)\frac{(\zeta_j-\eta_j)}{(1-c_1)\eta_j}\right) +\sum_i {\rm Li}_2\left(1-\frac{\zeta_i}{\eta_i}\right)\right.\\
	&\left.+\sum_{i}\log\left[\left(\frac{\varphi}{\psi}\right)'(\eta_i)\right]\log\left|\frac{\varepsilon}{\eta_i-\zeta_i}\right|+\log(\eta_i-\zeta_i)\log\left|\frac{\varepsilon}{\eta_i-\zeta_i}\right| \right.\\
	&\left.+\sum_i\frac12 \log^{2}\left(\frac{\varepsilon}{\eta_i-\zeta_i}\right)-\frac{\pi^2}{6} \right) + o_\varepsilon(1)\\
	&+\sum_i \frac{c_1+c_2-c_1c_2}{c_1c_2} \left( -\log^{2}\left( \varphi_p'(\eta_i) \frac{c_1}{c_1+c_2-c_1c_2} \right) -\log^{2}(\varepsilon)\right. \\
	&\left.-2\log(\varepsilon)\log\left[\left(\frac{\varphi}{\psi}\right)'(\eta_i)\right]-\frac{\pi^2}{3} \right)\\
	&+\sum_i \left(\log^{2}\left(\frac{c_1}{c_2}\lambda_{i}\right)-\pi^2\right)\left[ \frac{c_1+c_2-c_1c_2}{c_1c_2} - \frac1p \right] + o_\varepsilon(1).	
\end{align*}
which, in the limit of small $\varepsilon$, can be simplified as
\begin{align*}
	&\frac1{2\pi\imath}\oint_\Gamma\log^{2}\left( \frac{\varphi_p(z)}{\psi_p(z)} \right)  \left(\frac{\varphi_p'(z)}{\varphi_p(z)}-\frac{\psi_p'(z)}{\psi_p(z)} \right)\frac{\psi_p(z)}{c_2}dz \\
	&= -2 \left( \frac1p - \frac{c_1+c_2-c_1c_2}{c_1c_2} \right) \left(p \log(1-c_1) \log\left(\frac{c_1}{c_{2}(1-c_1)}\right)\right. \\ 
	&\left.+ \sum_{i}{\rm Li}_2\left(1-\frac{\zeta_i}{\lambda_i}\right)-{\rm Li}_2\left(1-\frac{\eta_i}{\lambda_i}\right)+\log(\lambda_i)\log\left(\frac{\eta_i-\lambda_i}{\lambda_i-\zeta_i}\right)\right. \\
	&\left. + \sum_i\sum_{j\neq i} F(\eta_i-\lambda_j,\zeta_i-\lambda_j;-\lambda_j)+p\log^{2}\left(\frac{c_1}{c_2(1-c_1)}\right)\right) \\
	&-2\frac{1-c_2}{c_2}\left(-\log(1-c_1) \log\left((1-c_1)(1-c_2)\right) + \sum_i\frac12\log^2(\eta_i)-\frac12\log^2(\zeta_i) \right. \\
	&\left. +\frac12\log^{2}\left((1-c_1)(1-c_2)\right)  \right) \\
	&-2\frac{c_1+c_2-c_1c_2}{c_1c_2} \left( \log(1-c_1)\sum_{j}\log\left(\frac{c_1}{c_1+c_2-c_1c_2}\varphi'(\eta_j)\frac{(\zeta_j-\eta_j)}{(1-c_1)\eta_j }\right) \right. \\
	&\left.+ \sum_i\sum_{j\neq i} F(\eta_i-\eta_j,\zeta_i-\eta_j;-\eta_j)\right.\\
	&\left.-\sum_i \log\left[\left(\frac{\varphi}{\psi}'(\eta_i)\right)\right]\log\left(\eta_i-\zeta_i\right)+\sum_i {\rm Li}_2\left(1-\frac{\zeta_i}{\eta_i}\right)-\frac12 \log^{2}\left(\eta_i-\zeta_i\right)\right) \\
	&-\sum_i \frac{c_1+c_2-c_1c_2}{c_1c_2} \left( -\log^{2}\left(\frac{\eta_i-\zeta_i}{(1-c_1)\eta_i} \right) \right. \\
	&\left.+2\log\left(\frac{\eta_i-\zeta_i}{(1-c_1)\eta_i} \right)\log\left(\frac{c_1}{c_1+c_2-c_1c_2}\varphi'(\eta_i)\frac{(\zeta_i-\eta_i)}{(1-c_1)\eta_i}\right) \right)\\
	&-\sum_i \log^{2}\left(\frac{c_1}{c_2}\lambda_{i}\right)\left[ \frac{c_1+c_2-c_1c_2}{c_1c_2} - \frac1p \right]. 	
\end{align*}
After further book-keeping and simplifications, we ultimately find:
\begin{align}
	\label{eq:log2t_exact}
	&\frac1{2\pi\imath}\oint_\Gamma\log^{2}\left( \frac{\varphi_p(z)}{\psi_p(z)} \right)  \left(\frac{\varphi_p'(z)}{\varphi_p(z)}-\frac{\psi_p'(z)}{\psi_p(z)} \right)\frac{\psi_p(z)}{c_2}dz \\
&= \frac{c_1+c_2-c_1c_2}{c_1c_2} \left[ \sum_{i=1}^p \left\{ \log^2\left( (1-c_1)\eta_i \right) - \log^2\left( (1-c_1)\lambda_i \right) \right\} \right.\nonumber \\
&\left.+ 2\sum_{1\leq i,j\leq p} \left\{ {\rm Li}_2\left( 1-\frac{\zeta_i}{\lambda_j}\right) - {\rm Li}_2\left( 1-\frac{\eta_i}{\lambda_j}\right) + {\rm Li}_2\left( 1-\frac{\eta_i}{\eta_j} \right) - {\rm Li}_2\left( 1-\frac{\zeta_i}{\eta_j}\right)\right\} \right] \nonumber \\
				&- \frac{1-c_2}{c_2} \left[ \log^2(1-c_2)-\log^2(1-c_1) + \sum_{i=1}^p \left\{ \log^2\left( \eta_i \right) - \log^2\left( \zeta_i \right) \right\}  \right] \nonumber \\
				&-\frac1p \left[ 2 \sum_{1\leq i,j\leq p} \left\{ {\rm Li}_2\left( 1-\frac{\zeta_i}{\lambda_j}\right) - {\rm Li}_2\left( 1-\frac{\eta_i}{\lambda_j}\right)  \right\} - \sum_{i=1}^p \log^2\left( (1-c_1)\lambda_i \right) \right] \nonumber
\end{align}
which provides an exact, yet rather impractical (the expression involves the evaluation of $O(p^2)$ dilogarithm terms which may be computationally intense for large $p$), final expression for the integral. 

\bigskip

At this point, it is also not easy to fathom why the retrieved expression would remain of order $O(1)$ with respect to $p$. In order to both simplify the expression and retrieve a visually clear $O(1)$ estimate, we next proceed to a large $p$ Taylor expansion of the above result. In particular, using the last item in Lemma~\ref{lem:prop_dilog}, we perform a (second order) Taylor expansion of all terms of the type ${\rm Li}_2(1-X)$ above in the vicinity of $\lambda_i/\lambda_j$.
This results in the following two relations 
\begin{align*}
	&\sum_{i,j}{\rm Li}_2\left(1-\frac{\zeta_i}{\lambda_j}\right)-{\rm Li}_2\left(1-\frac{\eta_i}{\lambda_j}\right)+{\rm Li}_2\left(1-\frac{\eta_i}{\eta_j}\right) - {\rm Li}_2\left(1-\frac{\zeta_i}{\eta_j}\right)\nonumber \\
	&=(\Delta_{\zeta}^{\eta})^{T}M(\Delta_{\lambda}^{\eta})+o_p(1) \\
	&\frac1p\sum_{i,j}{\rm Li}_2\left(1-\frac{\zeta_i}{\lambda_j}\right)-{\rm Li}_2\left(1-\frac{\eta_i}{\lambda_j}\right)\nonumber \\
	&=-\frac1p(\Delta_{\zeta}^{\eta})^{T}N1_p+o_p(1)
\end{align*}
with $\Delta_a^b$, $M$ and $N$ defined in the statement of Corollary~\ref{cor:log2t}.

With these developments, we deduce the final approximation
\begin{align*}
	&\frac1{2\pi\imath}\oint_\Gamma\log^{2}\left( - \frac{\varphi_p(z)}{\psi_p(z)} \right)  \left(\frac{\varphi_p'(z)}{\varphi_p(z)}-\frac{\psi_p'(z)}{\psi_p(z)} \right)\frac{\psi_p(z)}{c_2}dz \\
&=2\frac{c_1+c_2-c_1c_2}{c_1c_2}\left(\left(\Delta_{\zeta}^{\eta}\right)^{T}M\left(\Delta_{\lambda}^{\eta}\right)+\sum_i \frac{\log\left((1-c_1)\lambda_i\right)}{\lambda_i}\left(\eta_i-\lambda_i\right)\right) \\
&-\frac2p\left(\Delta_{\zeta}^{\eta}\right)^{T}N1_p+\frac1p\sum_i\log^{2}\left((1-c_1)\lambda_i\right)\\
&-2\frac{1-c_2}{c_2}\left( \frac12\log^{2}\left(1-c_2\right)-\frac12\log^{2}\left(1-c_1\right)+ \sum_{i}(\eta_i-\zeta_i)\frac{\log(\lambda_i)}{\lambda_i}  \right) + o_p(1).
\end{align*}
For symmetry, it is convenient to finally observe that $\log(1-c_1)\sum_i(\eta_i-\zeta_i)/\lambda_i\sim \log(1-c_1)\sum_i \log(\eta_i/\zeta_i)=-\log^2(1-c_1)$; replacing in the last parenthesis provides the result of Corollary~\ref{cor:log2t} for $c_1>0$.

\bigskip

To determine the limit as $c_1\to 0$, it suffices to remark that in this limit $\eta_i=\lambda_i+\frac{c_1}p\lambda_i+o(c_1)$ (this can be established using the functional relation $\varphi_p(\eta_i)=0$ in the small $c_1$ limit). Thus it suffices to replace in the above expression the vector $\eta-\zeta$ by the vector $\eta-\lambda$, the vector $\frac{c_1+c_2-c_1c_2}{c_1c_2}(\eta-\lambda)$ by the vector $\frac1p\lambda$, and taking $c_1=0$ in all other instances (where the limits for $c_1\to 0$ are well defined).

\section{Integration contour determination}
\label{sec:contour}

This section details the complex integration steps sketched in Appendix~\ref{sec:integral_form}. These details rely heavily on the works of \citep{CHO95} and follow similar ideas as in e.g., \citep{COU10b}.

Our objective is to ensure that the successive changes of variables involved in Appendix~\ref{sec:integral_form} move any complex contour closely encircling the support of $\mu$ onto a valid contour encircling the support of $\nu$; we will in particular be careful that the resulting contour, in addition to encircling the support of $\nu$, does not encircle additional values possibly bringing undesired residues (such as $0$). We will proceed in two steps, first showing that a contour encircling $\mu$ results on a contour encircling $\xi_2$ and a contour encircling $\xi_2$ results on a contour encircling $\nu$.

\bigskip

Let us consider a first contour $\Gamma_{\xi_2}$ closely around the support of $\xi_2$ (in particular not containing $0$). We have to prove that any point $\omega$ of this contour is mapped to a point of a contour $\Gamma_{\nu}$ closely around the support of $\nu$.

The change of variable performed in \eqref{eq:link_mnu_mxi2} reads, for all $\omega\in\CC\setminus{\rm Supp}(\xi_2)$,
\begin{align*}
	z\equiv z(\omega)=\frac{-\omega}{-(1-c_2^\infty )+c_2^\infty \omega m_{\xi_2}(\omega)} =\frac{-1}{m_{\tilde{\xi_2}}(\omega)}
\end{align*}
where we recall that $\tilde\xi_2=c_2^\infty \xi_2+(1-c_2^\infty ){\bm\delta}_0$. Since $\Im[\omega]\Im[m_{\tilde{\xi_2}}(\omega)]>0$ for $\Im[\omega]\neq 0$, we already have that $\Im[z]\Im[\omega]>0$ for all non-real $\omega$. 

\begin{figure}
\centering
	\begin{tabular}{cc}

	\begin{tikzpicture}[font=\footnotesize]
		\renewcommand{\axisdefaulttryminticks}{4} 
		\tikzstyle{every major grid}+=[style=densely dashed]       
		\tikzstyle{every axis legend}+=[cells={anchor=west},fill=white,
		at={(1,1.02)}, anchor=south east, font=\scriptsize ]
		\begin{axis}[
      	axis x line=center,
  		axis y line=center,
  		xlabel={$z$},
		ylabel={$\omega^\circ(z)$},
				width=.6\linewidth,
				xmin=-3,
				ymin=-3,
				xmax=12,
				ymax=12,
				grid=major,
				ymajorgrids=false,
				scaled ticks=true,
				xtick=\empty,
				 ytick=\empty,
			]
			\def\xmi{-3}
			\def\xma{12}
			\def\ymi{-10}
			\def\yma{40}
			\def\xmax{4.19};
			\def\ymax{3.275};
			\def\xmin{6.835};
			\def\ymin{8.332};
			\def\xpt{1.335};
			\def\ypt{1.247};
			\def\xpts{9.355};
			\def\ypts{10.03};
			\def\a{5};
			\def\b{6};
			\def\cof{0.05};
			\def\coff{10};
			\def\ymoy{\ymin/2+\ymax/2};
			\def\amoy{\a/2+\b/2};
			\addplot[gray,densely dashed] plot coordinates{(\xmi,\ymax)(\xma,\ymax)};
			\addplot[gray,densely dashed] plot coordinates{(\xmi,\ymin)(\xma,\ymin)};
						\addplot[gray,densely dashed] plot coordinates{(\a,\ymi)(\a,\yma)};
			\addplot[gray,densely dashed] plot coordinates{(\b,\ymi)(\b,\yma)};
			\addplot[blue,ultra thick] plot coordinates{(\a,0)(\b,0)};
			\addplot[blue,ultra thick] plot coordinates{(0,\ymax)(0,\ymin)};
			\addplot[green,dashed] plot coordinates{(0,\ypt)(\xpt,\ypt)(\xpt,0)};
			\addplot[green,dashed] plot coordinates{(0,\ypts)(\xpts,\ypts)(\xpts,0)};
\addplot[mark=*,only marks] coordinates {(0,\ypt)(\xpt,0)(\xpts,0)(0,\ypts)};
\node[] at (axis cs: -1,\ypt) {$\omega_{-}$};
\node[] at (axis cs: -1,\ypts) {$\omega_{+}$};
\node[] at (axis cs: \xpt,-1) {$z_{-}$};
\node[] at (axis cs: \xpts,-1) {$z_{+}$};
\node[blue,ultra thick] at (axis cs: -1,\ymoy){$\mathcal{S}_{\xi_2}$};
\node[blue,ultra thick] at (axis cs: \amoy,-1){$\mathcal{S}_\nu$};
			\addplot[red,thick,domain=-3:\a, samples=100]{x+\cof*1/(\b-\a)*( -x*(x*ln(\b-x)+\b) +x*(x*ln(\a-x)+\a))};
			\addplot[red,thick,domain=\b:12, samples=100]{x+\cof*1/(\b-\a)*( -x*(x*ln(x-\b)+\b) +x*(x*ln(x-\a)+\a))};
		\end{axis}
	\end{tikzpicture}
	&
	\begin{tikzpicture}[font=\footnotesize]
				\renewcommand{\axisdefaulttryminticks}{4} 
		\tikzstyle{every major grid}+=[style=densely dashed]       
		\tikzstyle{every axis legend}+=[cells={anchor=west},fill=white,
		at={(1,1.02)}, anchor=south east, font=\scriptsize ]
		\begin{axis}[
      	axis x line=center,
  		axis y line=center,
  		xlabel={$z$},
		ylabel={$\omega^\circ(z)$},
				width=.6\linewidth,
				xmin=-32,
				ymin=-30,
				xmax=45,
				ymax=120,
				grid=major,
				ymajorgrids=false,
				scaled ticks=true,
				xtick=\empty,
				 ytick=\empty,
			]
			\def\xmi{-32}
			\def\xma{45}
			\def\ymi{-30}
			\def\yma{180}
			\def\xmax{-10.72};
			\def\ymax{22.98};
			\def\xmin{21.7};
			\def\ymin{88.03};
			\def\xpt{-28.32};
			\def\ypt{13.89};
			\def\xpts{36.54};
			\def\ypts{95.04};
			\def\a{3};
			\def\b{7};
			\def\cof{0.05};
			\def\coff{10};
			\def\ymoy{\ymin/2+\ymax/2};
			\def\amoy{\a/2+\b/2};
			\addplot[gray,densely dashed] plot coordinates{(\xmi,\ymax)(\xma,\ymax)};
			\addplot[gray,densely dashed] plot coordinates{(\xmi,\ymin)(\xma,\ymin)};
						\addplot[gray,densely dashed] plot coordinates{(\a,\ymi)(\a,\yma)};
			\addplot[gray,densely dashed] plot coordinates{(\b,\ymi)(\b,\yma)};
			\addplot[blue,ultra thick] plot coordinates{(\a,0)(\b,0)};
			\addplot[blue,ultra thick] plot coordinates{(0,\ymax)(0,\ymin)};
			\addplot[green,dashed] plot coordinates{(0,\ypt)(\xpt,\ypt)(\xpt,0)};
			\addplot[green,dashed] plot coordinates{(0,\ypts)(\xpts,\ypts)(\xpts,0)};
\addplot[mark=*,only marks] coordinates {(0,\ypt)(\xpt,0)(\xpts,0)(0,\ypts)};
\node[] at (axis cs: -3,\ypt) {$\omega_{-}$};
\node[] at (axis cs: -3,\ypts) {$\omega_{+}$};
\node[] at (axis cs: \xpt,-10) {$z_{-}$};
\node[] at (axis cs: \xpts,-10) {$z_{+}$};
\node[blue,ultra thick] at (axis cs: -3,\ymoy){$\mathcal{S}_{\xi_2}$};
\node[blue,ultra thick] at (axis cs: \amoy,-10){$\mathcal{S}_\nu$};
			\addplot[red,thick,domain=-30:\a, samples=50]{x+\coff*1/(\b-\a)*( -x*(x*ln(\b-x)+\b) +x*(x*ln(\a-x)+\a))};
			\addplot[red,thick,domain=\b:\xma-5, samples=50]{x+\coff*1/(\b-\a)*( -x*(x*ln(x-\b)+\b) +x*(x*ln(x-\a)+\a))};
		\end{axis}
	\end{tikzpicture}
\end{tabular}
\caption{Variable change $z\mapsto \omega^\circ(z)=z+c_2^\infty \int \frac{zd\nu(t)}{z-t}$ for $c_2^\infty <1$ (left) and $c_2^\infty >1$ (right). $\mathcal{S}_\theta$ is the support of the probability measure $\theta$. For $0<\omega_-=\omega(z_-)<\inf \mathcal S_\nu$, the pre-image $z_-$ is necessarily negative for $c_2^\infty >1$.}
	\label{fig:change_variable}
\end{figure}
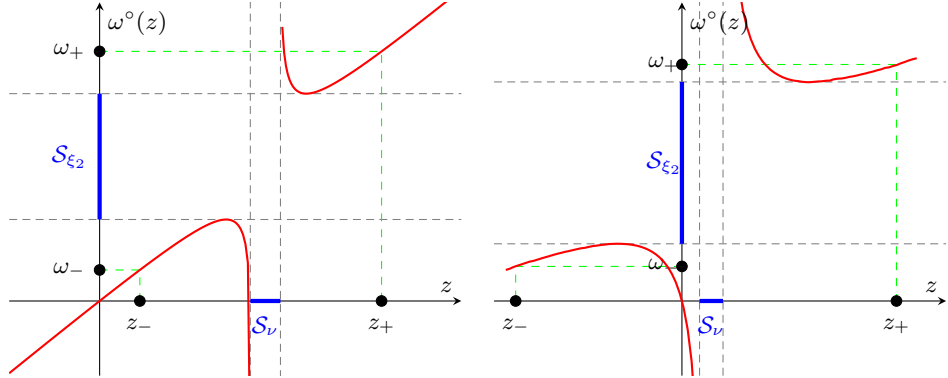

It therefore remains to show that real $\omega$'s (outside the support of $\xi_2$) project onto properly located real $z$'s (i.e., on either side of the support of $\nu$). This conclusion follows from the seminal work \citep{CHO95} on the spectral analysis of sample covariance matrices. The essential idea is to note that, due to $\eqref{eq:baisilverstein2}$, the relation $z(\omega)=-1/m_{\tilde{\xi_2}}(\omega)$ can be inverted as
\begin{align*}
	\omega\equiv \omega(z) &= -\frac{1}{m_{\tilde{\xi_2}}}+c_2^\infty \int\frac{td\nu(t)}{1+tm_{\tilde{\xi_2}}} = z + c_2^\infty \int\frac{td\nu(t)}{1-\frac{t}z}.
\end{align*}
In \citep{CHO95}, it is proved that the image by $\omega(\cdot)$ of $z(\RR\setminus{\rm Supp}(\xi_2))$ coincides with the increasing sections of the function $\omega^\circ:\RR\setminus{\rm Supp(\nu)}\to \RR$, $z\mapsto \omega(z)$. The latter being an explicit function, its functional analysis is simple and allows in particular to properly locate the real pairs $(\omega,z)$. Details of this analysis are provided in \citep{CHO95} as well as in \citep{HAC13}, which shall not be recalled here. The function $\omega^\circ$ is depicted in Figure~\ref{fig:change_variable}; we observe and easily prove that, for $c_2^\infty <1$, any two values $z_-<\inf({\rm Supp}(\nu))\leq \sup({\rm Supp}(\nu))<z_+$ have respectively images $\omega_-$ and $\omega_+$ satisfying $w_-<\inf({\rm Supp}(\xi_2))\leq \sup({\rm Supp}(\xi_2))<w_+$ as desired. This is however not the case for $c_2^\infty >1$ where $\{z_-,z_+\}$  enclose not only ${\rm Supp}(\nu)$ but also $0$  and therefore do not bring a valid contour. This essentially follows from the fact that $(\varphi_p/\psi_p)'(0)$ is positive for $c_2^\infty<1$ and negative for $c_2^\infty>1$.
\bigskip

\begin{figure}
\centering
	\begin{tabular}{cc}

	\begin{tikzpicture}[font=\footnotesize]
		\renewcommand{\axisdefaulttryminticks}{4} 
		\tikzstyle{every major grid}+=[style=densely dashed]       
		\tikzstyle{every axis legend}+=[cells={anchor=west},fill=white,
		at={(1,1.02)}, anchor=south east, font=\scriptsize ]
		\begin{axis}[
      	axis x line=center,
  		axis y line=center,
  		xlabel={$\omega$},
		ylabel={$u^\circ(\omega)$},
		width=.6\linewidth,
				xmin=-3,
				ymin=-3,
				xmax=8,
				ymax=10,
				grid=major,
				ymajorgrids=false,
				scaled ticks=true,
				xtick=\empty,
				 ytick=\empty,
			]
			\def\xmi{-3}
			\def\xma{12}
			\def\ymi{-10}
			\def\yma{16}
			\def\xmax{1.87};
			\def\ymax{1.592};
			\def\xmin{4.255};
			\def\ymin{5.221};
			\def\xpt{0.61};
			\def\ypt{0.601};
			\def\xpts{6.965};
			\def\ypts{7.61};
			\def\a{0.25};
			\def\b{0.5};
			\def\cof{0.05};
			\def\coff{10};
			\def\z{5.78};
			\def\ymoy{\ymin/2+\ymax/2};
			\def\amoy{0.5/\a+0.5/\b};
			\addplot[gray,densely dashed] plot coordinates{(\xmi,\ymax)(\xma,\ymax)};
			\addplot[gray,densely dashed] plot coordinates{(\xmi,\ymin)(\xma,\ymin)};
						\addplot[gray,densely dashed] plot coordinates{(1/\a,\ymi)(1/\a,\yma)};
			\addplot[gray,densely dashed] plot coordinates{(1/\b,\ymi)(1/\b,\yma)};
			\addplot[blue,ultra thick] plot coordinates{(1/\a,0)(1/\b,0)};
			\addplot[blue,ultra thick] plot coordinates{(0,\ymax)(0,\ymin)};
			\addplot[green,dashed] plot coordinates{(0,\ypt)(\xpt,\ypt)(\xpt,0)};
			\addplot[green,dashed] plot coordinates{(0,\ypts)(\xpts,\ypts)(\xpts,0)};
\addplot[mark=*,only marks] coordinates {(0,\ypt)(\xpt,0)(\xpts,0)(0,\ypts)};
\addplot[blue,mark=*,only marks] coordinates {(0,0)};
\node[] at (axis cs: -1,\ypt) {$u_{-}$};
\node[] at (axis cs: -1,\ypts) {$u_{+}$};
\node[] at (axis cs: \xpt,-1) {$\omega_{-}$};
\node[] at (axis cs: \xpts,-1) {$\omega_{+}$};
\node[blue,ultra thick] at (axis cs: -1,\ymoy){$\mathcal{S}_{\mu}$};
\node[blue,ultra thick] at (axis cs: \amoy,-1){$\mathcal{S}_{\xi_2}$};
			\addplot[red,thick,domain=-30:1/\b, samples=1000]{((1/x)+\cof*1/(\b-\a)*( -(\b*x+ln(1-\b*x))/(x^2) +(\a*x+ln(1-\a*x))/(x^2)))^(-1)};
			\addplot[red,thick,domain=1/\a:\z-\z/1000, samples=100]{((1/x)+\cof*1/(\b-\a)*( -(\b*x+ln(-1+\b*x))/(x^2) +(\a*x+ln(-1+\a*x))/(x^2)))^(-1)};
			\addplot[red,thick,domain=\z+\z/1000:\xma, samples=100]{((1/x)+\cof*1/(\b-\a)*( -(\b*x+ln(-1+\b*x))/(x^2) +(\a*x+ln(-1+\a*x))/(x^2)))^(-1)};
		\end{axis}
	\end{tikzpicture}
	&
	\begin{tikzpicture}[font=\footnotesize]
				\renewcommand{\axisdefaulttryminticks}{4} 
		\tikzstyle{every major grid}+=[style=densely dashed]       
		\tikzstyle{every axis legend}+=[cells={anchor=west},fill=white,
		at={(1,1.02)}, anchor=south east, font=\scriptsize ]
		\begin{axis}[
      	axis x line=center,
  		axis y line=center,
  		xlabel={$\omega$},
		ylabel={$u^\circ(\omega)$},
				width=.6\linewidth,
				xmin=-1,
				ymin=-2,
				xmax=4.5,
				ymax=7,
				grid=major,
				ymajorgrids=false,
				scaled ticks=true,
				xtick=\empty,
				 ytick=\empty,
			]
			\def\xmi{-1}
			\def\xma{4}
			\def\ymi{-2}
			\def\yma{7}
			\def\xmax{0.93};
			\def\ymax{0.7502};
			\def\xmin{2.373};
			\def\ymin{3.26};
			\def\xpt{0.42};
			\def\ypt{0.401};
			\def\xpts{3.776};
			\def\ypts{4.412};
			\def\a{1};
			\def\b{2};
			\def\cof{0.7};
			\def\coff{1.1};
			\def\coffe{0.1};
			\def\ymoy{\ymin/2+\ymax/2};
			\def\amoy{0.5*\a+0.5*\b};
			\def\eps{0.1}
			\def\e{-6.62}
			\addplot[gray,densely dashed] plot coordinates{(\xmi,\ymax)(\xma,\ymax)};
			\addplot[gray,densely dashed] plot coordinates{(\xmi,\ymin)(\xma,\ymin)};
						\addplot[gray,densely dashed] plot coordinates{(\a,\ymi)(\a,\yma)};
			\addplot[gray,densely dashed] plot coordinates{(\b,\ymi)(\b,\yma)};
			\addplot[blue,ultra thick] plot coordinates{(\a,0)(\b,0)};
			\addplot[blue,ultra thick] plot coordinates{(0,\ymax)(0,\ymin)};
			\addplot[green,dashed] plot coordinates{(0,\ypt)(\xpt,\ypt)(\xpt,0)};
			\addplot[green,dashed] plot coordinates{(0,\ypts)(\xpts,\ypts)(\xpts,0)};
\addplot[mark=*,only marks] coordinates {(0,\ypt)(\xpt,0)(\xpts,0)(0,\ypts)};
\addplot[blue,mark=*,only marks] coordinates {(0,0)};
\node[] at (axis cs: -0.5,\ypts+0.05) {$u_{+}$};
\node[] at (axis cs: -0.5,\ypt-0.05) {$u_{-}$};
\node[] at (axis cs: \xpt,-1) {$\omega_{-}$};
\node[] at (axis cs: \xpts,-1) {$\omega_{+}$};
\node[blue,ultra thick] at (axis cs: -0.5,\ymoy){$\mathcal{S}_{\mu}$};
\node[blue,ultra thick] at (axis cs: \amoy,-1){$\mathcal{S}_{\xi_2}$};
						\addplot[red,thick,domain=\xmi:\e-\eps, samples=100]{((1+\coffe-(\coffe/\coff))/x+(\coffe/\coff)*ln((\b-x)/(\a-x)))^(-1)};
						\addplot[red,thick,domain=\e+\eps:\a, samples=100]{((1+\coffe-(\coffe/\coff))/x+(\coffe/\coff)*ln((\b-x)/(\a-x)))^(-1)};
			\addplot[red,thick,domain=\b:\xma, samples=100]{((1+\coffe-(\coffe/\coff))/x+(\coffe/\coff)*ln((x-\b)/(x-\a)))^(-1)};
		\end{axis}
	\end{tikzpicture}
\end{tabular}
\caption{Variable change $u^\circ(\omega)=(\frac1\omega+c_1^\infty \int \frac1{t-\omega}d\xi_2(t) )^{-1}$ for $c_2^\infty <1$ (left) and $c_2^\infty >1$ (right). $\mathcal{S}_\theta$ is the support of the probability measure $\theta$.}
	\label{fig:change_variable2}
\end{figure}
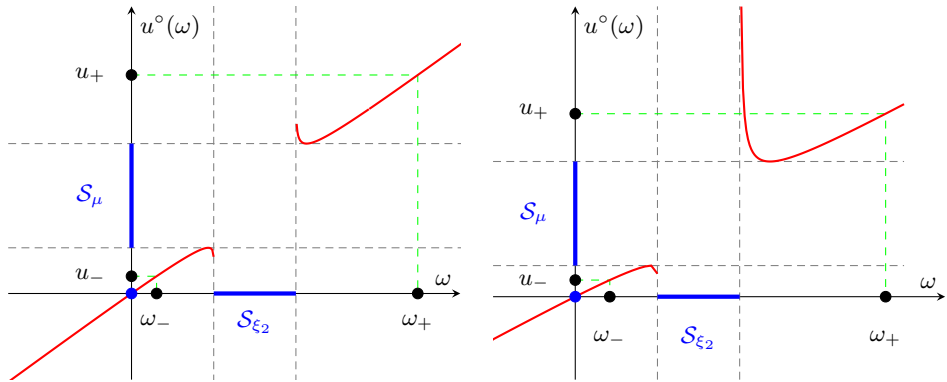

The same reasoning now holds for the second variable change. Indeed, note that here
\begin{align*}
	\omega&=u(1+c_1^\infty um_{\mu}(u))=u\left(1-c_1^\infty -\frac{c_1^\infty }{u}m_{\mu^{-1}}\left(\frac1u\right)\right)=-m_{\tilde{\mu}^{-1}}\left(\frac{1}{u}\right).
\end{align*}
Exploiting \eqref{eq:baisilverstein1} provides, as above, a functional inverse given here by
\begin{align*}
	u &\equiv u(\omega) = \left( \frac{1}{\omega}+c_1^\infty \int\frac{d\xi_2(t)}{t-\omega} \right)^{-1}
\end{align*}
the analysis of which follows the same arguments as above (see display in Figure~\ref{fig:change_variable2} of the extension to $u^\circ(\omega)=u(\omega)$ for all $\omega\in\RR\setminus {\rm Supp}(\xi_2)$). 

\section{The case \lowercase{$c_2^\infty>1$}}
\label{sec:case_c2>1}

Allowing for $c_2^\infty>1$ brings along some key difficulties. First recall from Appendix~\ref{sec:contour} that, for a contour $\Gamma$ surrounding ${\rm Supp}(\mu)$, if $c_2^\infty>1$, the image $(\varphi/\psi)(\Gamma)$ necessarily surrounds ${\rm Supp}(\nu)\cup \{0\}$ (while for $c_2^\infty<1$, $0$ is excluded from the interior of $(\varphi/\psi)(\Gamma)$ if $0$ is not contained within $\Gamma$). This implies that, if $f(z)$ has a singularity at $z=0$, the relation $\int fd\nu(t)=\frac1{2\pi\imath}\oint_{(\varphi/\psi)(\Gamma)}f(z)m_\nu(z)dz$ no longer holds. 

For $f(z)$ a polynomial in $z$, e.g., for $f(z)=z$, this poses no problem, thereby implying the validity of Corollary~\ref{cor:t} for all $c_2^\infty>0$.

For $f(z)=\log(z)^k$ ($k\geq 1$) and similarly for $f(z)=1/z$, to the best of our knowledge, there is no recovering from this technical difficulty. Notably, for $f(z)=\log^k(z)$, one cannot pass a closed path around ${\rm Supp}(\nu)\cup \{0\}$ without crossing a branch cut for the logarithm.\footnote{\label{foo:1}This problem is reminiscent of the simpler-posed, yet still open problem, consisting in evaluating $\frac1p\tr C^{-1}$ based on samples, say, $x_1,\ldots,x_p\sim \mathcal N(0,C)$, for $p>n$. While a consistent so-called G-estimator \cite{GIR87} does exist for all $p<n$, this is not the case when $p>n$.}

\bigskip

The case $f(z)=\log(1+sz)$ is more interesting. As the singularity for $\log(1+sz)$ is located at $z=-1/s<0$, one can pass a contour around ${\rm Supp}(\nu)\cup \{0\}$ with no branch cut issue. However, one must now guarantee that there exists a contour $\Gamma$ surrounding ${\rm Supp}(\mu)$ such that the leftmost real crossing of $(\varphi/\psi)(\Gamma)$ is located within $(-1/s,\inf\{{\rm Supp}(\nu)\})$. This cannot always be guaranteed. Precisely, one must ensure that there exists $u^-<\inf\{{\rm Supp}(\mu)\}$ such that $1+s(\varphi/\psi)(u^-)>0$. In the case $c_2^\infty <1$ where $(\varphi/\psi)(0)=0$, the increasing nature of $\varphi/\psi$ ensures that for all $u^-\in (0,\inf\{{\rm Supp}(\mu)\})$, $1+s(\varphi/\psi)(u^-)>1$ and the condition is fulfilled; however, for $c_2^\infty >1$, it is easily verified that $(\varphi/\psi)(0)<0$. As a consequence, a valid $u^-$ exists if and only if $1+s\lim_{u\uparrow \inf\{{\rm Supp}(\mu)\}}(\varphi/\psi)(u)>0$. 

When this condition is met, a careful calculus reveals that the estimators of Corollary~\ref{cor:log1st} are still valid when $c_2^\infty >1$, with additional absolute values in the logarithm arguments (those were discarded in the proof derivation of Corollary~\ref{cor:log1st} for $c_2^\infty <1$ as the arguments can be safely ensured to be positive). This explains the conclusion drawn in Remark~\ref{rem:log1st_c2>1}. For generic $\nu$, $\inf\{{\rm Supp}(\mu)\}$ is usually not expressible in explicit form (it can still be obtained numerically though by solving the fundamental equations \eqref{eq:link_mmu_mxi2} and \eqref{eq:link_mnu_mxi2}, or estimated by $\min_i\{\lambda_i,\lambda_i>0\}$ in practice). However, for $\nu={\bm \delta}_1$, i.e., for $C_1=C_2$, \cite[Proposition~2.1]{WAN17} provides the exact form of $\inf\{{\rm Supp}(\mu)\}$ and of $m_\mu(z)$; there, a simple yet cumbersome calculus leads to $\lim_{u\uparrow \inf\{{\rm Supp}(\mu)\}}(\varphi/\psi)(u)=(1-\sqrt{c_1^\infty +c_2^\infty -c_1^\infty c_2^\infty })/(1-c_1^\infty )$, which completes the results mentioned in Remark~\ref{rem:log1st_c2>1}.\footnote{Pursuing on the comments of Footnote~\ref{foo:1}, the fact that there exists a maximal value for $s$ allowing for a consistent estimate of $\int \log(1+st)d\nu(t)$ (and thus not of $\int \log(t)d\nu(t)$ which would otherwise be retrieved in a large $s$ approximation) is again reminiscent of the fact that consistent estimators for $\frac1p\tr (C+sI_p)^{-1}$ are achievable when $p>n$, however for not-too small values of $s$, thereby not allowing for taking $s\to 0$ in the estimate. This problem is all the more critical that $p/n$ is large (which we also do observe from Equation~\ref{eq:bound_s} for large values of $c_2$).}

\bibliographystyle{alpha}

\bibliography{/home/romano/Documents/PhD/phd-group/papers/rcouillet/tutorial_RMT/book_final/IEEEabrv.bib,/home/romano/Documents/PhD/phd-group/papers/rcouillet/tutorial_RMT/book_final/IEEEconf.bib,/home/romano/Documents/PhD/phd-group/papers/rcouillet/tutorial_RMT/book_final/tutorial_RMT.bib}

\end{document}